\documentclass[11pt,a4paper,english,reqno]{amsart}
\usepackage[T1]{fontenc}
\usepackage{lmodern,calrsfs}
\usepackage{mathtools,amsthm,amssymb,amscd}
\usepackage{graphicx,xcolor}
\usepackage{array,booktabs,longtable,tabularx,enumitem}
\usepackage{comment,framed}
\usepackage[all]{xy}
\xyoption{2cell}
\UseTwocells

\usepackage{tikz}
\usetikzlibrary{positioning,intersections,cd,matrix,arrows,backgrounds}
\usetikzlibrary{shapes,shapes.geometric,shapes.misc,decorations.pathreplacing}

\IfFileExists{mathx10.tfm}{%
  \DeclareFontFamily{U}{mathx}{}%
  \DeclareFontShape{U}{mathx}{m}{n}{<-> mathx10}{}%
  \DeclareSymbolFont{mathx}{U}{mathx}{m}{n}%
  \DeclareMathAccent{\widehat}{0}{mathx}{"70}%
  \DeclareMathAccent{\widecheck}{0}{mathx}{"71}%
}{\providecommand{\widecheck}{\check}}

\usepackage{microtype}
\usepackage[a4paper,hmargin=24mm,vmargin=20mm]{geometry}
\usepackage[colorlinks=true,allcolors=blue]{hyperref}

\allowdisplaybreaks[2]
\numberwithin{equation}{section}
\setlist[enumerate]{label=\textup{(\arabic*)},leftmargin=2em}

\theoremstyle{plain}
\newtheorem{theorem}{Theorem}[section]
\newtheorem{proposition}[theorem]{Proposition}
\newtheorem{lemma}[theorem]{Lemma}
\newtheorem{corollary}[theorem]{Corollary}

\newtheorem*{thm-nn}{Theorem}
\newtheorem*{prp-nn}{Proposition}
\newtheorem*{lem-nn}{Lemma}
\newtheorem*{cor-nn}{Corollary}
\newtheorem*{clm-nn}{Claim}
\newtheorem*{cnj-nn}{Conjecture}
\newtheorem*{prb-nn}{Problem}
\theoremstyle{definition}
\newtheorem{definition}[theorem]{Definition}
\newtheorem{example}[theorem]{Example}
\newtheorem{remark}[theorem]{Remark}

\newtheorem*{dfn-nn}{Definition}
\newtheorem*{convention}{Conventions}

\def\dsm#1,#2..#3{\bigoplus_{{#1}={#2}}^{#3}}
\def\sm#1,#2..#3{\sum_{{#1}={#2}}^{#3}}

\def\id{1\kern-.25em{\text{{\rm l}}}} %1\!\!{\text{{\rm l}}}   is good for 12pt

\def\isoto{\ \raise.3ex\hbox{$^{\sim}$}\kern-.8em\hbox{$\to$}\ } 
\def\bg{%
\family{cmr}\size{20}{12pt}\selectfont}

\def\bigzerou{%
\smash{\lower1.7ex\hbox{\bg 0}}}

\def\repr[#1;#2;#3;#4;#5]{
\left(
\begin{matrix}#1\\#2\end{matrix}
#3
\begin{matrix}#4\\#5\end{matrix}
\right)}
\def\Bi(#1,#2,#3){{}_{#1}{#2}_{#3}}
\def\Biup(#1,#2,#3){{}^{#1}{#2}^{#3}}
\def\Smor(#1,#2,#3){{}^{(#1)}{#2}^{(#3)}}

\hypersetup{hypertexnames=false,
 pdfauthor={Shengyong Pan},
 pdftitle={Hochschild operations and other invariants of Morita-type stable equivalences},
 pdfsubject={Homological invariants and relative stable equivalences}}

\title[Invariants of stable equivalences of Morita type]
{Hochschild operations and other invariants of stable equivalences of Morita type}
\author{Shengyong Pan}
\address{School of Mathematics and Statistics, Beijing Jiaotong University, Beijing 100044, China}
\email{shypan@bjtu.edu.cn}
\date{}

\begin{document}

\begin{abstract}
We study stable equivalences defined by bimodule tensor identities
over a commutative ring $R$. For the $R$-split exact structure,
these identities preserve positive Hochschild homology and
positive Hochschild cohomology with its cup product, Gerstenhaber
bracket and, in prime characteristic, restricted power. A
triangular comparison also identifies the relative singular
Hochschild cochain complexes as $B_\infty$-algebras.
For reduced pairs of finite-dimensional Frobenius algebras
over an algebraically closed field, with diagonalizable
Nakayama automorphisms, we construct compatible Frobenius
functionals and compare the Tate--Hochschild BV operators in
all degrees.

For flat algebras over $R$ we give a Hochschild-cohomological
criterion for comparing marked associative deformations.
It applies to biadjoint Morita-type pairs and, by reduction,
to all stable equivalences of Morita type between
finite-dimensional algebras over any field. We compare the
Bezrukavnikov--Ginzburg maps in characteristic zero by using
a common triangular deformation or by lifting a two-sided
tilting complex. We also compare marked module deformations
over local Artin rings and their augmented universal rings.

Further applications concern connected outer automorphism
group schemes over perfect fields, the kernels of the
quotient maps on ordinary algebraic $K$-theory, and the
nonconnective $K$-theory of enhanced relative singularity
categories. We construct biadjoint representatives over
Henselian Noetherian local rings and relative stable tensor
functors for almost $\nu$-stable derived equivalences.
We also answer a question of Abe and Hoshino about preservation
of the Gorenstein condition in a fixed codimension.
\end{abstract}

\maketitle
\setcounter{tocdepth}{1}{\small\tableofcontents}
\section{Introduction}
\label{sec:introduction}

Let $R$ be a commutative ring and let $A$ and $B$ be $R$-algebras.
A pair of bimodules defining a stable equivalence of Morita type over a
field gives mutually inverse tensor functors on the stable module
categories. Over $R$, the same identities give an equivalence of the
relative stable categories for the exact structure of $R$-split
sequences. In this paper, we use these identities to compare Hochschild
operations, deformation functors, outer automorphism groups, algebraic
$K$-theory, and homological dimensions.

\begin{convention}
Unless otherwise stated, $R$ is a commutative Noetherian ring.
Every $R$-algebra $A$ is finite as an $R$-module, and its
structure map has image in $Z(A)$. All bimodules are $R$-central,
that is, their two $R$-actions agree.

We write $A\text{-}\mathrm{Mod}$ for the category of all left
$A$-modules, $A\text{-}\mathrm{mod}$ for its full subcategory of
finitely generated modules, and $A\text{-}\mathrm{proj}$ for
the category of finitely generated projective left $A$-modules.
We set $\mathcal D(A)=D(A\text{-}\mathrm{Mod})$ and write
$K^b(A\text{-}\mathrm{proj})$ for the bounded homotopy category
of $A\text{-}\mathrm{proj}$. Ordinary injective resolutions
are taken in $A\text{-}\mathrm{Mod}$.

We identify right $A$-modules with left $A^{\mathrm{op}}$-modules.
For an $A$-$B$-bimodule $X$, the notation $X^{\mathrm{op}}$
denotes the $B^{\mathrm{op}}$-$A^{\mathrm{op}}$-bimodule with
action $b^{\mathrm{op}}xa^{\mathrm{op}}=axb$.
\end{convention}

Let $A$ and $B$ be $R$-algebras. We use the following Morita-type
bimodule identities.
\begin{equation}
\begin{gathered}
 {}_AM_B\text{ and }{}_BN_A\text{ are finitely generated projective
 on each side},\\
 M\otimes_BN\simeq A\oplus P,\qquad
 N\otimes_AM\simeq B\oplus Q,\\
 P\in (A\otimes_RA^{\mathrm{op}})\text{-}\mathrm{proj},\qquad
 Q\in (B\otimes_RB^{\mathrm{op}})\text{-}\mathrm{proj}.
\end{gathered}
\tag{MT}\label{eq:MT}
\end{equation}
 We refer to \eqref{eq:MT} as the
direct-summand identities. We also consider the stabilized identities,
which assert that the two tensor products are isomorphic to the regular
bimodules in the relative stable bimodule categories.

For each left $A$-module $X$, the module $P\otimes_AX$ is a direct
summand of a finite direct sum of copies of $A\otimes_RX$. Thus
$P\otimes_AX$ is relatively projective, and similarly for $Q$.
It follows that the two tensor functors induce mutually inverse
equivalences of the relative stable categories
\cite{RouquierZimmermann2003,Linckelmann2025}. If $R$ is a field, these
are the usual stable module categories. If $A$ is projective over $R$,
relative and ordinary projectivity agree on $A$-modules that are
projective over $R$.

\subsection{Hochschild invariants}

For finite-dimensional algebras over a field, Liu and Xi proved
invariance of positive Hochschild homology
\cite[Theorem~4.4]{LiuXi2005}, and Pan and Zhou proved invariance of the
stable Hochschild cohomology ring \cite[Theorem~1.2]{PanZhou2010}.
In Section~\ref{sec:hochschild-without-adjunction}, we give the
corresponding ordinary and relative comparisons over a commutative
base. In the proof of the ring isomorphism we show that the
off-diagonal components vanish for positive-degree classes in the
image of the tensor functor. Proposition~\ref{prop:pan-zhou-compression}
uses this vanishing to correct the compression argument in
\cite[Theorem~1.2]{PanZhou2010}.

In Theorem~\ref{cor:gerstenhaber-stable-identities}, we show that
under the stabilized identities the positive relative Hochschild
cohomologies of $A$ and $B$ are isomorphic, with their cup products,
Gerstenhaber brackets, and restricted powers in prime characteristic.
The proof uses the fact that the two characteristic maps into
$\operatorname{Ext}^n_{A\otimes_R B^{\mathrm{op}},\mathrm{rel}}(M,M)$
have the same image. The restrictions from the triangular algebra
with off-diagonal bimodule $M$ respect the brace operations and
induce the required comparison.

Briggs and Rubio y Degrassi proved invariance of the restricted Lie
structure for self-injective algebras and in sufficiently high
degrees for Iwanaga--Gorenstein algebras \cite{BriggsRubio2022}.
Wang proved the corresponding Gerstenhaber invariance for singular
Hochschild cohomology \cite[Theorem~6.2]{Wang2021Gerstenhaber}.
Our comparison applies in every positive ordinary degree over a
field and in every positive relative degree over a commutative ring.
An example of rank seven in
Subsection~\ref{sec:stable-center-bracket-counterexample} shows that
the projective center can fail to be a Gerstenhaber Lie ideal.

For algebras over a field, Chen, Li and Wang proved that singular
equivalences of Morita type with level preserve the $B_\infty$
homotopy type of the singular Hochschild complex
\cite[Theorem~9.6 and Proposition~9.7]{ChenLiWangBInfinity}.
In Theorem~\ref{thm:relative-singular-braces}, we prove the relative
version over a commutative ring using the $R$-split exact structure.
Corollary~\ref{cor:singular-braces-bv} identifies the cohomology map
of the singular comparison with the Tate transfer in the Frobenius
case considered here.

\subsection{The BV operator and deformations}

K\"onig, Liu and Zhou compared BV structures for symmetric algebras
under the hypotheses of \cite[Theorem~10.7]{KoenigLiuZhou2012}.
For a Frobenius algebra with diagonalizable Nakayama automorphism,
Itagaki, Sanada and Usui constructed a BV structure on
Tate--Hochschild cohomology \cite{ISUTateBV}.
In Theorem~\ref{thm:tate-bv-compatible}, we prove that the comparison
associated with a reduced pair preserves this BV operator in every
degree if the Frobenius functionals are compatible. Under the
hypotheses of Section~\ref{sec:tate-bv},
Corollary~\ref{cor:tate-bv-normalization} shows that one may prescribe
one functional and choose the other by multiplying an initial
functional by a central unit. This preserves the Nakayama
automorphism. The proof in negative degrees uses the twisted
Hochschild chain maps. At degree $-1$ we also identify the two
kernels of the norm maps in the complete complexes.

In Theorem~\ref{thm:stable-associative-deformations}, we compare
marked associative deformations of flat algebras over a commutative
ring by taking the two corners of a triangular algebra. The
criterion uses Hochschild cohomology with square-zero coefficients
in degrees one, two and three. Biadjoint Morita-type pairs satisfy
this criterion. By reduction to blocks, the comparison applies to
all stable equivalences of Morita type between finite-dimensional
algebras over any field. For the relation with Liao--Seol admissible
triples and the corresponding comparison in the singular category
of bimodules, see Subsection~\ref{subsec:relative-triangular-comparison}.
In characteristic zero, the corner restrictions also carry the
Bezrukavnikov--Ginzburg map of a common triangular deformation
to the maps for its two corners
(Proposition~\ref{prop:bg-triangular-deformations}).
For a standard derived equivalence between flat algebras over a
commutative ring, we lift its two-sided tilting complex to compare
marked deformations over nilpotent coefficient extensions.
The comparison extends to complete coefficient rings by lifting
compatible systems of complexes.
The lifted kernels give the comparison of the two
Bezrukavnikov--Ginzburg maps in the homotopy category of dg algebras
(Theorem~\ref{thm:bg-standard-derived}).

We also consider module deformations over a commutative local
Artin ring. Theorem~\ref{thm:artin-module-deformation-morita}
extends the field-level universal deformation ring comparison to
lifts marked over the Artin base, using the Gorenstein projectivity
of the residue module. For a module free over its annihilator
quotient, the corresponding functor is pro-representable over
that quotient under the usual stable endomorphism hypothesis.
Example~\ref{ex:artin-module-deformation-syzygy} computes the
comparison for modules of ranks one and two over a mixed
characteristic base.
Example~\ref{ex:artin-module-deformation-nonrepresentable} shows
that pro-representability can fail upon returning to the original
test category, even for test rings free over the original base.

\subsection{Outer automorphism groups and algebraic $K$-theory}

Rouquier proved invariance of the connected outer automorphism
group for self-injective algebras
\cite[Theorem~4.15]{Rouquier2011Automorphisms} and treated the
group-scheme structure in \cite[Theorem~4.3]{RouquierGroupesPreprint}.
For a reduced pair the error bimodules are projective and injective
over their enveloping algebras \cite[Lemma~4.1(3)]{CPX2012}. We use
this property to split the error summands in algebraic families
and over local Artin rings. Thus, over a perfect field, a stable
equivalence of Morita type induces an isomorphism between the
identity components of the outer automorphism group schemes
(Corollary~\ref{thm:outer-group-scheme}). The tangent isomorphism in
Theorem~\ref{thm:formal-outer} agrees with the comparison on
$\mathrm{HH}^1$. Consequently, it preserves integrability at
each order. For the earlier comparisons of integrable derivations
and congruence subgroups, we refer to
\cite[Theorem~5.1 and Lemma~5.2]{LinckelmannIntegrable} and
\cite[Theorem~2]{BriggsRubioIntegrable}.

Let $A$ and $B$ be finite-dimensional algebras over a field.
We prove that a stable equivalence of Morita type identifies the
kernels of the maps $K_i(A)\to K_i(A/\operatorname{rad}A)$ and
$K_i(B)\to K_i(B/\operatorname{rad}B)$
(Theorem~\ref{thm:radical-kernel-invariance}). If both quotient maps
have ring sections, the homotopy fibres are equivalent
(Theorem~\ref{thm:radical-spectra}). In particular, for split
algebras the ordinary $K$-theory spectra are equivalent if and only
if the algebras have the same number of isomorphism classes of
simple modules. Example~\ref{ex:k-tensor-syzygy} shows that the map
on ordinary $K$-theory induced by a defining bimodule can fail to
be an isomorphism.

Over a commutative ring, we show that the stabilized identities
induce an equivalence of the enhanced relative singularity
categories for the $R$-split exact structure, and hence an
equivalence on their nonconnective $K$-theory
(Theorem~\ref{thm:relative-singular-k-invariance}). Over a field,
the enhancement comparison also follows from the tensor-induced
dg lift of Chen, Liu and Wang
\cite[Proposition~3.1 and Lemma~3.2]{ChenLiuWang2023}.

K\"onig, Liu and Zhou proved invariance of stable cyclic
homology over a field \cite[Theorem~9.9]{KoenigLiuZhou2012}.
We give a relative construction in
Proposition~\ref{prop:relative-stable-cyclic}. We also apply the
radical filtration to dg additive correspondences and obtain
the comparison of the radical mixed-complex fibres
(Corollary~\ref{cor:radical-mixed-complex}).
Theorem~\ref{thm:relative-nilpotent-additive} gives a version for
nilpotent ideals with split filtrations over the base.

\subsection{Adjoints, homological dimensions, and tensor realization}

Over Henselian Noetherian local rings, we construct biadjoint
representatives for pairs satisfying the Morita-type identities,
allowing relatively projective error bimodules and algebras that
are finite over the base. For related constructions over fields
and Artin rings, we refer to \cite[Lemma~4.1]{CPX2012} and
\cite[Lemma~3.11]{XZ2026}. The direct-summand identities and
biadjointness give the comparisons of individual homological
dimensions by the arguments of
\cite[Lemma~4.2 and Proposition~4.5]{bel-LiuXi2007} and
\cite[Theorem~3.2]{ChenRen2022}. For the original pair, the retraction
argument of \cite{bel-Tang2026} gives equality of the left
self-injective dimensions and of the right self-injective
dimensions. The same retractions give equality of individual
projective dimensions. In
Section~\ref{sec:individual-gorenstein-dimensions}, we give a
criterion for the original pair to preserve individual Gorenstein
projective dimensions.

For an almost $\nu$-stable derived equivalence over an Artin ring,
Hu and Xi constructed an ordinary stable equivalence and, over
a field, proved that it is induced by a stable equivalence of
Morita type \cite[Theorems~3.7 and~5.3]{HuXi2010}.
Theorem~\ref{anu-thm-raw-semt} extends the bimodule identities
to Artin algebras that are projective over the base and gives a
relative stable tensor realization.
Theorem~\ref{thm:gpd-almost-nu} records the equality of individual
Gorenstein projective dimensions under the ordinary stable functor
as a consequence of known results. For ordinary stable categories,
we also give a local criterion under which a central tensor functor
inducing a stable equivalence is a Morita equivalence
(Corollary~\ref{cor:rigidity-global}).

Abe and Hoshino ask, after \cite[Example~4.6]{AbeHoshino2009},
whether concentration of the base-ring Ext groups suffices to
preserve the Gorenstein condition of a fixed codimension under
a derived equivalence. We answer this question affirmatively
in Theorem~\ref{thm:abe-hoshino-answer}. The proof reduces the
tilting complex along regular parameters and uses the
self-injectivity comparison over an Artinian local base.

The paper is organized as follows. In Section~\ref{sec:relative-and-abe},
we recall relative stable categories and fix the bimodule notation.
In Sections~\ref{sec:hochschild-without-adjunction}
and~\ref{sec:relative-gerstenhaber}, we compare Hochschild homology,
cohomology, and operations, explain the correction to the
compression argument, and compare the relative singular
Hochschild cochains. In Section~\ref{sec:tate-bv}, we construct
compatible Frobenius functionals and identify the singular
comparison with the Tate--Hochschild BV transfer. In Section~\ref{sec:associative-deformations}, we
study marked associative deformations and the
Bezrukavnikov--Ginzburg maps. Subsection~\ref{subsec:standard-derived-deformations}
also compares graded deformations under Koszul duality.
In Subsection~\ref{subsec:artin-module-deformations}, we study
module deformations and pro-representability over local Artin rings.
In Sections~\ref{sec:outer-groups} and~\ref{sec:k-theory}, we apply
the bimodule identities to outer automorphisms, algebraic
$K$-theory, stable cyclic homology, and additive radical fibres.

In Section~\ref{sec:adjoint-normalization}, we construct biadjoint
representatives and apply the tensor comparisons to homological
dimensions in Sections~\ref{sec:bel-correction}
and~\ref{sec:individual-gorenstein-dimensions}.
Section~\ref{sec:ordinary-stable-rigidity} is devoted to ordinary
stable tensor functors. In Sections~\ref{sec:relative-frobenius-derived}
and~\ref{anu-sec}, we discuss relative stable equivalences arising
from derived equivalences, including the construction of Rouquier
and Zimmermann. In Subsection~\ref{subsec:abe-hoshino-answer},
we prove the codimension statement for endomorphism algebras
of tilting complexes. Finally, in Appendices~\ref{sec:abe-boundary}
and~\ref{sec:separable-boundary}, we give examples concerning
projectivity over the base and separable tensor functors.

\section*{Acknowledgments}
This work was supported by the Beijing Natural Science Foundation (Grant Nos.~1262017, 1252011).

\section{Relative stable categories}
\label{sec:relative-and-abe}

In this section, we recall the relative stable categories used
in the construction of Rouquier and Zimmermann
(Section~\ref{sec:relative-frobenius-derived}). Throughout, $R$ is
a commutative Noetherian ring, all algebras are finite over $R$, all
modules are finitely generated, and all bimodules are $R$-central.

For an $R$-algebra $A$, write
$A^e=A\otimes_RA^{\mathrm{op}}$ and put
\[
 \mathcal P_R(A)
   =\operatorname{add}\{A\otimes_RV:V\in R\text{-}\mathrm{mod}\},
 \qquad
 \mathcal I_R(A)
   =\operatorname{add}\{\operatorname{Hom}_R(A,V):
                                    V\in R\text{-}\mathrm{mod}\}.
\]
These are, respectively, the projective and injective objects in
$A\text{-}\mathrm{mod}$ with respect to the $R$-split exact structure.
The canonical epimorphism
$A\otimes_RX\to X$ and monomorphism
$X\to\operatorname{Hom}_R(A,X)$ are $R$-split.  Hence this exact
category has enough projectives and enough injectives.  We define the additive quotient
\begin{equation}
 A\text{-}\underline{\mathrm{mod}}_R
    =(A\text{-}\mathrm{mod})/[\mathcal P_R(A)],
 \label{eq:relative-stable-category}
\end{equation}
where $[\mathcal P_R(A)]$ denotes the ideal of maps factoring
through objects of $\mathcal P_R(A)$.  We call
\eqref{eq:relative-stable-category} the $R$-relative stable category;
see \cite[Section~2]{Linckelmann2025}.  Its definition does not require
$\mathcal P_R(A)=\mathcal I_R(A)$.

Throughout this paper, ``relative'' refers to the displayed
$R$-split exact structure and the class $\mathcal P_R(A)$.
If $R$ is a field,
then $\mathcal P_R(A)=A\text{-}\mathrm{proj}$, so
\eqref{eq:relative-stable-category} is the ordinary stable module
category. The statements below use explicit bimodule identities
and the specified class of relative projective modules.

\subsection{Bimodule identities and relative stable equivalences}

\begin{proposition}\label{prop:relative-morita-bimodules}
Let ${}_AM_B$ and ${}_BN_A$ be projective on both sides.  Suppose that
there are bimodule isomorphisms
\begin{equation}
 M\otimes_BN\simeq A\oplus P,
 \qquad N\otimes_AM\simeq B\oplus Q,
 \qquad P\in A^e\text{-}\mathrm{proj},\quad
        Q\in B^e\text{-}\mathrm{proj}.
 \label{eq:relative-raw-mt}
\end{equation}
Then $N\otimes_A-$ and $M\otimes_B-$ preserve $R$-split short
exact sequences and relative projective modules.  They induce
inverse additive equivalences
\[
 A\text{-}\underline{\mathrm{mod}}_R
       \simeq B\text{-}\underline{\mathrm{mod}}_R.
\]
The same conclusions hold if \eqref{eq:relative-raw-mt} is replaced
by the stabilized identities
\begin{equation}
 \begin{split}
 (M\otimes_BN)\oplus P_1&\simeq A\oplus P_2,\\
 (N\otimes_AM)\oplus Q_1&\simeq B\oplus Q_2,
 \end{split}
 \qquad
 P_i\in A^e\text{-}\mathrm{proj},\quad
 Q_i\in B^e\text{-}\mathrm{proj}.
 \label{eq:relative-stabilized-mt}
\end{equation}
\end{proposition}

\begin{proof}
Since $N_A$ is finitely generated projective, the composite of
$N\otimes_A-$ with restriction to $R$-modules is a direct summand
of a finite direct sum of restriction functors.  An $R$-linear
splitting therefore induces a splitting after tensoring: compose
the splitting on the finite direct sum with the inclusion and
projection of this direct summand.  Thus $N\otimes_A-$ preserves
$R$-split short exact sequences.  The same holds for $M\otimes_B-$.

For $V\in R\text{-}\mathrm{mod}$,
\[
 N\otimes_A(A\otimes_RV)\simeq N\otimes_RV
\]
is relatively $R$-projective as a $B$-module, since ${}_BN$ is a
direct summand of a finite direct sum of copies of $B$.  Thus the
functors preserve relative projectives and descend to the
quotients.

For every $X\in A\text{-}\mathrm{mod}$, the projectivity of $P$
over $A^e$ implies
\[
 P\otimes_AX\in
   \operatorname{add}(A\otimes_RX)\subseteq\mathcal P_R(A).
\]
Indeed, tensor a split inclusion of $P$ into a finite free
$A^e$-module with $X$ over $A$.  Hence the natural decomposition
\[
 M\otimes_BN\otimes_AX\simeq X\oplus(P\otimes_AX)
\]
induces a natural isomorphism with the identity on the relative
stable category.  The other composite is treated similarly.
Under \eqref{eq:relative-stabilized-mt}, tensoring again gives
these natural isomorphisms, since the additional relatively
projective summands vanish in the quotient.
\end{proof}

\section{Hochschild homology and cohomology without adjunction}
\label{sec:hochschild-without-adjunction}

Let $R$ be a commutative ring and let $A$ and $B$ be unital
central $R$-algebras. We work in the categories of all modules and
bimodules over the indicated rings. Put
$A^e=A\otimes_R A^{\mathrm{op}}$ and $B^e=B\otimes_R B^{\mathrm{op}}$.
We assume that $M$ and $N$ are finitely generated projective on
each side and satisfy the following conditions.

\begin{equation}
\begin{gathered}
 {}_AM_B\text{ and }{}_BN_A\text{ are projective on both sides},\\
 M\otimes_B N\simeq A\oplus P,
 \qquad N\otimes_A M\simeq B\oplus Q,\\
 P\in A^e\text{-}\mathrm{proj},
 \qquad Q\in B^e\text{-}\mathrm{proj}.
\end{gathered}
\tag{HH-MT}\label{eq:hh-mt}
\end{equation}
Here ``projective bimodule'' means projective over the indicated
enveloping algebra.  The results below require neither adjointness nor
preservation of injectives for the tensor functors.  The algebras need
not be projective over $R$, and the tensor functors need not induce
inverse equivalences of the ordinary stable module categories.

\subsection{Two conventions over a non-field base}

We distinguish the ordinary derived functors over the enveloping
algebra from Hochschild (co)homology computed by the $R$-relative bar
construction.  Write
\[
 \mathrm{HH}^{n}_{\mathrm{abs}}(A/R)
   =\operatorname{Ext}_{A^e}^{n}(A,A),
 \qquad
 \mathrm{HH}^{\mathrm{abs}}_{n}(A/R;X)
   =\operatorname{Tor}^{A^e}_{n}(A,X),
\]
and omit $X$ when $X=A$.  In the Tor expression, the right $A^e$-action
on $A$ is $x(a\otimes b^{\mathrm{op}})=bxa$.
The subscript ``abs'' refers to ordinary Ext and Tor over $A^e$; it does
not mean Hochschild homology with the base changed to $\mathbb Z$.

For an $R$-algebra $C$, the $R$-split exact structure on
$C$-modules consists of the short exact sequences which split after
restriction to $R$.  Its projective objects are the direct summands of
modules $C\otimes_R V$.  Denote the associated relative derived functors
by $\operatorname{Ext}_{C,\mathrm{rel}}^n$ and
$\operatorname{Tor}^{C,\mathrm{rel}}_n$, and put
\[
 \mathrm{HH}^{n}_{\mathrm{rel}}(A/R)
   =\operatorname{Ext}_{A^e,\mathrm{rel}}^{n}(A,A),
 \qquad
 \mathrm{HH}^{\mathrm{rel}}_{n}(A/R;X)
   =\operatorname{Tor}^{A^e,\mathrm{rel}}_{n}(A,X).
\]
The augmented bar complex
\[
 \cdots\longrightarrow A\otimes_R A\otimes_R A\otimes_R A
 \longrightarrow A\otimes_R A\otimes_R A
 \longrightarrow A\otimes_R A\longrightarrow A
\]
is an $R$-split resolution by relatively projective $A^e$-modules.
Consequently, these relative groups are the groups computed
by the usual $R$-relative Hochschild chain and cochain complexes.
Their Yoneda product agrees with the cup product.

If $A$ is projective over $R$, each term of this bar resolution is also
projective over $A^e$.  In this case the two theories agree, including
their cohomology products.  They may differ when $A$ is not projective
over $R$.  Neither notation denotes Shukla cohomology or Hochschild
cohomology defined by a derived enveloping algebra.

In statements valid for both conventions, $\varepsilon$ denotes either
$\mathrm{abs}$ or $\mathrm{rel}$.  The notation
$\operatorname{Ext}_{C,\mathrm{abs}}^n$ and
$\operatorname{Tor}^{C,\mathrm{abs}}_n$ means ordinary Ext and Tor.

\subsection{Tensor functors and relative exactness}

\begin{lemma}\label{lem:hh-tensor-properties}
Under \eqref{eq:hh-mt}, the bimodules $P$ and $Q$ are projective on
both sides.  The functors
\[
 \begin{aligned}
 F&=N\otimes_A-\otimes_A M:A^e\text{-}\mathrm{Mod}
       \longrightarrow B^e\text{-}\mathrm{Mod},\\
 G&=M\otimes_B-\otimes_B N:B^e\text{-}\mathrm{Mod}
       \longrightarrow A^e\text{-}\mathrm{Mod}
 \end{aligned}
\]
are exact, preserve ordinary projective bimodules, preserve $R$-split
short exact sequences, and preserve relatively $R$-projective
bimodules.  Their extensions to either the ordinary bounded derived
categories or the bounded derived categories of the $R$-split exact
structures are faithful.
\end{lemma}

\begin{proof}
Since ${}_BN$ is projective, ${}_A(M\otimes_BN)$ is a direct summand of
a finite direct sum of copies of ${}_AM$.  The analogous argument on
the right proves that $M\otimes_BN$, and hence $P$, is projective on
both sides.  The argument for $Q$ is identical.  One-sided
projectivity of $M$ and $N$ proves exactness of $F$ and $G$.

There is a bimodule isomorphism
\[
 F(A^e)\simeq {}_BN\otimes_R M_B.
\]
Since ${}_BN$ and $M_B$ are projective, the right-hand side is a direct
summand of a finite direct sum of copies of $B\otimes_R B$, and is
therefore $B^e$-projective.  This proves that $F$ preserves ordinary
projectives; the proof for $G$ is the same.  Only split inclusions are
tensored, so flatness over $R$ is unnecessary.

If $T_C$ is a finitely generated projective right
$C$-module, the underlying $R$-module functor $T\otimes_C-$ is a
retract of a finite direct sum of the restriction functor
$C\text{-}\mathrm{Mod}\to R\text{-}\mathrm{Mod}$.  It therefore
takes an $R$-split short exact sequence to an $R$-split short exact
sequence.  A splitting is obtained by composing a splitting on the
finite direct sum with the inclusion and projection defining the
retract.  The same argument applies to tensoring on the other side.
Applying it twice proves that $F$ and $G$ preserve $R$-split exact
sequences.
Moreover,
\[
 F(A^e\otimes_R V)\simeq N\otimes_R V\otimes_R M
\]
is a direct summand of a finite direct sum of copies of
$B^e\otimes_R V$.  Thus $F$ and $G$ also preserve relative
projectives.

Put $U=M\otimes_BN\simeq A\oplus P$.  The fixed bimodule
decomposition gives an isomorphism of exact functors
\begin{equation}
 GF\simeq\operatorname{Id}\oplus K_A,
 \qquad
 K_A(X)=(P\otimes_AX)\oplus(X\otimes_AP)
           \oplus(P\otimes_AX\otimes_AP).
 \label{eq:hh-gf-splitting}
\end{equation}
Similarly, $FG\simeq\operatorname{Id}\oplus K_B$.
These isomorphisms descend to both derived categories.  If $F(f)=0$,
then $GF(f)=0$, whose component on the identity summand is $f$.
Thus $f=0$.  The same argument proves faithfulness of $G$.
The induced derived functors need not be full.
\end{proof}

\subsection{Positive Hochschild homology}

We first establish a cyclic Tor comparison by means of projective
resolutions.

\begin{lemma}\label{lem:hh-cyclic-tor}
Let ${}_AM_B$ and ${}_BN_A$ be bimodules, with $N$ projective on both
sides, and set $C=A\otimes_R B^{\mathrm{op}}$.  Regard $N$ as a right
$C$-module by $n(a\otimes b^{\mathrm{op}})=bna$.  For either
$\varepsilon$ and every $n\geq0$ there are isomorphisms
\begin{equation}
 \mathrm{HH}^{\varepsilon}_n(A/R;M\otimes_BN)
 \simeq \operatorname{Tor}^{C,\varepsilon}_n(N,M)
 \simeq \mathrm{HH}^{\varepsilon}_n(B/R;N\otimes_AM).
 \label{eq:hh-cyclic-tor}
\end{equation}
\end{lemma}

\begin{proof}
For ordinary Tor, choose a projective
$C$-resolution $L_\bullet\to M$.  Since ${}_BN$ and $N_A$ are
projective, both
\[
 L_\bullet\otimes_BN\longrightarrow M\otimes_BN,
 \qquad
 N\otimes_AL_\bullet\longrightarrow N\otimes_AM
\]
are exact augmented complexes.  If $L_i=C$, their terms become,
respectively, $A\otimes_RN$ and $N\otimes_RB$.  The former is
$A^e$-projective because $N_A$ is projective, and the latter is
$B^e$-projective because ${}_BN$ is projective.  Taking direct
summands proves the same assertion for every projective $L_i$.
Thus these are ordinary projective resolutions over the respective
enveloping algebras.

For every $A$--$B$-bimodule $L$ there are natural trace isomorphisms
\[
 A\otimes_{A^e}(L\otimes_BN)
       \longrightarrow N\otimes_C L,
       \qquad a\otimes(l\otimes n)\longmapsto na\otimes l,
\]
and
\[
 B\otimes_{B^e}(N\otimes_AL)
       \longrightarrow N\otimes_C L,
       \qquad b\otimes(n\otimes l)\longmapsto n\otimes lb.
\]
Their inverses send $n\otimes l$ to $1\otimes(l\otimes n)$ and
$1\otimes(n\otimes l)$, respectively.  Applied degreewise to
$L_\bullet$, these are isomorphisms of chain complexes.  Taking
homology gives \eqref{eq:hh-cyclic-tor} in the ordinary convention.

For relative Tor, choose an $R$-split resolution of $M$
by relatively projective $C$-modules.  Tensoring with $N$ over either
$A$ or $B$ preserves its $R$-split exactness by the retract argument
in Lemma~\ref{lem:hh-tensor-properties}.  A term $C\otimes_RV$
becomes $A\otimes_RN\otimes_RV$ or
$N\otimes_RB\otimes_RV$, which is relatively projective over the
appropriate enveloping algebra.  The same trace isomorphisms now
compare the relative resolutions and prove the relative assertion.
\end{proof}

\begin{theorem}\label{thm:hh-positive-homology}
Assume \eqref{eq:hh-mt}.  For every $n>0$ and for each convention
$\varepsilon$ there is an $R$-module isomorphism
\[
 \mathrm{HH}^{\varepsilon}_n(A/R)
       \simeq \mathrm{HH}^{\varepsilon}_n(B/R).
\]
No adjunction or $R$-projectivity hypothesis on the algebras is
required.
\end{theorem}

\begin{proof}
An ordinary projective enveloping-algebra module is also relatively
projective.  Therefore, for $n>0$,
\[
 \operatorname{Tor}^{A^e,\varepsilon}_n(A,P)=0,
 \qquad
 \operatorname{Tor}^{B^e,\varepsilon}_n(B,Q)=0.
\]
Using the two decompositions in \eqref{eq:hh-mt} and
Lemma~\ref{lem:hh-cyclic-tor}, we obtain
\begin{align*}
 \mathrm{HH}^{\varepsilon}_n(A/R)
 &\simeq \mathrm{HH}^{\varepsilon}_n(A/R;A\oplus P)\\
 &\simeq \mathrm{HH}^{\varepsilon}_n(A/R;M\otimes_BN)\\
 &\simeq \mathrm{HH}^{\varepsilon}_n(B/R;N\otimes_AM)\\
 &\simeq \mathrm{HH}^{\varepsilon}_n(B/R;B\oplus Q)
  \simeq \mathrm{HH}^{\varepsilon}_n(B/R).
\end{align*}
\end{proof}

\subsection{The image of positive cohomology under tensoring}

Projectivity of $Q$ does not imply
$\operatorname{Ext}_{B^e}^n(B,Q)=0$ for $n>0$.
Thus projectivity alone does not eliminate the off-diagonal entries in
$\operatorname{Ext}_{B^e}^n(B\oplus Q,B\oplus Q)$.
We prove the required vanishing for classes in the image of $F$.

We write $\alpha:X\to Y[n]$ for an Ext class in either
the ordinary derived category or the derived category of the
$R$-split exact structure.  Equivalently, the arguments may be expressed
in terms of Yoneda extensions and exact functors.

\begin{lemma}
\label{lem:hh-annihilation}
Let $E$ be an $A^e$-projective bimodule which is projective as a left
$A$-module.  Let $n>0$ and
$\alpha\in\mathrm{HH}_{\varepsilon}^{n}(A/R)$.
For bimodule maps $u:A\to E$ and $v:E\to A$ one has
\[
 u[n]\alpha=0,\qquad \alpha v=0.
\]
Consequently the morphism
$\operatorname{diag}(\alpha,0):A\oplus E\to(A\oplus E)[n]$
commutes, with the appropriate shift, with every degree-zero
endomorphism of $A\oplus E$.
\end{lemma}

\begin{proof}
Tensoring on the right with $E$ is an exact functor because
${}_AE$ is projective.  It is also exact for the $R$-split exact
structure, by the retract argument above.  The map $u$ defines a
natural transformation
\[
 -\otimes_Au:\operatorname{Id}=-\otimes_AA
                         \longrightarrow -\otimes_AE
\]
between these exact functors on $A^e$-modules.  Applying this
transformation degreewise to complexes and passing to the derived
category gives a natural transformation of derived functors.
Applied to $\alpha$, naturality yields
\[
 u[n]\alpha=(\alpha\otimes_AE)u.
\]
The class $\alpha\otimes_AE$ belongs to
$\operatorname{Ext}_{A^e,\varepsilon}^n(E,E)=0$, since $E$ is
projective in either exact structure.  This proves the first
vanishing.  The second follows directly from
$\operatorname{Ext}_{A^e,\varepsilon}^n(E,A)=0$.

Finally, every bimodule endomorphism of $A$ is multiplication by an
element of $Z(A)$.  Such multiplication extends to a natural
endomorphism of the identity functor on $A^e$-modules, and therefore
commutes with every Ext class.  In a block matrix endomorphism of
$A\oplus E$, the two off-diagonal commutation identities follow from
the vanishings above, and the upper-left identity follows from
centrality.
\end{proof}

\begin{lemma}\label{lem:hh-image-compression}
Fix the decomposition $F(A)=B\oplus Q$, with inclusion
$i_B:B\to F(A)$ and projection $p_B:F(A)\to B$.  For every $n>0$
and $\alpha\in\mathrm{HH}_{\varepsilon}^{n}(A/R)$, there is a
unique $\phi_n(\alpha)\in\mathrm{HH}_{\varepsilon}^{n}(B/R)$
such that
\begin{equation}
 F(\alpha)=i_B[n]\phi_n(\alpha)p_B,
 \qquad
 \phi_n(\alpha)=p_B[n]F(\alpha)i_B.
 \label{eq:hh-image-compression}
\end{equation}
In particular, the map $F(\alpha)$ has zero off-diagonal entries in
this decomposition.  The analogous assertion holds for $G$.
\end{lemma}

\begin{proof}
By \eqref{eq:hh-gf-splitting},
\[
 GF(A)\simeq A\oplus E_A,
 \qquad E_A=P\oplus P\oplus(P\otimes_AP).
\]
The bimodule $E_A$ is projective on both sides and is
$A^e$-projective.  For the latter assertion, $P\otimes_AP$ is a
direct summand of a finite direct sum of copies of $A\otimes_RP$,
which is $A^e$-projective because $P_A$ is projective.

Since $K_A$ is exact and $K_A(A)=E_A$ is projective, the morphism
$K_A(\alpha)$ lies in
$\operatorname{Ext}_{A^e,\varepsilon}^n(E_A,E_A)=0$.
Thus in the natural decomposition above,
\[
 GF(\alpha)=\operatorname{diag}(\alpha,0).
\]
Let $e=i_Bp_B\in\operatorname{End}_{B^e}(F(A))$.
Lemma~\ref{lem:hh-annihilation} says that $GF(\alpha)$ commutes
with $G(e)$.  Since $G$ is faithful on the relevant derived
category, this implies
\[
 F(\alpha)e=e[n]F(\alpha).
\]
The restriction of $F(\alpha)$ to its source summand $Q$ is zero,
because $Q$ is projective and $n>0$.  Hence
$F(\alpha)=F(\alpha)e=e[n]F(\alpha)$, which is exactly
\eqref{eq:hh-image-compression}.  Uniqueness follows by composing
with $p_B[n]$ and $i_B$.
\end{proof}

\subsection{The positive cohomology algebra and the stable degree zero}

Let $Z^{\mathrm{pr}}(A)$ denote the ideal of $Z(A)$ consisting of
bimodule endomorphisms $A\to A$ which factor through an ordinary
projective $A^e$-module.  In either cohomology convention this ideal
annihilates positive cohomology.  Indeed, if
$z:A\to T\to A$ with $T$ ordinary projective and $\alpha$ has
positive degree, then $\alpha z=0$, because
$\operatorname{Ext}_{A^e,\varepsilon}^{>0}(T,A)=0$.
Centrality of the degree-zero action also gives $z\alpha=0$.
We may therefore form the graded algebra
\begin{equation}
 \overline{\mathrm{HH}}_{\varepsilon}^{*}(A/R)
   =\mathrm{HH}_{\varepsilon}^{*}(A/R)/Z^{\mathrm{pr}}(A),
 \label{eq:hh-stable-ring}
\end{equation}
where the quotient changes only degree zero.
The same degree-zero ideal is obtained by allowing factorization
through all relatively $R$-projective bimodules.  Indeed, the
multiplication $A^e\to A$ is an $R$-split epimorphism.  A map from a
relatively projective bimodule to $A$ lifts across it, so every such
factorization of an endomorphism of $A$ already factors through
the ordinary free bimodule $A^e$; see
Lemma~\ref{lem:unit-relative-projective-center}.
This equality of the projective-center ideals does not identify
ordinary and relative Hochschild cohomology in positive degrees
when $A$ is not projective over $R$.

\begin{theorem}\label{thm:hh-cohomology-algebra}
Assume \eqref{eq:hh-mt}.  For each
$\varepsilon\in\{\mathrm{abs},\mathrm{rel}\}$ there is an
isomorphism of graded $R$-algebras
\[
 \overline{\mathrm{HH}}_{\varepsilon}^{*}(A/R)
       \simeq
 \overline{\mathrm{HH}}_{\varepsilon}^{*}(B/R).
\]
In particular,
\[
 \bigoplus_{n>0}\mathrm{HH}_{\varepsilon}^{n}(A/R)
       \simeq
 \bigoplus_{n>0}\mathrm{HH}_{\varepsilon}^{n}(B/R)
\]
as positively graded, generally nonunital $R$-algebras.  Thus the
isomorphism preserves every positive cohomology group and all products
of positive-degree classes.
\end{theorem}

\begin{proof}
Choose split maps
\[
 i_B:B\to F(A),\quad p_B:F(A)\to B,
 \qquad
 i_A:A\to G(B),\quad p_A:G(B)\to A
\]
from the decompositions in \eqref{eq:hh-mt}.  In positive degree
define $\phi$ by \eqref{eq:hh-image-compression}; define $\psi$
using $G$ and the second pair of split maps.  Exact functors
preserve Yoneda products.  Since every positive $F(\alpha)$ is
supported on the $B$-summand, compression is multiplicative on
this image.  Thus $\phi$ is a homomorphism on positive cohomology,
and likewise for $\psi$.

In degree zero define
\[
 \phi_0([a])=[p_BF(a)i_B],
 \qquad
 \psi_0([b])=[p_AG(b)i_A].
\]
These maps are well-defined modulo the projective centers because
$F$ and $G$ preserve ordinary projective bimodules.  For two
degree-zero endomorphisms, the difference between compressing
their product and multiplying their compressions factors through
$Q$, or respectively $P$.  It therefore vanishes in the quotient.
For a product with a positive-degree factor, multiplicativity
follows from the block-diagonal identity
\eqref{eq:hh-image-compression}.  Thus $\phi$ and $\psi$ are
homomorphisms of the graded algebras in
\eqref{eq:hh-stable-ring}.

It remains to prove that they are inverse.  Denote the natural
split maps of $GF(A)=A\oplus E_A$ in
\eqref{eq:hh-gf-splitting} by
\[
 i:A\longrightarrow GF(A),\qquad
 p:GF(A)\longrightarrow A.
\]
There is another split pair
\[
 j=G(i_B)i_A:A\longrightarrow GF(A),
 \qquad
 q=p_AG(p_B):GF(A)\longrightarrow A,
 \qquad qj=1_A.
\]
For $\alpha$ of positive degree, put
$\gamma=\psi\phi(\alpha)$.  The natural decomposition gives
$GF(\alpha)=i[n]\alpha p$.
Applying Lemma~\ref{lem:hh-image-compression} first to $F$ and
then to $G$ gives the other expression
$GF(\alpha)=j[n]\gamma q$.
Compressing with $q[n]$ and $j$ yields
\[
 \gamma=(qi)[n]\alpha(pj).
\]
Set $u=qi$ and $v=pj$.  These are elements of $Z(A)$, so
$u[n]\alpha=\alpha u$.  Also
\[
 uv=qi\,pj=qj-q(1-ip)j
          =1_A-z,
 \qquad z=q(1-ip)j\in Z^{\mathrm{pr}}(A),
\]
because $1-ip$ is the projection onto the ordinary projective
summand $E_A$.  Consequently
\[
 \gamma=\alpha uv=\alpha(1_A-z)=\alpha.
\]
Interchanging $A$ and $B$ proves $\phi\psi=\operatorname{Id}$
in positive degree.

For a degree-zero endomorphism $a$ one still has
$\psi_0\phi_0([a])=[qGF(a)j]$.
The natural decomposition now writes $GF(a)$ as
$\operatorname{diag}(a,K_A(a))$.  Its second component factors
through $E_A$, and hence contributes zero modulo
$Z^{\mathrm{pr}}(A)$.  Therefore
\[
 [qGF(a)j]=[uav]=[a(1_A-z)]=[a].
\]
The same argument on $B$ proves the inverse identity in degree
zero and completes the proof.
\end{proof}

\begin{remark}
In general, compression does not identify
$\operatorname{Ext}_{B^e}^{>0}(B\oplus Q,B\oplus Q)$ with
$\operatorname{Ext}_{B^e}^{>0}(B,B)$, and $Q$ need not be injective.
Lemma~\ref{lem:hh-image-compression} asserts that the off-diagonal
components vanish for the classes $F(\alpha)$.  Faithfulness is
supplied by the natural direct-summand decompositions of the derived
functors.  To compare Gerstenhaber brackets, we also need the
common-image and triangular-cochain arguments of
Section~\ref{sec:relative-gerstenhaber}.  Preservation of a
Batalin--Vilkovisky operator is not a consequence of this argument.
\end{remark}

\subsection{Degree zero and transfer maps}

\begin{example}
\label{ex:hh-degree-zero}
Let $R$ be any nonzero commutative ring, put $A=R$ and
$B=R\times R$, and take $M=(R,R)$ with coordinatewise right
$B$-action and scalar left $A$-action.  Take $N=(R,R)^{\mathrm t}$
with coordinatewise left $B$-action and scalar right $A$-action.
These bimodules are projective on both sides.  One has
\[
 M\otimes_BN\simeq R\oplus R=A\oplus P,
 \qquad P=R\in A^e\text{-}\mathrm{proj}.
\]
Writing $e_1,e_2$ for the two central idempotents of $B$, one also
has
\[
 N\otimes_AM\simeq B\oplus Q,
 \qquad
 Q=(Be_1\otimes_Re_2B)\oplus(Be_2\otimes_Re_1B),
\]
and $Q$ is $B^e$-projective.  Thus \eqref{eq:hh-mt} holds, even
with both algebras projective over $R$.  Nevertheless
\[
 \mathrm{HH}^{0}(A/R)=R,\qquad
 \mathrm{HH}^{0}(B/R)=R\oplus R,
\]
and the same difference occurs for $\mathrm{HH}_{0}$, since the
two algebras are commutative.  The distinction holds in both
conventions.  Their stable centers are both zero, consistently
with Theorem~\ref{thm:hh-cohomology-algebra}.
\end{example}

For finite-dimensional self-injective algebras, Pogorza{\l}y's
comparison \cite{Pogorzaly2001Hochschild} preserves the Yoneda
product in positive degrees. Its stable-morphism formulation has
\[
 \underline{\operatorname{End}}_{A^e}(A)
       =Z(A)/Z^{\mathrm{pr}}(A)
\]
in degree zero, whereas ordinary Hochschild cohomology has
$\mathrm{HH}^0(A)=Z(A)$. Thus the stable formulation must not be
read as invariance of the ordinary center. The distinction is also
explicit in \cite[Section~1]{Pogorzaly2014Orbit}, where the
syzygy orbit algebra is described as a quotient of ordinary
Hochschild cohomology. Neither the Yoneda-product comparison nor
the orbit-algebra results alone give compatibility with the
Gerstenhaber bracket.

Liu and Xi proved positive homology invariance for
finite-dimensional algebras over a field in
\cite[Theorem~4.4]{LiuXi2005}; the cyclic Tor argument above is
a direct proof of the corresponding comparison over the present base
ring.  Xi's theorem
\cite[Theorem~4.2]{Xi2008} assumes Artin algebras projective over
the base and a stable equivalence of adjoint type.  For arbitrary
finite-dimensional algebras over a field, Pan and Zhou subsequently
proved preservation of the stable Hochschild cohomology graded
ring without an adjoint-type hypothesis
\cite[Theorem~1.2]{PanZhou2010}.  The preceding proof treats an arbitrary
commutative base and keeps ordinary enveloping-algebra Ext distinct from
relative bar cohomology.  Its specialization to a field recovers these
known invariance results.

Bouc constructed homology transfers associated with suitable
one-sided projective bimodules \cite{Bouc1997}; their composition,
additivity, and positive-degree vanishing for projective bimodules
are recorded in \cite[Proposition~2.1]{KoenigLiuZhou2012}.  Those
properties already give positive homology invariance over a
field without adjunction.  Linckelmann's cohomology transfer in
\cite[Section~2]{Linckelmann1999} is constructed for symmetric
algebras over a commutative ring, using adjunction maps and
symmetrizing forms.  The compression maps used in the proof above
are defined directly from \eqref{eq:hh-mt} and do not require
those additional data.  Identifying them with a normalized
cohomological transfer requires a separate comparison.

\begin{corollary}
\label{cor:hh-projective-base}
If $A$ and $B$ are projective as $R$-modules and satisfy
\eqref{eq:hh-mt}, their ordinary enveloping-algebra and relative
bar definitions agree.  In the customary notation, for every
$n>0$,
\[
 \mathrm{HH}_n(A/R)\simeq\mathrm{HH}_n(B/R),
 \qquad
 \mathrm{HH}^{n}(A/R)\simeq\mathrm{HH}^{n}(B/R),
\]
and the isomorphisms in cohomology form an isomorphism of the
positive graded algebras, extending to an isomorphism after
replacing degree zero by the stable center.  No adjoint-type
hypothesis is needed.
\end{corollary}

\subsection{A correction to the compression argument of Pan--Zhou}
\label{sec:pan-zhou-correction}

The graded-algebra invariance theorem of Pan--Zhou
\cite[Theorem~1.2]{PanZhou2010} holds under its original hypotheses.
Claims~1 and~3 in its proof require correction for arbitrary graded
endomorphisms of a bimodule with a projective direct summand in the
derived category. We prove the multiplicativity and inverse identities
on the tensor images used in the proof.

Assume \eqref{eq:hh-mt} and use the functors $F,G$ above. Let
$i_B:B\to F(A)$ and $p_B:F(A)\to B$ be the split maps for the
summand $B$ of $F(A)\simeq B\oplus Q$. Define
\[
 \phi_n(\alpha)=p_B[n]F(\alpha)i_B\qquad(n>0),
 \qquad
 \phi_0([a])=[p_BF(a)i_B].
\]
Define $\psi$ using $G$ and the split maps $i_A,p_A$ for the
summand $A$ of $G(B)\simeq A\oplus P$.

\begin{proposition}\label{prop:pan-zhou-compression}
Assume \eqref{eq:hh-mt} and let
$\varepsilon\in\{\mathrm{abs},\mathrm{rel}\}$. The compression maps
\[
 \overline{\mathrm{HH}}_{\varepsilon}^{*}(A/R)
 \xrightarrow{\ \phi\ }
 \overline{\mathrm{HH}}_{\varepsilon}^{*}(B/R),
 \qquad
 \overline{\mathrm{HH}}_{\varepsilon}^{*}(B/R)
 \xrightarrow{\ \psi\ }
 \overline{\mathrm{HH}}_{\varepsilon}^{*}(A/R)
\]
are well-defined homomorphisms of graded $R$-algebras and
satisfy $\psi\phi=1$ and $\phi\psi=1$.
\end{proposition}

\begin{proof}
By Lemma~\ref{lem:hh-image-compression}, for $n>0$ one has
$F(\alpha)=i_B[n]\phi_n(\alpha)p_B$. Let $\alpha$ and $\beta$ have degrees $n$ and $m$. If $m+n>0$,
inserting $i_Bp_B$ between the two tensor images leaves the
compressed composite unchanged. In fact, if $n>0$ then
$F(\alpha)i_B$ factors through $i_B[n]$, and if $m>0$ then
$p_B[m]F(\beta)$ factors through $p_B$. Hence
\[
 p_B[m+n]F(\beta)[n]F(\alpha)i_B
  =\bigl(p_B[m]F(\beta)i_B\bigr)[n]\,
       p_B[n]F(\alpha)i_B .
\]
This equality also uses the indicated representatives when
one degree is zero. For two degree-zero factors $a,b$, the
difference is
\[
 p_BF(b)F(a)i_B-(p_BF(b)i_B)(p_BF(a)i_B)
   =p_BF(b)(1-i_Bp_B)F(a)i_B.
\]
It factors through $Q$. Also, $F$ sends projective bimodules to
projective bimodules, so $\phi_0$ is well-defined modulo
$Z^{\mathrm{pr}}(B)$. The same reasoning applies to $\psi$.
The identities $p_Bi_B=1_B$ and $p_Ai_A=1_A$ show that both
compressions preserve the units.

For the successive compressions, let $i,p$ be the split maps for
$A$ in the natural decomposition $GF(A)\simeq A\oplus E_A$ of
\eqref{eq:hh-gf-splitting}. Put
\[
 j=G(i_B)i_A,\qquad q=p_AG(p_B),\qquad
 u=qi,\qquad v=pj.
\]
Since $qj=1_A$, one has
\[
 uv=1_A-z,\qquad z=q(1-ip)j\in Z^{\mathrm{pr}}(A).
\]
Indeed $1-ip$ factors through the projective summand $E_A$.
For $n>0$, the natural decomposition gives
$GF(\alpha)=i[n]\alpha p$, and hence
\[
 \psi_n\phi_n(\alpha)
   =q[n]GF(\alpha)j=u[n]\alpha v
   =\alpha(1_A-z)=\alpha.
\]
Here $u,v$ are bimodule endomorphisms of $A$. Their action is
central, and $Z^{\mathrm{pr}}(A)$ annihilates positive cohomology.
In degree zero, the component of $GF(a)$ on $E_A$ contributes
only a map factoring through $E_A$. Thus
\[
 \psi_0\phi_0([a])=[uav]=[a(1_A-z)]=[a].
\]
Interchanging $A$ and $B$ proves the other inverse identity.
\end{proof}

\section{The relative Gerstenhaber bracket}
\label{sec:relative-gerstenhaber}

Let $R$ be an arbitrary commutative ring. We use the categories of
all modules, as in Section~\ref{sec:hochschild-without-adjunction}.
Write $\mathrm{HH}^n_R(A)$ for
$\mathrm{HH}^n_{\mathrm{rel}}(A/R)$, computed by the relative
Hochschild cochain complex. All relative Ext groups in this section
are taken with respect to the $R$-split exact structure. The
Gerstenhaber bracket is defined by insertion of cochains.

We first show that the isomorphism in positive degrees from
Theorem~\ref{thm:hh-cohomology-algebra} preserves the Gerstenhaber
bracket. We identify its graph with the restriction image from a
triangular algebra, using the injectivity and equality of the
images of the two characteristic maps into a mixed bimodule Ext
group. We then treat the stabilized identities and the restricted
power in prime characteristic.

\begin{theorem}
\label{thm:relative-gerstenhaber}
Assume \eqref{eq:hh-mt}.  The compression isomorphism
\[
 \phi_n(\alpha)=p_B[n]
       (N\otimes_A\alpha\otimes_A M)i_B,
 \qquad n>0,
\]
from Theorem~\ref{thm:hh-cohomology-algebra} satisfies
\[
 \phi_{m+n-1}([\alpha,\beta]_A)
   =[\phi_m(\alpha),\phi_n(\beta)]_B
 \quad
 (\alpha\in\mathrm{HH}^m_R(A),\ 
  \beta\in\mathrm{HH}^n_R(A),\ m,n>0).
\]
Consequently $\mathrm{HH}^{>0}_R(A)$ and
$\mathrm{HH}^{>0}_R(B)$ are isomorphic as nonunital Gerstenhaber
algebras.  This conclusion requires neither projectivity over $R$
nor adjointness or preservation of injectives for the tensor functors.
The functors need not induce an equivalence of the ordinary stable
module categories.

If $A$ and $B$ are projective over $R$, the same assertion holds for
ordinary $\operatorname{Ext}_{A^e}^*(A,A)$ and
$\operatorname{Ext}_{B^e}^*(B,B)$ in positive degrees.  In particular,
it holds for classical stable equivalences of Morita type between
finite-dimensional algebras over an arbitrary field.
\end{theorem}

\subsection{The two characteristic actions have the same image}

Put $C=A\otimes_R B^{\mathrm{op}}$.  Exact tensoring defines maps
\begin{align*}
 a_n:\mathrm{HH}^n_R(A)&\longrightarrow
          \operatorname{Ext}_{C,\mathrm{rel}}^n(M,M),
 &a_n(\alpha)&=\alpha\otimes_A M,\\
 b_n:\mathrm{HH}^n_R(B)&\longrightarrow
          \operatorname{Ext}_{C,\mathrm{rel}}^n(M,M),
 &b_n(\beta)&=M\otimes_B\beta.
\end{align*}
These are the maps on Yoneda extensions, or equivalently on the
derived categories of the indicated exact structures.
All tensor functors here preserve $R$-split conflations by the retract
argument of Lemma~\ref{lem:hh-tensor-properties}.

\begin{lemma}
\label{lem:gerstenhaber-common-image}
For every $n>0$, the maps $a_n,b_n$ are injective and
\begin{equation}
 a_n=b_n\phi_n.
 \label{eq:gerstenhaber-common-image}
\end{equation}
In particular, their images coincide, and
\[
 \ker(a_n-b_n)
   =\{(\alpha,\phi_n(\alpha)):
                        \alpha\in\mathrm{HH}^n_R(A)\}.
\]
The common image may be a proper submodule of
$\operatorname{Ext}_{C,\mathrm{rel}}^n(M,M)$.
\end{lemma}

\begin{proof}
The composite of $-\otimes_A M$ and $-\otimes_B N$ is
$\operatorname{Id}\oplus(-\otimes_A P)$.  Hence the first functor induces injective maps on Ext groups.
This proves injectivity of $a_n$.  Tensoring on the left proves
injectivity of $b_n$.

Consider
\[
 L=N\otimes_A-:C\text{-}\mathrm{Mod}
                        \longrightarrow B^e\text{-}\mathrm{Mod}.
\]
Its composite with $M\otimes_B-$ contains the identity functor as
a natural direct summand, since $M\otimes_BN=A\oplus P$.
Thus $L_*$ is injective on all relevant relative Ext groups.
For $\alpha\in\mathrm{HH}^n_R(A)$, set $\beta=\phi_n(\alpha)$.
Lemma~\ref{lem:hh-image-compression} gives
\[
 L_*(a_n(\alpha))
 =N\otimes_A\alpha\otimes_A M
 =\begin{pmatrix}\beta&0\\0&0\end{pmatrix}
 \quad\text{on }N\otimes_AM=B\oplus Q.
\]
The fixed decomposition of this bimodule also gives
\[
 L_*(b_n(\beta))
 =(N\otimes_AM)\otimes_B\beta
 =\begin{pmatrix}\beta&0\\0&Q\otimes_B\beta\end{pmatrix}
 =\begin{pmatrix}\beta&0\\0&0\end{pmatrix}.
\]
The last equality follows from
$\operatorname{Ext}_{B^e,\mathrm{rel}}^n(Q,Q)=0$ for $n>0$,
since $Q$ is relatively projective.
Injectivity of $L_*$ proves \eqref{eq:gerstenhaber-common-image}.
Since $\phi_n$ is an isomorphism, the images coincide.  Finally,
$a_n(\alpha)=b_n(\beta)$ implies
$b_n(\phi_n(\alpha)-\beta)=0$, and injectivity of $b_n$ proves
the stated kernel formula.
\end{proof}

\subsection{A relative triangular comparison}
\label{subsec:relative-triangular-comparison}

Let
\[
 T=\begin{pmatrix}A&M\\0&B\end{pmatrix},\qquad
 e=\begin{pmatrix}1&0\\0&0\end{pmatrix},\qquad f=1-e.
\]
Guccione--Guccione \cite[Theorem~1]{GuccioneGuccione2001}
established the relative triangular long exact sequence over a
commutative base.  We describe its maps and verify the compatibility
of corner restriction with the Gerstenhaber bracket.

Keller \cite[Section~4.5]{Keller2003Higher} uses a triangular
comparison whose corner restriction maps preserve the
$B_\infty$ structure.
For triangular algebras projective over the base, Hermann
\cite[Section~6]{Hermann2016Recollements} proves compatibility
with the strict Gerstenhaber structure.  A corresponding triangular
comparison in stable derived categories is used in
\cite{BriggsRubio2022}.  The additional step here is the
common-image identity \eqref{eq:gerstenhaber-common-image}, which
identifies the graph of the positive-degree compression map.
The relative bar construction below supplies the required corner
comparison without projectivity over the base.

\medskip
\noindent\textit{Keller admissible triples and the singular comparison.}
Liao and Seol \cite{LiaoSeol2023} give a direct-sum version
of Keller's triangular comparison for dg algebras. Their
admissibility conditions include weakly cone-nilpotent
quasi-isomorphisms from the two algebras to the appropriate
one-sided dg endomorphism algebras of the mixed bimodule,
as well as exactness of the partial Hochschild complexes.
The two restrictions then induce Gerstenhaber isomorphisms
on direct-sum Hochschild cohomology. They apply the comparison
to the enveloping algebra and the Chevalley--Eilenberg algebra
of a finite-dimensional Lie algebra to obtain the
Duflo--Kontsevich isomorphism.

For ordinary algebras and a defining pair $(M,N)$ for a
stable equivalence of Morita type, the action map
$A\to\operatorname{End}_{B^{\mathrm{op}}}(M)$ is a map
between algebras concentrated in degree zero. Its being
a quasi-isomorphism would make $M$ a Morita bimodule.
Indeed, $B_B$ is a summand of $N\otimes_A M$, which is a
summand of a finite direct sum of copies of $M_B$ by
projectivity of $N_A$. Thus $M_B$ is a progenerator, and
the endomorphism algebra condition gives a Morita equivalence.
The deformation comparison of
Theorem~\ref{thm:stable-associative-deformations} uses the
two triangular restrictions in degrees one, two and three
with all square-zero coefficients.

There is also a comparison in the singular category of
bimodules. Let $A,B$ be finite-dimensional over a field and
put $C=A\otimes_k B^{\mathrm{op}}$. Write
$\operatorname{End}_{\mathrm{sg},D}^{\mathrm{dg}}(X)$ for the
endomorphism complex of $X$ in the dg singularity category of
$D$. The exact tensor functors give a span of quasi-isomorphisms
of dg algebras
\[
 \operatorname{End}_{\mathrm{sg},A^e}^{\mathrm{dg}}(A)
 \longrightarrow \operatorname{End}_{\mathrm{sg},C}^{\mathrm{dg}}(M)
 \longleftarrow \operatorname{End}_{\mathrm{sg},B^e}^{\mathrm{dg}}(B).
\]
For the first arrow, $-\otimes_A M$ and $-\otimes_B N$
preserve perfect bimodule complexes and descend to the dg
quotients. On the full subcategories of objects represented
by complexes perfect over the left copy of $A$, the
composites are isomorphic to the identities in the singular
categories. Indeed, the error $X\otimes_A P$ is a retract of
a finite direct sum of $X\otimes_k A$, which is perfect over
$A^e$ when $X$ is perfect on the left. The same argument uses
$Q$ for the other composite. This gives isomorphisms on the
cohomology of the endomorphism complexes of $A$ and $M$ in
all degrees. For the second arrow, use the functors
$M\otimes_B-$ and $N\otimes_A-$ on complexes perfect over the
right copy of $B$. For the brace operations we use the
singular Hochschild cochain comparison in
Section~\ref{subsec:relative-singular-braces}.

\begin{lemma}
\label{lem:relative-triangular-gerstenhaber}
There is a long exact sequence
\begin{equation}
 \begin{aligned}
 \cdots\longrightarrow\mathrm{HH}^n_R(T)
  &\xrightarrow{\ r_n=(r_A,r_B)\ }
      \mathrm{HH}^n_R(A)\oplus\mathrm{HH}^n_R(B)\\
  &\xrightarrow{\ a_n-b_n\ }
      \operatorname{Ext}_{C,\mathrm{rel}}^n(M,M)
      \longrightarrow\mathrm{HH}^{n+1}_R(T)
      \longrightarrow\cdots .
 \end{aligned}
 \label{eq:relative-triangular-les}
\end{equation}
The two corner maps $r_A,r_B$ preserve the cup product and the
Gerstenhaber bracket.  These assertions do not require $T$, $A$, $B$
or $M$ to be projective over $R$.
\end{lemma}

\begin{proof}
Consider the short exact sequence of $T$-bimodules
\begin{equation}
 \begin{aligned}
 0\longrightarrow Te\otimes_A M\otimes_B fT
  &\xrightarrow{\ (+,-)\ }
       (Te\otimes_A eT)\oplus(Tf\otimes_B fT)\\
  &\xrightarrow{\ \mathrm{mult}\ }
       T\longrightarrow0.
 \end{aligned}
 \label{eq:triangular-bimodule-sequence}
\end{equation}
The kernel is the difference of the two copies of the upper-right
corner $M$.  The sequence splits as a sequence of $R$-modules.
The induction functors
\[
 \begin{split}
 L_A(X)&=Te\otimes_A X\otimes_A eT,\\
 L_B(Y)&=Tf\otimes_B Y\otimes_B fT,\\
 L_C(Z)&=Te\otimes_A Z\otimes_B fT
 \end{split}
\]
are exact and preserve $R$-split conflations, since $M$ is
projective on each side.
They also preserve relative projectives.  For example,
\[
 L_A(A^e\otimes_R V)\simeq Te\otimes_R V\otimes_R eT
\]
is a direct summand of a finite direct sum of copies of
$T^e\otimes_R V$, since $Te$ and $eT$ are projective on the
relevant $T$-sides.  The arguments for $L_B$ and $L_C$ are the same.
The induction--restriction adjunctions therefore identify
the Ext terms obtained from \eqref{eq:triangular-bimodule-sequence}
with those in \eqref{eq:relative-triangular-les}.  The two maps to the $M$-term are induced by tensoring a Hochschild
extension with $M$ on the right and on the left, respectively.
They occur with opposite signs, giving the map $a_n-b_n$.

To compare brackets, let $E=R\times R$ act on $T$ by
$(r,s)\mapsto re+sf$.  This map need not be injective.
The $E$-relative bar resolution of $T$ is a projective resolution in
the $R$-split exact structure.  Indeed, its terms decompose along
the idempotents $e,f$ into direct sums of bimodules
\[
 Te_i\otimes_R V\otimes_R e_jT,
 \qquad e_i,e_j\in\{e,f\},
\]
which are relatively $R$-projective, and the augmented bar complex
has an $R$-linear contracting homotopy.  Equivalently, separability of $E$ over $R$ realizes this resolution
as the direct summand of the $R$-relative bar resolution indexed by
composable sequences of corners.  No $R$-flatness is required.

Thus
\[
 C_E^n(T,T)=\operatorname{Hom}_{E^e}
                      (T^{\otimes_E n},T)
\]
computes $\mathrm{HH}_R^*(T)$.  Its inclusion into the usual
$R$-relative Hochschild cochain complex is a quasi-isomorphism and
preserves cup products and brace operations.  Restriction to a pure
$A$-corner, or to a pure $B$-corner, also preserves these operations:
$E$-bilinearity ensures that the value of a cochain on a pure corner
lies in the same corner, so all the insertion formulas restrict
term by term.  These restrictions induce the maps $r_A,r_B$ in
\eqref{eq:relative-triangular-les}.  In particular they preserve
the bracket.  When the algebras are projective over the base, this compatibility
also appears in \cite[pp.~372--373]{BendiffalahGuin2006}.
The bar argument proves the relative version required here.
\end{proof}

\subsection{Compatibility with the Gerstenhaber bracket}

\begin{proof}[Proof of Theorem~\ref{thm:relative-gerstenhaber}]
Lemmas~\ref{lem:gerstenhaber-common-image} and
\ref{lem:relative-triangular-gerstenhaber} show that, for $n>0$,
\begin{equation}
 \operatorname{im}r_n
   =\{(\alpha,\phi_n(\alpha)):
                       \alpha\in\mathrm{HH}^n_R(A)\}.
 \label{eq:relative-gerstenhaber-graph}
\end{equation}
Take $\alpha\in\mathrm{HH}^m_R(A)$ and
$\beta\in\mathrm{HH}^n_R(A)$, with $m,n>0$.  Choose
$x\in\mathrm{HH}^m_R(T)$ and $y\in\mathrm{HH}^n_R(T)$ such that
\[
 r_m(x)=(\alpha,\phi_m(\alpha)),\qquad
 r_n(y)=(\beta,\phi_n(\beta)).
\]
The bracket compatibility of the two restriction maps gives
\[
 r_{m+n-1}([x,y])
   =\bigl([\alpha,\beta],
           [\phi_m(\alpha),\phi_n(\beta)]\bigr).
\]
Since $m+n-1>0$, this pair belongs to the graph in
\eqref{eq:relative-gerstenhaber-graph}.  Its second component must
therefore be $\phi_{m+n-1}([\alpha,\beta])$, as required.
The identical argument with $xy$ proves cup compatibility, in
agreement with Theorem~\ref{thm:hh-cohomology-algebra}.

The argument uses the existence of lifts of graph elements and does
not require $r_n$ to be injective, even for $n=1$.  When $A$ and $B$
are $R$-projective, their relative bar resolutions are ordinary
enveloping-algebra projective resolutions and the two cohomology
conventions, including the bracket, agree.
\end{proof}

\subsection{Relative projective error bimodules}

\begin{proposition}
\label{prop:gerstenhaber-relative-errors}
The relative conclusion of Theorem~\ref{thm:relative-gerstenhaber}
remains valid when ordinary enveloping-algebra projectivity of $P,Q$
is replaced by
\[
 P\mid A^e\otimes_R V,\qquad
 Q\mid B^e\otimes_R W
\]
for $R$-modules $V,W$, provided that the decompositions
$M\otimes_BN\simeq A\oplus P$ and
$N\otimes_AM\simeq B\oplus Q$ still hold and $M,N$ are
projective on both sides.  The corresponding graded-ring
isomorphism has degree zero given by the quotients of $Z(A)$ and
$Z(B)$ by maps factoring through relatively $R$-projective
bimodules.
\end{proposition}

\begin{proof}
The decompositions imply that $P,Q$ are projective on each side,
as in Lemma~\ref{lem:hh-tensor-properties}.  The tensor functors
$F,G$ preserve relative projectives.  The additional term $P\otimes_AP$ in the composite error is
also relatively projective: if $P\mid A^e\otimes_R V$, then
\[
 P\otimes_AP\mid A\otimes_R V\otimes_R P,
\]
and the latter is a direct summand of a finite direct sum of copies
of $A^e\otimes_R V$, since $P_A$ is projective.

Lemma~\ref{lem:hh-annihilation} uses only projectivity in the
chosen exact structure and one-sided projectivity.  Its proof,
the compression proof, and the positive-degree inverse proof in
Theorem~\ref{thm:hh-cohomology-algebra} therefore work with
relative projectives throughout.  The central error in the inverse
calculation factors through a relative projective and hence
annihilates positive relative Ext.  This proves the ring isomorphism
in positive degrees.  The common-image and graph arguments then
give bracket compatibility.
\end{proof}

\subsection{Stable bimodule identities over arbitrary bases}

\begin{theorem}\label{cor:gerstenhaber-stable-identities}
Let $R$ be any commutative ring and let the central bimodules
$ {}_AM_B,{}_BN_A$ be finitely generated projective on both sides.
Suppose $M\otimes_BN\simeq A$ and $N\otimes_AM\simeq B$ in the
respective additive quotients by relatively $R$-projective
bimodules. Then positive relative Hochschild homology is invariant,
and $\mathrm{HH}^{>0}_R(A)\simeq\mathrm{HH}^{>0}_R(B)$.
The latter isomorphism preserves products and brackets and, if $R$
is an $\mathbb F_p$-algebra, the restricted power on shifted-even
classes. If the quotient identities use ordinary projective
bimodules, positive ordinary Hochschild homology is also invariant.
\end{theorem}

\begin{proof}
A stable isomorphism $X\simeq Y$ lifts to an isomorphism
$X\oplus E\simeq Y\oplus E'$ with relative projectives $E,E'$.
Indeed, if $gf-1_X=ba$ through a relative projective $E'$,
then $(f,a):X\to Y\oplus E'$ has retraction $(g,-b)$.
The complementary summand has zero identity in the quotient
and is relative projective. The same reasoning applies to
ordinary projectives. Write $U=M\otimes_BN$, $V=N\otimes_AM$,
and choose decompositions
\[
 U\oplus P_1\simeq A\oplus P_2,\qquad
 V\oplus Q_1\simeq B\oplus Q_2.
\]
Positive Tor vanishes on the corresponding projective errors,
so the homology assertions follow from Lemma~\ref{lem:hh-cyclic-tor}.

All remaining Ext groups in the proof are relative. Tensoring an
$R$-split conflation with an induced bimodule preserves an
$R$-linear splitting. By taking direct summands, tensoring with
a relative projective is therefore exact for the $R$-split
structures. The tensor product of two induced bimodules is
induced at the two endpoints, so the tensor product of two
relative projectives is relative projective. The same conclusion
holds when one factor is relative projective and the other is
projective on each side, by expressing the latter as a summand
of a free module on the remaining side.

Write $F$ and $G$ for the two-sided tensor functors of
Lemma~\ref{lem:hh-tensor-properties}. For $n>0$ and
$\alpha\in\mathrm{HH}^n_R(A)$ we show that $x=F(\alpha)\oplus0$
on $D=F(A)\oplus Q_1$ is supported on $B$ under
$D\simeq B\oplus Q_2$. Expanding
\[
 (U\oplus P_1)\otimes_A-\otimes_A(U\oplus P_1)
 \simeq (A\oplus P_2)\otimes_A-\otimes_A(A\oplus P_2)
\]
gives an isomorphism $GF\oplus H_1\simeq\operatorname{Id}\oplus H_2$
of exact functors. Each term of $H_1$ or $H_2$ has a relatively
projective outer factor and the other outer factor is
side-projective or relative projective. In particular $H_i(A)$
is relative projective. The reverse expansion gives
$FG\oplus K_1\simeq\operatorname{Id}\oplus K_2$. Each $K_i$
takes side-projective bimodules and relative projective
bimodules to relative projectives. Since $K_1(D)$ is relative projective,
$K_1$ acts by zero on $\operatorname{Ext}^n_{B^e,\mathrm{rel}}(D,D)$.
If $G(\xi)=0$ for a class in this group, then
$(FG\oplus K_1)(\xi)=0$. Its component on the identity functor
is $\xi$, which proves injectivity.

After adjoining zero on $H_1(A)$, the morphism $G(x)$ becomes
$\operatorname{diag}(\alpha,0)$ on
$A\oplus H_2(A)\oplus G(Q_1)$. Here $G(Q_1)$ is relative projective.
Lemma~\ref{lem:hh-annihilation} applies to this diagonal morphism,
since its proof uses only exactness of tensoring with a relative
projective in the $R$-split structure. The diagonal morphism commutes, with the appropriate shift,
with $G(e)\oplus0$, where $e$ is the idempotent of $D$ with
image $B$. The injectivity just proved gives $xe=e[n]x$.
Since $Q_2$ is relative projective, the component of $x$
with source $Q_2$ vanishes. Thus $x$ is the diagonal
extension of a class $\beta\in\mathrm{HH}^n_R(B)$.

Consider the two characteristic maps $a_n,b_n$ on $M$.
The map $L_*$ induced by $L=N\otimes_A-$ on positive self-Ext
of $M$ is injective: its composite with $M\otimes_B-$,
after adjoining $P_1\otimes_A-$, contains the identity,
and $P_i\otimes_A M$ is relative projective. Applying
$L$ to $a_n(\alpha)$ and adjoining zero on $Q_1$ gives $x$.
Applying it to $b_n(\beta)$ and adjoining zero on $Q_1$ gives
$(V\oplus Q_1)\otimes_B\beta$, since $Q_1\otimes_B\beta=0$.
Under the stabilized decomposition this is also
$\operatorname{diag}(\beta,0)$. Hence
$a_n(\alpha)=b_n(\beta)$. For the other inclusion, tensor
on the right by $N$ and use the compression for $G$ constructed
by interchanging $A,B$ above. Both characteristic maps are
injective: tensor with $N$ on the appropriate side, adjoin
the zero component on the stabilizing projective, and use
the regular summand of the resulting tensor action. Thus their
common image defines a unique isomorphism $\phi_n$.

The triangular exact sequence identifies the image of corner
restriction in each positive degree with the graph of $\phi_n$.
The brace-compatible restrictions give the product and bracket
assertions as in Theorem~\ref{thm:relative-gerstenhaber}.
Over an $\mathbb F_p$-algebra regard the relative cochains as
a differential graded brace algebra over $\mathbb F_p$.
The restricted operation on shifted-even classes is represented by
the $p$-fold iterated unary brace
\cite[Theorem~2 and Appendix~A]{BriggsRubio2022}.
Corner restriction preserves this operation. An input of
degree $n>0$ has output degree $p(n-1)+1>0$. Applying the
operation to a triangular lift and using the graph in the
output degree proves
$\phi_{p(n-1)+1}(\alpha^{[p]})=\phi_n(\alpha)^{[p]}$.
\end{proof}

For actual identities this is the compression map of
Theorem~\ref{thm:relative-gerstenhaber}. In particular, the
restricted-power result applies to all classical stable
equivalences of Morita type in positive characteristic. In
degree one the operation is the $p$-fold composition of derivations.

\begin{remark}
Hermann's Oberwolfach report \cite[Section~4, pp.~1363--1364]
{Hermann2015Oberwolfach} cites the stable Hochschild ring result and
asks about invariance of the Lie structure in appropriate degrees.
The proposed approach uses a bilax functor and an Ext-vanishing
condition on the complementary summand of the image of the unit.
This question is distinct from the one answered by Volkov and
Witherspoon \cite{VolkovWitherspoon2024}.  They address
Question~5.2.15 of Hermann's monograph: whether the bracket defined
by loops in extension categories gives a Gerstenhaber algebra
structure.  Their result applies when the monoidal unit has a
projective power-flat resolution.
Theorem~\ref{thm:relative-gerstenhaber} answers the positive-degree
version of the 2015 invariance question under the displayed
bimodule hypotheses.  It removes the self-injectivity and
Gorenstein assumptions for the ordinary positive-degree assertion.
The generalized singular or dg equivalences considered in
\cite{BriggsRubio2022,Wang2021Gerstenhaber} are different inputs
from the bimodule identities used here.

Our proof obtains bracket compatibility from the common-image
identity \eqref{eq:gerstenhaber-common-image} and corner restriction
of cochains.  It does not construct a bilax structure on the given
tensor functor.  The required vanishing holds for the off-diagonal
component of a class in the image of $F$; the group
$\operatorname{Ext}_{B^e,\mathrm{rel}}^{>0}(B,Q)$ itself need not vanish.
\end{remark}

\begin{remark}
The full stable Hochschild ring, including its modified degree zero,
need not inherit the Gerstenhaber bracket from ordinary Hochschild
cohomology.  Proposition~\ref{prop:projective-center-not-gerstenhaber-ideal}
gives an explicit algebra, free of rank seven over its base, for
which the projective center is not a Gerstenhaber Lie ideal.  Thus
the positive-degree range above is essential.  For symmetric
algebras over a field, a projective-center Lie-ideal statement is available in
\cite[Lemma~10.6]{KoenigLiuZhou2012}.

The graph argument gives an isomorphism on positive cohomology.
It does not give a brace or $B_\infty$ quasi-isomorphism between
the Hochschild cochain complexes: the triangular restriction maps
may have nonzero kernels. Section~\ref{sec:tate-bv} compares
the Batalin--Vilkovisky operator under compatibility of the
Frobenius functionals.
\end{remark}

\subsection{Relative singular Hochschild cochains}
\label{subsec:relative-singular-braces}

We use the right singular model of Chen, Li and Wang
\cite[Sections~8.1 and~9.2--9.4]{ChenLiWangBInfinity}. For a
central algebra $D$ over a commutative ring $R$, write
$\mathsf S_R(D)$ for the bounded derived category of all
$D$-modules with the $R$-split exact structure, modulo the
thick subcategory generated by its relative projectives.
Write $\Omega_R^j$ for syzygies in this exact structure.

\begin{theorem}\label{thm:relative-singular-braces}
Let $A,B$ be central $R$-algebras and let $M,N$ be central
$A$-$B$ and $B$-$A$ bimodules, respectively, finitely generated
projective on each side. Suppose that
\[
 M\otimes_BN\simeq\Omega_R^\ell A,\qquad
 N\otimes_AM\simeq\Omega_R^\ell B
\]
in the respective additive quotients by $R$-relative projective
bimodules, for some $\ell\geq0$. Put
$T=\left(\begin{smallmatrix} A&M\\0&B\end{smallmatrix}\right)$.
Then the two diagonal restrictions
\[
 C^\bullet_{\mathrm{sg},\mathrm{right},R}(A)
 \ \xleftarrow{\ r_A\ }\
 C^\bullet_{\mathrm{sg},\mathrm{right},R\times R}(T)
 \ \xrightarrow{\ r_B\ }\
 C^\bullet_{\mathrm{sg},\mathrm{right},R}(B)
\]
are strict $B_\infty$ quasi-isomorphisms. On cohomology the map
from $B$ to $A$ takes $b$ to the unique $a$ such that
$a\otimes_A M=M\otimes_B b$ in
$\operatorname{Hom}_{\mathsf S_R(A\otimes_R B^{\mathrm{op}})}
(M,M[d])$.
\end{theorem}

\begin{proof}
Put $E=R\times R$ and use its two vertex idempotents acting
on $T$. For a diagonal algebra $D$ put $\bar D=D/(R1_D)$.
Use normalized bar cochains with values in
$D_p=\bar D^{\otimes_R p}\otimes_R D$ in degree $-p$,
with the left action given by the bar differential. The right
singular complex is their colimit under the transition
that prepends a suspended bar entry. For $T$ use
$T/\operatorname{im}(E\to T)$ and tensor products over $E$.
Insertion of the unit contracts the augmented normalized bars
over their base. The terms for $T$ decompose into induced
relative projectives $Te_i\otimes_R U\otimes_R e_jT$.
The face and insertion formulas of
\cite[Section~8.1]{ChenLiWangBInfinity} use tensor products,
unit maps and Koszul signs, and give the bar differential,
transitions and brace operations over $R$. On a pure
corner all intermediate endpoints are at that corner, so both
restrictions commute with each operation and with the transition.

Let $C=A\otimes_R B^{\mathrm{op}}$. Use the total bar
resolution $\operatorname{Bar}_R(A)\otimes_A M
\otimes_B\operatorname{Bar}_R(B)$ of $M$ for the mixed input
words. Use also the resolution formed by tensoring the two
\emph{augmented} bars with $M$ in between, with the term $M$
as augmentation, as in \cite[Section~9.3]{ChenLiWangBInfinity}.
Its $p$th term, with suspensions omitted, is
\[
 (A\otimes_R\bar A^{\otimes p}\otimes_R M)
 \oplus (M\otimes_R\bar B^{\otimes p}\otimes_R B)
 \oplus\bigoplus_{i+j=p-1}
 A\otimes_R\bar A^{\otimes i}\otimes_R M
     \otimes_R\bar B^{\otimes j}\otimes_R B.
\]
Both resolutions are $R$-split, using the bar contractions
linear on the sides tensored with $M$. Their terms are
relative projectives, by the induced form of the mixed
summands and the one-sided projectivity of $M$ in the two
outer summands. The $p$th syzygy of the second resolution
has underlying module
\[
 M_p=(\bar A^{\otimes p}\otimes_R M)
 \oplus\bigoplus_{i+j=p-1}
       \bar A^{\otimes i}\otimes_R M\otimes_R
                          \bar B^{\otimes j}\otimes_R B.
\]
Indeed, in the cokernel of the next differential the bar
relations move a leftmost coefficient to a representative
beginning with the left unit. The mixed summand gives the
same reduction for a word beginning at $M$. Applying the
differential and projecting onto the displayed module gives
the inverse to these unit representatives. These are the
$R$-linear formulas of
\cite[Lemma~9.8 and eqs.~(9.9)--(9.10)]{ChenLiWangBInfinity}.
The pure left and right maps into this resolution give
$A_p\otimes_A M\to M_p$ and $M\otimes_B B_p\to M_p$.
The former is pure left inclusion. The latter uses the right
bar action, namely the alternating sum obtained by merging
adjacent entries of $M\otimes_R\bar B^{\otimes p}\otimes_R B$,
including its two outer faces.

A composable triangular input word is pure or contains
exactly one entry from $M$. Thus at stage $p$ restriction onto
the two pure complexes is a degreewise split epimorphism
with kernel the mixed complex shifted by $[-1]$. The mixed
term computes relative Ext from $M$ to $M_p$ in degree $-p$,
using the first resolution. To compute the connecting map,
extend pure cochains by zero. In a mixed word, the only
surviving faces for these extensions are the outer actions
on a word ending or beginning in $M$. They use the two
syzygy maps above with opposite signs. The transition on
the kernel prepends a suspended left bar factor and is zero
on a word beginning in $M$. This induces the connecting
map of the mixed syzygy sequence. Indeed, lift a cochain
with values in $M_p$ by the left-unit representatives,
apply the coboundary and project to $M_{p+1}$. The internal
faces cancel in pairs and the left face gives precisely the
prepended factor with the suspension sign. The same check
applies to the pure terms. This is the calculation of
\cite[eqs.~(9.20)--(9.22) and Lemma~9.12]{ChenLiWangBInfinity};
all splittings and calculations above take place over $R$.

In the $R$-split exact structure, for modules $X,Y$ one has
\[
 \operatorname{Hom}_{\mathsf S_R(D)}(X,Y[d])
  \simeq\varinjlim_p
       \operatorname{Ext}^{d+p}_{D,\mathrm{rel}}(X,\Omega_R^pY)
\]
with connecting transitions. Indeed, $Y\to\Omega_R^pY[p]$
has perfect cone. These maps are cofinal among the right
denominators out of $Y$, since Hom from a bounded complex of
relative projectives, or its shift or retract, to
$\Omega_R^pY[p]$ vanishes for $p$ large. This vanishing gives
both factorization across any denominator and eventual
equality of factorizations in the calculus of fractions.
Taking the colimit, the connecting map in degree $d$ is,
up to the sign of the kernel identification, the difference
of the two tensor characteristic maps to
$\operatorname{Hom}_{\mathsf S_R(C)}(M,M[d])$.

For the first characteristic map use right tensor with $M$
and then with $N$. These functors preserve relative conflations
and relative projectives. The two orders of composition give
the shift $[-\ell]$ on the Hom groups between $A$ and its
shifts and between $M$ and its shifts, respectively. To see
this, stabilize each syzygy identity by adding relative
projectives as in the proof of
Theorem~\ref{cor:gerstenhaber-stable-identities}. A relative
projective error has the form of a summand of
$D\otimes_R U\otimes_R D$. Tensoring it on the left with
an object projective over the left endpoint makes a
relative projective bimodule. The bar syzygy sequences
split on their left side, and after tensoring their middle
terms on these objects are also relative projective.
The induced shift identifications are natural on these
Hom groups, since all tensors and connecting maps descend
to the quotient. For the syzygies this follows by their
one-sided split presentations as summands of induced bar
terms and the $R$-split tensor argument. Thus both composites
on Hom are bijective, and so is the first characteristic
map. For the second use left tensor and the right endpoint
instead.

Filtered colimits are exact. Since both characteristic
maps are bijective, the cohomology of the triangular
complex is their graph, and restriction to either corner
is an isomorphism. This proves the theorem.
\end{proof}

\begin{remark}
Over a field the theorem is the comparison in
\cite[Theorem~9.6 and Proposition~9.7]{ChenLiWangBInfinity}.
The relative comparison uses the split exact structure and
gives its two strict corner maps directly on the bar models.
\end{remark}

\subsection{The bracket and the projective center}
\label{sec:stable-center-bracket-counterexample}

Replacing degree zero by the stable center need not give a quotient
Gerstenhaber algebra.  The following example is free over an arbitrary
nonzero commutative base ring.

\begin{proposition}
\label{prop:projective-center-not-gerstenhaber-ideal}
Let $R$ be a nonzero commutative ring.  There exists an $R$-algebra $A$,
free of rank seven as an $R$-module, together with
\[
 c\in Z^{\mathrm{pr}}(A)
 \qquad\text{and}\qquad
 \alpha\in\mathrm{HH}^{2}_R(A),
\]
such that
\[
 [\alpha,c]\ne 0
 \quad\text{in }\mathrm{HH}^{1}_R(A).
\]
Here $A^e_R=A\otimes_R A^{\mathrm{op}}$,
$Z^{\mathrm{pr}}(A)$ consists of the bimodule endomorphisms of $A$ that
factor through projective $A^e_R$-modules, and
$\mathrm{HH}^{*}_R(A)$ is Hochschild cohomology relative to
$R$-split exact sequences.  Since $A$ is free over $R$, this cohomology agrees with
$\operatorname{Ext}^{*}_{A^e_R}(A,A)$.

Consequently, the graded ring
\[
 \mathrm{HH}^{0}_{\mathrm{st},R}(A)
     =Z(A)/Z^{\mathrm{pr}}(A),
 \qquad
 \mathrm{HH}^{n}_{\mathrm{st},R}(A)
     =\mathrm{HH}^{n}_R(A)\quad(n>0)
\]
need not inherit the Gerstenhaber bracket from
$\mathrm{HH}^{*}_R(A)$.
\end{proposition}

\begin{proof}
Let $Q$ have vertices $0,1$, a loop $x$ at vertex $0$, and arrows
\[
 a:0\longrightarrow 1,
 \qquad
 b:1\longrightarrow 0.
\]
We multiply paths from left to right.  Equivalently, our Peirce
convention is
\[
 e_0x=x=xe_0,
 \qquad
 e_0a=a=ae_1,
 \qquad
 e_1b=b=be_0.
\]
Define
\[
 A=RQ/(x^3,xa,bx,ba),
 \qquad c=ab.
\]
The residue classes of the surviving paths form the $R$-basis
\[
 e_0,\ e_1,\ x,\ a,\ b,\ x^2,\ c.
\]
Let $E=Re_0\oplus Re_1\cong R\times R$, and let $J$ be the arrow
ideal, with basis $x,a,b,x^2,c$.  Thus $A=E\oplus J$ as
$E$-bimodules.  The only nonzero products of two elements from this
basis of $J$ are
\[
 x\cdot x=x^2,
 \qquad a\cdot b=c.
\]
In particular, $J^3=0$.  Both $x$ and $c$ belong to $Z(A)$.  Since $R$ is nonzero and $c$
is a basis element, $c\ne0$.

We first show that $c$ lies in the projective center.
Give $A\otimes_R A$ its outer bimodule structure,
\[
 r\cdot(u\otimes v)\cdot s=ru\otimes vs.
\]
This bimodule is free of rank one over $A^e_R$.  We have
\[
 Ja=0=bJ,
 \qquad e_0a=a,
 \qquad be_0=b,
 \qquad e_1a=0=be_1.
\]
It follows that
\[
 ra\otimes b=a\otimes br
 \qquad(r\in A).
\]
Hence
\[
 \iota:A\longrightarrow A\otimes_R A,
 \qquad \iota(r)=ra\otimes b,
\]
is a bimodule homomorphism.  The multiplication map
$\mu:A\otimes_R A\to A$ is also a bimodule homomorphism, and
\[
 (\mu\iota)(r)=rab=rc.
\]
Thus multiplication by $c$ factors through a free $A^e_R$-module, so
$c\in Z^{\mathrm{pr}}(A)$.

We next construct a Hochschild two-cocycle.  The algebra
$E\cong R\times R$ is separable over $R$, with
separability element $e_0\otimes e_0+e_1\otimes e_1$.  The normalized
$E$-relative Hochschild cochains are
\[
 C_E^n(A)=\operatorname{Hom}_{E\text{-}E}
       \bigl(J^{\otimes_E n},A\bigr).
\]
They compute $\mathrm{HH}^{*}_R(A)$.  Indeed, the normalized
$E$-relative bar resolution has terms
\[
 A\otimes_E J^{\otimes_E n}\otimes_E A,
\]
and the usual comparison with the $R$-relative bar resolution is a
homotopy equivalence of relative projective resolutions.  In this
example, each Peirce component of every $J^{\otimes_E n}$ is free
over $R$.  Each displayed term is consequently a finite direct sum
of bimodules $Ae_i\otimes_R e_jA$, which are projective over
$A^e_R$.  Thus the same resolution also computes ordinary
$\operatorname{Ext}^{*}_{A^e_R}(A,A)$.  This also follows from the $R$-freeness of $A$.

Define an $E$-bimodule homomorphism
\[
 f:J\otimes_E J\longrightarrow A
\]
by
\[
 f(x,a)=a,
 \qquad f(x,c)=c,
\]
and set its value to zero on all other tensors of basis elements
of $J$.  These values respect the Peirce components and therefore
define an $E$-relative cochain.
If $\pi:A\to J$ is the projection along $E$, the corresponding
ordinary $R$-bilinear Hochschild cochain is
\[
 \widetilde f(u,v)=f\bigl(\pi(u),\pi(v)\bigr).
\]
We use $f$ for this extension as well.

For $u,v,w\in J$, the Hochschild differential is
\[
 (df)(u,v,w)
 =u f(v,w)-f(uv,w)+f(u,vw)-f(u,v)w.
\]
The first term vanishes because every nonzero value of $f$ is a
multiple of $a$ or $c$, and $Ja=Jc=0$.  The second term vanishes
because $uv$ is a linear combination of $x^2$ and $c$, whereas the
first argument of a nonzero value of $f$ must be $x$.  On triples of
basis elements, the third term is nonzero only at $(x,a,b)$, where it
equals $f(x,c)=c$.  The fourth term is nonzero only at that same
triple, where it equals $-f(x,a)b=-c$.  These two values cancel.
Therefore $df=0$ in the normalized $E$-relative complex.  Equivalently,
the extension of $f$ by $\pi$ is a cocycle in the ordinary Hochschild
cochain complex: triples involving $E$ satisfy the differential
identity by $E$-bilinearity and balancing.  Put
\[
 \alpha=[f]\in\mathrm{HH}^{2}_R(A).
\]

Finally, we compute the bracket.  Regard $c$ as a Hochschild
zero-cochain.  With the usual Gerstenhaber
sign convention, insertion of a zero-cochain into a two-cochain gives
\[
 [f,c](u)=f(c,u)-f(u,c).
\]
Consequently, the one-cochain $D=[f,c]$ satisfies
\[
 D(x)=-c,
 \qquad
 D(e_0)=D(e_1)=D(a)=D(b)=D(x^2)=D(c)=0.
\]
It is a derivation, either by the Hochschild bracket identity or by
direct verification using the multiplication table above.  Every
inner derivation annihilates the central element $x$.  Since
$D(x)=-c\ne0$, the derivation $D$ is not inner.  Thus
\[
 [\alpha,c]=[D]\ne0
 \quad\text{in }\mathrm{HH}^{1}_R(A).
\]

Finally, the canonical projection to the stable Hochschild graded
ring kills $c$ and is the identity in degree one.  If the Gerstenhaber
bracket descended along this projection, it would send
$[\alpha,c]$ to the bracket with the zero class of $c$, hence to
zero.  Its actual image is the nonzero class $[D]$.  Hence the bracket does not descend to the stable quotient.
\end{proof}

\begin{remark}
\label{rem:positive-range-essential}
Taking $R$ to be $\mathbb Z/4\mathbb Z$ gives a finite free Artin
$R$-algebra.  Taking $R=k[[t]]$ gives a finite free algebra over a
positive-dimensional complete Noetherian local ring.  Thus the
obstruction persists for algebras projective over either kind of
base.  Adjoining the ordinary stable center in degree zero does
not, in general, extend a Gerstenhaber invariance theorem in
positive degrees.
The quotient considered here changes only degree zero of ordinary
Hochschild cohomology and must be distinguished from singular or
Tate--Hochschild cohomology.  The all-degree invariance theorem of
\cite{Wang2021Gerstenhaber} concerns the latter theory; see
Section~\ref{sec:relative-gerstenhaber}.
\end{remark}

\section{The Tate--Hochschild BV operator for Frobenius algebras}
\label{sec:tate-bv}

Let $k$ be algebraically closed. For a finite-dimensional Frobenius
algebra $R$, fix a Frobenius functional $\lambda_R$ with
$\lambda_R(ab)=\lambda_R(b\nu_R(a))$ and put $\sigma_R=\nu_R^{-1}$.
We assume that $\nu_R$ is diagonalizable. Recall that the
BV operator of Itagaki, Sanada and Usui
\cite[Theorem~4.5 and Corollary~4.7]{ISUTateBV} is dual to
the twisted Connes operator in positive degrees. In degree
$-m-1$ it is given by the signed Connes operator on the
eigenvalue-one part of the right-twisted chains
$C_m(R,R_{\sigma_R})$. When $m=0$, its domain in cohomology
is the kernel of the norm from twisted $\mathrm{HH}_0$ to
the center.

Let $A,B$ have no semisimple direct factors and let $(M,N)$ be a
reduced pair defining a stable equivalence of Morita type. Put
$E=\operatorname{End}_{B^{\mathrm{op}}}(M)$ and let $\rho:A\to E$
be the left action. Choose a diagonalizable linear isomorphism
$\theta:M\to M$ with
\[
 \theta(amb)=\sigma_A(a)\theta(m)\sigma_B(b).
\]
For existence, the bimodule dualities in
\cite[Lemma~4.1(1)]{CPX2012} give
$D(A)\otimes_A M\simeq M\otimes_B D(B)$, where
$D=\operatorname{Hom}_k(-,k)$. Using $D(R)\simeq R_{\nu_R}$
and inverting the resulting semilinear map gives the required twist.
Its semisimple part has the same twist since multiplication by
eigenvectors sends generalized eigenspaces to the corresponding
shifted generalized eigenspaces.

Choose a right $B$-projective dual basis $(m_i,f_i)$ with
$\theta(m_i)=s_i m_i$ and $f_i\theta=s_i\sigma_B f_i$.
To obtain one, map a free module with basis of weights $s_i$
onto $M$ using eigenvectors and take the eigenvalue-one component
of a $B$-linear splitting under conjugation. Set
\[
 \lambda_E(e)=\sum_i s_i^{-1}\lambda_B(f_i(e m_i)).
\]
This is the restriction of the matrix functional
$\lambda_B\operatorname{tr}(D_s^{-1}-)$, with
$D_s=\operatorname{diag}(s_i)$, to the equivariant idempotent
corner defining $M$. The matrix functional is Frobenius and its
inverse Nakayama automorphism is entrywise $\sigma_B$ followed
by conjugation by $D_s$. The off-diagonal Peirce components
pair trivially with the corner, so the restricted pairing is
nondegenerate. The same construction on $B\oplus M$ restricts
to $\lambda_B$ and $\lambda_E$ in its two corners.

\begin{theorem}\label{thm:tate-bv-compatible}
In the preceding setting suppose
\begin{equation}\label{eq:tate-bv-compat}
 \lambda_A(a)=\lambda_E(\rho(a))
     =\sum_i s_i^{-1}\lambda_B(f_i(a m_i))\qquad(a\in A).
\end{equation}
Then the comparison $\Phi:\widehat{\mathrm{HH}}^*(B)
 \xrightarrow{\sim}\widehat{\mathrm{HH}}^*(A)$ induced by
$M\otimes_B-\otimes_BN$ and the stable regular summand is
an isomorphism of BV algebras in all degrees.
\end{theorem}

\begin{proof}
Define $\tau:E\to A$ by
$\lambda_A(a\tau(e))=\lambda_E(\rho(a)e)$.
Associativity and the twisted trace identity
$\lambda_R(uv)=\lambda_R(\sigma_R(v)u)$ show that $\tau$
is an $A$-bimodule map. Condition \eqref{eq:tate-bv-compat}
gives $\tau\rho=1_A$.

Put $S=A^e$, $T=E^e$, $L=T\otimes_S-$ and
$\pi=\tau\otimes\tau^{\mathrm{op}}$. Equip $R^e$ with
$\lambda_{R^e}=\lambda_R\otimes\lambda_R$, so
$\lambda_S\pi=\lambda_T$. The natural map
\[
 T\otimes_S X\longrightarrow\operatorname{Hom}_S(T,X),
 \qquad t\otimes x\longmapsto(t'\mapsto\pi(t't)x)
\]
is an isomorphism. For $X=S$, compose with $\lambda_S$ and use
nondegeneracy of $\lambda_T$. The general assertion follows
from finite projectivity of $T$ on the left over $S$.
Thus $L$ is both left and right adjoint to restriction.

Let $p:L(A)=E\otimes_A E\to E$ be multiplication and let
$j:E\to L(A)$ correspond to $\tau$ under the right adjunction.
The map $\rho$ is the unit of $-\otimes_A M\dashv-\otimes_B N$
on the regular bimodule, using $N=\operatorname{Hom}_{B^{\mathrm{op}}}(M,B)$.
It is a stable isomorphism. Under the Morita equivalence between
$E$ and $B$, $p$ corresponds to the counit $N\otimes_A M\to B$
of the same adjunction with left endpoint $B$, so it is also a
stable isomorphism. These adjunctions descend to equivalences
on the stable categories of bimodules projective on each side.
Here $M_B$ is a progenerator since $B$ is a summand of
$N\otimes_A M$ and $N_A$ is finite projective.

Let $u:A\to\operatorname{Res}L(A)$ be the induction unit and
$v:\operatorname{Res}L(A)\to A$ the other counit, given by
$\pi\otimes1$. Then $(\operatorname{Res}p)u=\rho$,
$vu=1_A$ and $v\operatorname{Res}j=\tau$.
Since $u$ is a stable isomorphism, $\operatorname{Res}(pj)$
is the identity stably. Restriction detects this equality:
precomposition with $p$ is bijective on stable maps from $E$
and, under adjunction, is restriction followed by precomposition
with $\rho$. Thus $j=p^{-1}$ stably. From $E$ to $A$ the
comparison is restriction followed by compression with $\rho$
and $\tau$. After the Morita step this is the given bimodule
comparison, since the Morita functor from $B$ to $E$ is
$M\otimes_B-\otimes_B N$ before restriction. Its inverse
from $A$ to $E$ is $[g]\mapsto p[n]L([g])j$ in degree $n$.

For a left $S$-module $X$, write $X^\vee=\operatorname{Hom}_S(X,S)$.
The dual of a termwise finite projective right resolution
$Q_\bullet$ of $X^\vee$ gives the negative part of a complete
resolution of $X$, spliced to a projective resolution by evaluation.
There are natural identifications
\[
 T\otimes_S Q_i^\vee\simeq(Q_i\otimes_S T)^\vee,\qquad
 (L(X))^\vee\simeq X^\vee\otimes_S T.
\]
The first is finite projective duality and the second follows from
adjunction and left $S$-projectivity of $T$. Under these maps,
induction in negative Tate Ext is change of rings on Tor.
Precomposition with $j$ and postcomposition with $p$ act by
$j^\vee$ and $p$ on the two Tor variables. Indeed, lift the
coefficient maps to projective resolutions and apply finite
projective duality termwise. The identifications commute with
evaluation at the splice, giving the comparison of evaluation
kernels in degree $-1$ \cite[Proposition~2.3]{ISUTateBV}.

Write $\gamma_x^R$ for the image of $x$ in
$R_{\sigma_R}\simeq\operatorname{Hom}_{R^e}(R,R^e)$.
The formula in \cite[Lemma~2.1]{ISUTateBV} gives
$\lambda_{R^e}(\gamma_x^R(y))=\lambda_R(xy)$.
By adjunction,
$\pi j^\vee(\gamma_x^A\otimes1)=\gamma_x^A\tau$.
On the other hand,
\[
 \lambda_S\bigl(\pi\gamma_{\rho(x)}^E(y)\bigr)
 =\lambda_E(\rho(x)y)
 =\lambda_A(x\tau(y))
 =\lambda_S\bigl(\gamma_x^A(\tau(y))\bigr).
\]
Both maps to $S$ are $S$-linear. Multiplying the input by
elements of $S$ and using nondegeneracy proves equality of the maps.
The adjunction therefore gives
$j^\vee(\gamma_x^A\otimes1)=\gamma_{\rho(x)}^E$.
Also $p(1\otimes a\otimes1)=\rho(a)$. Consequently the negative
comparison from $A$ to $E$ is induced on twisted Hochschild
chains by applying $\rho$ to every factor, including the
coefficient. For degree $-1$ this map preserves the norm kernel
by the evaluation comparison above. Explicitly, the two
maps on the domain and codomain of the norm are induced by
$j^\vee\otimes p$ and $h\mapsto pL(h)j$. The identity
$(j^\vee L(\gamma))(y)=L(\gamma)(j(y))$ makes the square with
the evaluation maps commute before passing to the kernel.

For the Morita step, let $e$ be a $\sigma_D$-fixed full idempotent
of a Frobenius algebra $D$, put $H=eDe$ and
$f=e\otimes e^{\mathrm{op}}$. The Morita functor on enveloping
modules is $D^e f\otimes_{H^e}-$. Multiplication identifies the
induction of $H$ with $D$, and restriction gives
$D^\vee f\simeq H^\vee$. For $x\in H$, $\gamma_x^D$ is
supported on $f$, and its restriction to $H$ takes values in
$fD^e f$. Comparing the restricted Frobenius functionals shows
that this restriction is $\gamma_x^H$. Finite projective duality
on the induced resolutions therefore gives the corner inclusions
on both Tor variables, including the map on the norm kernel.
Apply this calculation to the two full corners of
$D=\operatorname{End}_{B^{\mathrm{op}}}(B\oplus M)$.

The map $\rho$ and the two corner inclusions, viewed as functors
into the category with objects $B$ and $M$, preserve identities
and twists. Thus they preserve faces, degeneracies and cyclic
rotation in the paracyclic model. On the eigenvalue-one part
they commute with Connes. The corner maps induce Morita
isomorphisms on twisted Hochschild homology, also on the
eigenvalue-one part by diagonalizability. This proves
commutation with the BV operator for $n\le-1$.

For positive degrees use left-twisted chains with coefficients
in $ {}_{\sigma_R}R$ paired with cochains by
$\lambda_R(a_0 f(a_1,\ldots,a_n))$. The same paracyclic diagram
gives a transpose comparison commuting with the BV operator.
To identify it, let $\beta_E$ correspond to $\beta\in\mathrm{HH}^n(B)$
under Morita equivalence and put $\alpha=\Phi(\beta)$. For $n>0$
the equality of characteristic actions in
Lemma~\ref{lem:gerstenhaber-common-image} says
$\rho_*(\alpha)=\rho^*(\beta_E)$ in $\mathrm{HH}^n(A,E)$.
Pair this identity with a twisted chain and use
$\lambda_A=\lambda_E\rho$. The Morita step respects the pairing
by evaluation on corner chains and restriction of the functional.
Thus the transpose is $\Phi$ in positive degrees.
In degree zero it is compression by $\tau$ after the Morita step,
which induces the stable center comparison. This proves
commutation from degree one to degree zero. The operator is
zero in degree zero. The comparison preserves the product by
\cite[Theorem~6.2]{Wang2021Gerstenhaber}.
\end{proof}

\begin{remark}
In the left-twisted degree-zero chain model, write
$e=\sum_i e\,m_i f_i$. Moving the morphism $M\to B$ given
by $f_i$ across the twisted cyclic rotation gives the term
$s_i^{-1}f_i(e m_i)$ in the $B$ corner. Thus $\lambda_E$
is the pullback of $\lambda_B$ through the twisted Morita
comparison and is independent of the equivariant dual basis.
The existence of compatible functionals is established in
Corollary~\ref{cor:tate-bv-normalization}.
\end{remark}

The choice of functional matters even for the identity equivalence.
For instance, suppose $\operatorname{char}k\ne2$ and let $C$ be the
radical-square-zero algebra of the oriented two-cycle.
Let $\lambda_C$ take value $1$ on each arrow and $0$ on each
vertex idempotent. On $A=C\otimes_k k[x]/(x^2)$ put
$\lambda(c_0+c_1x)=\lambda_C(c_1)$ and
$\lambda'(a)=\lambda((1+x)a)$. The Nakayama automorphism of both
functionals exchanges the two vertices and the two arrows of $C$
and fixes $x$. For the invariant derivation $d=x\partial_x$,
the degree-one operator gives $\Delta_\lambda[d]=1$ and
$\Delta_{\lambda'}[d]=1-x$, by evaluating
$\lambda((\Delta_\lambda d)a)=\lambda(d(a))$ and the analogous
formula for $\lambda'$. The two values already differ in
$\mathrm{HH}^0(A)$. They differ in the stable center as well:
multiplication by $x$ on $S\otimes_k k[x]/(x^2)$ is nonzero for
a simple $C$-module $S$. After restriction to $C$ this module is
$S^{\oplus2}$, whose endomorphisms factoring through projectives
are zero. An element of the projective center acts on every
module by a map factoring through a projective, and restriction
from $A$ to $C$ preserves projectives.

\begin{corollary}
\label{cor:tate-bv-normalization}
Keep the algebra and bimodule hypotheses of
Theorem~\ref{thm:tate-bv-compatible}. Fix the functional
$\lambda_B$, an initial functional $\lambda_A^0$ defining the
prescribed $\nu_A$, and any diagonalizable semilinear
isomorphism $\theta$ as above. Then $\lambda_E\circ\rho$ is a
Frobenius functional on $A$ with Nakayama automorphism $\nu_A$.
More precisely, there is a central unit $z\in Z(A)^\times$ such
that
\[
 \lambda_E(\rho(a))=\lambda_A^0(az)\qquad(a\in A).
\]
Consequently, for a compatible choice of Frobenius
functionals the given stable equivalence induces an isomorphism
of Tate--Hochschild BV algebras in all degrees.
\end{corollary}

\begin{proof}
Define $\tau:E\to A$ by
$\lambda_A^0(a\tau(e))=\lambda_E(\rho(a)e)$.
The twisted trace identities and $\rho\sigma_A=\sigma_E\rho$
show that $\tau$ is an $A$-bimodule map. This step works
before imposing \eqref{eq:tate-bv-compat}.

Write $M^*=\operatorname{Hom}_{B^{\mathrm{op}}}(M,B)$ and
identify $E$ with $M\otimes_B M^*$ by the evaluation of rank-one
endomorphisms. For $m\in M$ and $f\in M^*$, the equivariant
dual basis gives
\begin{align*}
 \lambda_A^0(\tau(m\otimes f))
   &=\sum_i s_i^{-1}\lambda_B\bigl(f_i(m)f(m_i)\bigr)\\
   &=\sum_i s_i^{-1}\lambda_B\bigl(f(m_i)\nu_B(f_i(m))\bigr)\\
   &=\lambda_B\bigl(f(\theta^{-1}m)\bigr).
\end{align*}
The last equality follows from
$\theta^{-1}m=\sum_i s_i^{-1}m_i\nu_B(f_i(m))$.
Since $B$ is Frobenius, the map $M^*\to D(M)$ given by
$f\mapsto\lambda_B\circ f$ is a vector space isomorphism.
Similarly $\operatorname{Hom}_A(M,A)\to D(M)$ given by
composition with $\lambda_A^0$ is an isomorphism.
The preceding identity and invertibility of $\theta$ show that
the $B$--$A$-bimodule map
\[
 M^*\longrightarrow\operatorname{Hom}_A(M,A),
 \qquad f\longmapsto\bigl(m\mapsto\tau(m\otimes f)\bigr)
\]
is an isomorphism. Hence $\tau$ is the evaluation counit for
$-\otimes_B M^*\dashv-\otimes_A M$. The map $\rho$
is the unit of $-\otimes_A M\dashv-\otimes_B M^*$.

With left endpoint $A$, the two tensor functors induce
inverse equivalences on the stable categories of bimodules
projective on each side, by \eqref{eq:MT}. They preserve
projective bimodules, so both adjunctions descend. Thus both
$\tau$ and $\rho$ are isomorphisms in the stable category.
Their composite is multiplication by $z=\tau(1_E)\in Z(A)$
and is stably invertible. The regular bimodule $A$ has no
projective direct summand, since such a summand would be a
separable direct algebra factor. By
Lemma~\ref{lem:projective-factorization-radical}, $z$ is an
actual unit. It follows that $a\mapsto\lambda_E(\rho(a))=\lambda_A^0(az)$
is nondegenerate. Since $z$ is central, this functional has the
same Nakayama automorphism as $\lambda_A^0$.
Taking $\lambda_A=\lambda_E\rho$ now gives
\eqref{eq:tate-bv-compat}, and
Theorem~\ref{thm:tate-bv-compatible} applies.
\end{proof}

\subsection{Comparison with the Tate transfer}
\label{subsec:singular-braces-bv}

\begin{corollary}\label{cor:singular-braces-bv}
Under the hypotheses of Theorem~\ref{thm:tate-bv-compatible},
the isomorphism in $\operatorname{Ho}(B_\infty)$ of
Theorem~\ref{thm:relative-singular-braces}, in the direction from
$B$ to $A$, induces $\Phi$ under the identifications of singular
Hochschild cohomology with Tate--Hochschild cohomology. In
particular, its cohomology map commutes with the Tate BV
operators in every degree, including the norm kernel in degree
$-1$.
\end{corollary}

\begin{proof}
The enveloping algebras are Frobenius. The stable category of
each enveloping algebra identifies with its singularity category,
and the stable shift identifies the Hom groups of the regular
bimodule with Tate Ext. These identifications are induced by
passing from the stable module to the singular quotient and so
commute with the exact tensor functors on bimodules. Let $b$ be
a degree $d$ class for $B$ and let $a$ be its image under the
singular corner comparison. Its defining equation is
$a\otimes_A M=M\otimes_B b$. Tensor it on the right by $N$.
Put $V=M\otimes_B N$ and let $i:A\to V$ and $p:V\to A$
be the stable regular inclusion and projection used to define
$\Phi$. Naturality of tensoring in the coefficient bimodule
gives
\[
 a=p[d](a\otimes_A V)i
   =p[d](M\otimes_B b\otimes_B N)i=\Phi(b).
\]
For the first equality one may use actual maps representing
$i,p$ and the natural transformations induced by them on the
relative derived categories and then on the quotients. This
works for every $d$. The negative chain comparison and the
evaluation square at the splice in
Theorem~\ref{thm:tate-bv-compatible} now give the stated BV
compatibility for this very map. The assertion compares the
cohomology operators.
\end{proof}

\section{Deformations of algebras and modules}
\label{sec:associative-deformations}

In this section the base is a commutative ring $R$, and the
flat $R$-algebras need not be finite over $R$ unless indicated otherwise.
For an $R$-flat algebra $D$ and a surjection $S\to R$ of commutative rings with nilpotent
kernel, let $\operatorname{Def}_{D/R}(S)$ be the set of isomorphism
classes of $S$-flat associative unital $S$-algebras $\mathcal D$
with a marking $\mathcal D\otimes_S R\simeq D$. Isomorphisms
are required to preserve the marking. Morphisms of coefficient
rings are maps over $R$, and base change defines
$\operatorname{Def}_{D/R}$. When $R=k$ is a field, we also write
$\operatorname{Def}_D$ on local Artin $k$-algebras with residue field $k$.

Let $I$ be an ideal of $S$ with $S/I=R$. A compatible system of
deformations over $S/I^n$ includes specified identification maps
between successive reductions. If $S$ is $I$-adically complete,
a formal deformation is an $I$-adically complete $S$-algebra
whose reductions are flat over $S/I^n$, with a marking modulo
$I$. For formal deformations, base change is completed base change.

\subsection{Stable equivalences of Morita type}

Let $A,B$ be $R$-flat algebras and let $M$ be an $R$-flat
$A$-$B$-bimodule. Set $C=A\otimes_R B^{\mathrm{op}}$ and
$U=\left(\begin{smallmatrix}A&M\\0&B\end{smallmatrix}\right)$.
For an $R$-module $J$, write
\[
 \mathrm{HH}^{n}(D/R;J\otimes_R D)
   =\operatorname{Ext}_{D^e}^{n}(D,J\otimes_R D),
 \qquad E_J^n=\operatorname{Ext}_{C}^{n}(M,J\otimes_R M).
\]
Let $a_J^n,b_J^n$ denote the characteristic maps from the
cohomologies of $A$ and $B$, respectively, to $E_J^n$, given
by derived tensoring with $M$. The coefficient objects
$(J\otimes_R A)\otimes_A^{\mathbf L}M$ and
$M\otimes_B^{\mathbf L}(J\otimes_R B)$ both identify with
$J\otimes_R M$ by flatness over $R$. Denote
by $r_A^n,r_B^n$ the two corner restrictions from the cohomology
of $U$ with coefficients in $J\otimes_R U$.

\begin{theorem}\label{thm:stable-associative-deformations}
Suppose that, for every $R$-module $J$, both corner restrictions
are surjective in degree one, bijective in degree two and
injective in degree three. Then taking corners induces natural
bijections
\[
 \operatorname{Def}_{A/R}
  \xleftarrow{\ \sim\ }\operatorname{Def}_{U/R}
  \xrightarrow{\ \sim\ }\operatorname{Def}_{B/R}.
\]
The same comparison holds for compatible systems of marked
deformations over $S/I^n$. If $S$ is $I$-adically complete
and $I$ is finitely generated, it holds for formal deformations.

The restriction hypothesis holds if both $a_J^n,b_J^n$ are
surjective for $n=1$, bijective for $n=2$ and injective for
$n=3$, for every $J$. In this case the comparison on
$R\ltimes J$ is $(b_J^2)^{-1}a_J^2$. In particular, these
conditions hold if $M$ belongs to a biadjoint pair $(M,N)$
satisfying \eqref{eq:hh-mt}. For such a pair the comparison
on $R\ltimes J$ is the degree-two compression map with
coefficients in $J$.
\end{theorem}

\begin{proof}
We first record the obstruction calculation for a step
$S'\twoheadrightarrow S$ in a nilpotent extension of $R$,
where the kernel $J$ is annihilated by $\ker(S'\to R)$.
In the triangular case we work with lifts of the two vertex
idempotents. They exist across nilpotent ideals, and two
choices of complementary lifts are conjugate by a unit
congruent to $1$. The lower-left Peirce component vanishes
because it reduces to zero and the ideal is nilpotent.
Thus the choice of vertices does not change the set of
marked isomorphism classes.

Let $\mathcal D$ be an $S$-flat lift of $D$. Take a nonpositive
semifree dg algebra resolution $L\to\mathcal D$. In the
triangular case take it over $S\times S$, with generators in
the two diagonal components and the upper-right component.
Such a resolution is obtained by choosing degree-zero generators
in the three components and adjoining generators in negative
degrees to kill the successive kernels in cohomology. The two
corners are semifree resolutions of the corner algebras, since
a path leaving a vertex cannot return to it.

Lift the graded tensor algebra underlying $L$ to $S'$ and lift
its differential on the generators to a derivation $\widetilde d$.
The square $\widetilde d^{\,2}$ takes values in $J\otimes_S L$.
It is a degree-two cocycle in the complex
$\operatorname{Der}(L,J\otimes_S L)$, with derivations taken
over the coefficient ring and, in the triangular case, vanishing
on the vertices. The square of an odd derivation is a derivation
in every characteristic, since the two cross terms cancel in the
Leibniz formula. Also $[\widetilde d,\widetilde d^{\,2}]=0$,
and replacing $\widetilde d$ by $\widetilde d+h$ changes the
square by $[d,h]$, since $J^2=0$. Reducing the source to $R$ and
using flatness of $\mathcal D$, we compute this complex as
$\operatorname{Der}(L\otimes_S R,J\otimes_R D)$ up to
quasi-isomorphism. To see this, the bimodule of noncommutative
differentials of a semifree algebra is semifree, filtered by
the generators, and $J\otimes_S L\to J\otimes_R D$ is a
quasi-isomorphism. For $L_0=L\otimes_S R$ and $E=R$ or
$E=R\times R$, the sequence
\[
 0\longrightarrow\Omega^1_{L_0/E}
 \longrightarrow L_0\otimes_E L_0
 \longrightarrow L_0\longrightarrow0
\]
identifies cohomology in derivation degrees two and one with
$\mathrm{HH}^{3}(D/R;J\otimes_R D)$ and
$\mathrm{HH}^{2}(D/R;J\otimes_R D)$, respectively. Here $E$
is separable over $R$, and $R$-flatness of $D$ identifies the
derived enveloping algebra with $D^e$. In degree zero the closed
derivations of $L_0$ with these coefficients are the derivations
of $D$ into $J\otimes_R D$ vanishing on $E$. These
identifications commute with taking corners.

A vanishing obstruction allows one to correct $\widetilde d$
to a differential. The resulting free complex over $S'$ has
cohomology concentrated in degree zero and its degree-zero
cohomology is $S'$-flat. Indeed, for any $S'$-module $V$, the
filtration $0\subseteq JV\subseteq V$ has subquotients over
$S$. Tensoring with the free complex and using
$S$-flatness of $\mathcal D$ gives cohomology concentrated in
degree zero. This gives both assertions and the required
identification after reduction. Conversely, for any flat lift
one can lift the augmentation of $L$ to it together with the
differential, successively on the generators. At degree zero
one lifts the chosen images. For a generator of degree $-q$, with $q\geq1$, lift its old
differential to degree $1-q$ of the partially constructed
resolution. The defect of this lift as a cycle of the augmented
complex lies in degree $2-q$ of the reduction kernel, where
we place the augmentation target in degree one. This defect is
a cycle in the exact augmented complex
$J\otimes_S(L\to\mathcal D)$. A preimage in degree $1-q$
corrects the lift. These degrees involve only generators
already constructed. For two flat lifts, an isomorphism over
the fixed reduction likewise lifts to a chain map. Choose
the image of each generator to reduce to that generator of $L$.
At degree zero correct its augmentation, and in negative
degrees correct the defect of the chain-map equation using
the same exact reduction kernel of the target resolution. Two lifted differentials differ by a degree-one
derivation with values in $J\otimes_S L$. Conjugation by $1+h$,
with $h$ a degree-zero derivation with values there, changes this
difference by a coboundary. Hence the obstruction vanishes
precisely when a flat lift exists, and the lift classes with
fixed reduction then form a torsor under
$\mathrm{HH}^{2}(D/R;J\otimes_R D)$. Automorphisms inducing the
identity on the reduction have the form $1+\delta$, where
$\delta$ is a derivation of $D$ into $J\otimes_R D$ vanishing on $E$.

For either corner, injectivity in degree three gives existence
of a triangular lift when the corner lifts. Bijectivity in
degree two gives a unique lift class with the specified reduction
and corner class. The corner maps on derivations vanishing on
the vertices are surjective. Indeed, surjectivity in degree one
lifts derivation classes modulo inner derivations, and an
inner derivation of either corner lifts by using the same
diagonal element in $J\otimes_R U$. This also gives a relative
statement for isomorphisms. Given an isomorphism between two
triangular reductions and a compatible isomorphism of the
lifted corners, use the former to identify the reductions.
Injectivity in degree two gives a triangular isomorphism lifting
that identification. The difference from the specified corner
isomorphism is a derivation as above, and lifting it corrects
the triangular isomorphism. Induction along the powers of a
nilpotent kernel proves the bijections. Naturality follows
because taking a corner commutes with base change. The same
lifting argument produces compatible triangular lifts and
compatible comparison isomorphisms over $S/I^n$.
For finitely generated $I$, the limit has exactly the given
reductions by \cite[Tag \texttt{09B8}]{StacksProject}.

To verify the characteristic-map criterion, use the triangular
long exact sequence with coefficients in $J\otimes_R U$,
whose middle map is $a_J^n-b_J^n$
\cite[Section~4.5]{Keller2003Higher}. The sequence with these
coefficients follows by splitting the dg bar cochains into the
two corners and the mixed strings. In a triangular semifree
model each mixed string contains one off-diagonal entry.
The two-sided bar resolution of this entry gives the complex
computing $E_J^\bullet$, shifted by one in the mixed summand.
The connecting map is the difference of the two tensor actions,
as in Lemma~\ref{lem:relative-triangular-gerstenhaber}.
Equivalently, the homotopy fibre of the restriction to $A$ is
the shifted cone of the characteristic map from $B$, and
conversely. The conditions on $b_J^n$ and $a_J^n$ therefore
give the stated conditions on $r_A^n$ and $r_B^n$, respectively.
In degree two the image in the two corners is the graph of
$(b_J^2)^{-1}a_J^2$.

Finally, suppose $(M,N)$ is a biadjoint pair satisfying
\eqref{eq:hh-mt}. Use the adjunction
$-\otimes_B N\dashv-\otimes_A M$ on bimodules. Under the
resulting Ext adjunction, $a_J^n$ is pullback along the counit
\[
 \epsilon_A:M\otimes_B N\longrightarrow A
\]
in the first variable of $\operatorname{Ext}_{A^e}^n(-,J\otimes_R A)$.
On bimodules projective on both sides the two functors preserve
projective bimodules and are inverse equivalences after factoring
out projectives. This follows from \eqref{eq:hh-mt}, since
tensoring a projective error bimodule with such an object
again gives a projective bimodule. The adjunction descends to
the quotient, so $\epsilon_A$ is invertible there. A map in
the first variable factoring through a projective induces
zero in positive Ext. Hence $a_J^n$ is an isomorphism for
$n>0$, and the adjunction on the other side proves the same
for $b_J^n$. If $a_J^2(\alpha)=b_J^2(\beta)$, apply
$N\otimes_A-$ and take the regular $B$-component in
$N\otimes_A M=B\oplus Q$. On the left this is the compression
of $\alpha$, and on the right it is $\beta$, as in
Lemma~\ref{lem:gerstenhaber-common-image}.
\end{proof}

\begin{cor-nn}
A stable equivalence of Morita type between finite-dimensional
algebras over any field induces a natural isomorphism of the
functors of marked associative algebra deformations. Its
tangent map is the degree-two compression map.
\end{cor-nn}

\begin{proof}
A marked deformation of a product is the product of the
marked deformations of its factors. Indeed the block idempotents
lift, and the off-diagonal Peirce components vanish by
nilpotence. A separable factor has a unique marked lift class
at each order, since its positive Hochschild cohomology
with arbitrary coefficients vanishes. Apply Corollary~\ref{cor:adjoint-relative-stable-identities}
with the ordinary projectivity class. It gives a reduced pair
supported on the products of the nonseparable blocks. After
restriction to those products, its two composites contain the
regular bimodules, and its tensor functors are biadjoint.
Thus Theorem~\ref{thm:stable-associative-deformations} applies.
The reduced mixed bimodule, extended by zero on the separable
factors, gives the triangular comparison for the original
algebras. On a separable factor all positive cohomology with
the indicated coefficients vanishes, so the characteristic
map criterion on the full triangular algebra follows as well.
The comparison in positive cohomology is unchanged by omitting
projective summands from the pair or by this block reduction.
Indeed, the omitted tensor terms are projective, and the
separable factors have zero positive cohomology. Compression
therefore identifies the resulting comparison with that of
the original pair.
\end{proof}

For the rest of this subsection let $k$ be a field of characteristic
zero, and take the finite-dimensional triangular algebra $U$ of the
preceding corollary. Let $\mathcal O=k[[t_1,\ldots,t_d]]$ and let
$\mathfrak m$ be its
maximal ideal and let $V=(\mathfrak m/\mathfrak m^2)^*$.
For an adically complete flat deformation $\mathcal A$ with
special fibre $A$, Bezrukavnikov and Ginzburg associate a map
in the homotopy category of dg algebras
\[
 \operatorname{Deform}_{\mathcal A}:
 \operatorname{Sym}(V[-2])\longrightarrow
 \mathbf R\operatorname{Hom}_{A^e}(A,A).
\]
We use their conventions for the augmentation and the special
fibre identification \cite[Theorem~1.2.3]{BezrukavnikovGinzburg2007}.

\begin{proposition}\label{prop:bg-triangular-deformations}
Corresponding formal deformations $\mathcal A,\mathcal B$
in Theorem~\ref{thm:stable-associative-deformations} can be
realized as the corners of a common formal deformation
$\mathcal U$ of the triangular algebra $U$. For such a choice
the two Bezrukavnikov--Ginzburg maps are the composites of
$\operatorname{Deform}_{\mathcal U}$ with the two morphisms on
derived dg endomorphism algebras induced by taking corners.
\end{proposition}

\begin{proof}
The formal comparison in Theorem~\ref{thm:stable-associative-deformations}
gives $\mathcal U$ with the two prescribed marked corners.
The reductions are finite free over the Artin quotients of
$\mathcal O$, and hence $\mathcal U$ is finite free over
$\mathcal O$.

To compare the dg maps, use a common base resolution
$\mathcal R\to k$ and exterior dg algebra $\Lambda$ in the
construction of \cite[(3.1.4)--(3.1.6)]{BezrukavnikovGinzburg2007}.
For a deformation $\mathcal D$ the dg functor defining the
map is given by the kernel
\[
 \Theta_{\mathcal D}(W)
   =\big((\mathcal R\otimes_{\mathcal O}\mathcal R)
                           \otimes_{\mathcal O}\mathcal D\big)
                      \otimes_\Lambda W.
\]
If $e$ is either of the lifted vertices, taking $e(-)e$ on
$\Theta_{\mathcal U}(W)$ gives
$\Theta_{e\mathcal Ue}(W)$ as a natural dg identification.
Indeed $e$ acts on the algebra factor and commutes with the
base resolution and the exterior algebra action. Taking a
corner is exact and commutes with the augmentations to the
special fibre. In particular it carries the distinguished
object obtained from the trivial $\Lambda$-module to the
corresponding distinguished object for the corner. Passing to
derived dg endomorphisms gives
\[
\begin{split}
 \operatorname{Deform}_{\mathcal A}&=r_A\operatorname{Deform}_{\mathcal U},\\
 \operatorname{Deform}_{\mathcal B}&=r_B\operatorname{Deform}_{\mathcal U}.
\end{split}
\]
Here $r_A,r_B$ are induced by corner restriction on dg bimodule
categories. Their cohomology maps are the diagonal restrictions
used above.
\end{proof}

\subsection{Standard derived equivalences}
\label{subsec:standard-derived-deformations}

We keep the coefficient categories and the conventions for
formal deformations above. Let $A,B$ be $R$-flat algebras.

Fix a two-sided tilting complex $X$ of $B$-$A$-bimodules
with inverse $Y$. Such a pair can be chosen if the categories
of perfect complexes over $A$ and $B$ are $R$-linearly equivalent,
or if the unbounded derived categories of all modules are
$R$-linearly equivalent, by derived Morita theory
\cite{RickardDerived1989,anu-Rickard1991,Keller1994Deriving}.
For finite-dimensional algebras over a field one may start with
the bounded derived categories of finite-dimensional modules.
For the comparison of Bezrukavnikov--Ginzburg maps put
\[
 \Phi_X=X\otimes_A^{\mathbf L}(-)\otimes_A^{\mathbf L}Y.
\]
The identification $\Phi_X(A)\simeq B$ induces the comparison
$\Phi_{X,*}$ on derived dg endomorphism algebras of the diagonals.
Over a field of characteristic zero, Keller's Hochschild
comparison \cite[Section~4.6]{Keller2003Higher} also gives the
classical infinitesimal comparison by the dg Lie criterion.

\begin{proposition}\label{prop:standard-derived-deformations}
Let $R,A,B,X,Y$ be as above.
\begin{enumerate}
\item On surjections $S\to R$ with nilpotent kernel, $X$ induces a
natural isomorphism
\[
 \operatorname{Def}_{A/R}\simeq\operatorname{Def}_{B/R}.
\]
For corresponding deformations $\mathcal A,\mathcal B$, there is
a two-sided tilting complex over $\mathcal B$ and $\mathcal A$
whose special fibre is $X$.
\item For any ideal $I$ with $S/I=R$, the correspondence extends
to compatible systems of marked flat deformations over the
$S/I^n$. If $S$ is $I$-adically complete and $I$ is finitely
generated, taking inverse limits gives a natural bijection
of marked formal deformations. If both formal algebras are
$S$-flat, the two-sided tilting complex exists over $S$.
In particular, this holds if $S$ is Noetherian and complete.
\end{enumerate}
\end{proposition}

\begin{proof}
Choose a bounded complex $T$ of finitely generated projective right
$A$-modules representing $X_A$. Its endomorphism complex has
cohomology $B$ in degree zero, with the identification supplied
by the left action on $X$.

First let $I=\ker(S\to R)$ be nilpotent. Lift $T$ successively
through $S/I^j$. At a step with kernel $J=I^j/I^{j+1}$ the
projective terms lift by idempotent lifting. The square of a
graded lift of the differential is a degree-two cocycle in
$J\otimes_R\operatorname{End}_A^\bullet(T)$
\cite{Rickard1991Lifting,Lowen2005Obstruction}.
This complex has cohomology $J\otimes_R B$ in degree zero only,
because the endomorphism complex is bounded and $R$-flat and
represents the flat module $B$. Subtracting a degree-one cochain
with the indicated coboundary gives a differential. Iteration
yields a bounded projective lift $P$ over $\mathcal A$.

For a bounded complex $C$ of flat $S$-modules, suppose that
$R\otimes_S C$ has cohomology only in degree zero and that this
cohomology is $R$-flat. If $L$ is an $R$-module, then
$L\otimes_S C$ has cohomology only in degree zero, by derived
tensor over $R$. For an arbitrary $S$-module $L$, the filtration
by $I^jL$ gives the same concentration, since tensoring with $C$
is termwise exact. Taking $L=S$ and then arbitrary $L$ proves
that $C\simeq H^0(C)$ and that $H^0(C)$ is $S$-flat. In particular,
\[
 H^0(S'\otimes_S C)\simeq S'\otimes_S H^0(C)
\]
for every commutative $S$-algebra $S'$. The cycles and boundaries
are flat by descending along the bounded complex using the exact
sequences of cycles, boundaries and cohomology.

If $Q$ is another bounded projective lift, apply this argument
to $\operatorname{Hom}_{\mathcal A}^\bullet(P,Q)$. Its terms are
summands of finite sums of $\mathcal A$, and its special fibre
represents $B$ in degree zero. Thus a specified homotopy
equivalence of the special fibres lifts to a chain map.
Derived reduction to $R$ is conservative on $\mathcal D(S)$.
Indeed, if $R\otimes_S^{\mathbf L}C'=0$, then by associativity
$(I^j/I^{j+1})\otimes_S^{\mathbf L}C'=0$ for every $j$, and the
finite filtration of $S$ gives $C'=0$. Applying this to the
cone proves that the lifted map is a homotopy equivalence.

Let $H=\operatorname{End}_{\mathcal A}^\bullet(P)$ and put
$\mathcal B=H^0(H)$. The preceding flatness argument gives
flatness and base change for $\mathcal B$. The
quasi-isomorphisms $H\leftarrow\tau_{\leq0}H\to\mathcal B$
give a $\mathcal B$-$\mathcal A$-complex $\mathcal X$ with
underlying right complex $P$. These maps commute with base
change by flatness of the cycles, boundaries and cohomology.
The special fibre has the same underlying right complex and
the same left action on its derived endomorphisms as $X$.
After tensoring on the right with $Y$, it becomes a complex
with cohomology only in degree zero and with the two regular
$B$-actions there. This identifies the special fibre with $X$.

Put $\mathcal Y=\mathbf R\operatorname{Hom}_{\mathcal A}
(\mathcal X,\mathcal A)$. Right perfection gives the
isomorphism $\mathcal B\simeq
\mathcal X\otimes_{\mathcal A}^{\mathbf L}\mathcal Y$.
The other evaluation map is an isomorphism after derived
reduction to $R$, hence an isomorphism by conservativity.
These kernels lift $X,Y$ and are mutually inverse.
For the same $\mathcal A$, any such lift has underlying right
complex a lift of $T$, so its endomorphism algebra gives the
same marked class of $\mathcal B$. The construction for $Y$
therefore gives the inverse correspondence, and the formula
for $H^0(H)$ gives naturality. For the related lifting
correspondence on orders over a complete discrete valuation
ring, see \cite[Proposition~3.6]{Eisele2012Lifting}.

Now put $S_n=S/I^n$ and let $(\mathcal A_n)$ be a compatible
system over the $S_n$. The obstruction calculation gives
compatible projective lifts $P_n$ with compatible differentials.
Put $H_n=\operatorname{End}_{\mathcal A_n}^\bullet(P_n)$ and
$\mathcal B_n=H^0(H_n)$. These form a compatible flat system
by (1). For two choices of the system of lifts, the Hom
comparison gives compatible marked homotopy equivalences by
successive lifting. They induce an isomorphism of the systems
of marked algebras. The base change formula in (1) proves
naturality on these systems.

For inversion, put $K_n=\tau_{\le0}H_n$, let $(U_n)$ be a system
of projective lifts over $(\mathcal B_n)$ of the right $B$-complex
of $Y$, and put
$W_n=\operatorname{Hom}_{\mathcal A_n}^\bullet(P_n,\mathcal A_n)$.
At each stage $W_n$ is the right complex of the inverse tilting
kernel from (1). Thus $U_n\simeq W_n$ over $K_n$, via the
quasi-isomorphism $K_n\to\mathcal B_n$. These isomorphisms can
be chosen compatibly with the markings. Indeed both complexes
over $\mathcal B_n$ are perfect lifts of the same tilting
complex. The flatness argument for their Hom complex gives
surjectivity on degree-zero morphisms under reduction, and a
lift of an isomorphism is an isomorphism by conservativity.
The left actions of $\mathcal A_n$ on the $W_n$ give compatible
marked isomorphisms
\[
 \mathcal A_n\simeq\operatorname{End}_{\mathcal D(\mathcal B_n)}(U_n).
\]
Interchanging the roles of $X$ and $Y$ gives the other inverse
identity on systems.

If $S$ is complete and $I$ finitely generated, the complete
algebras $\mathcal A,\mathcal B$ are the inverse limits of the
systems, with reductions $\mathcal A_n,\mathcal B_n$, respectively
\cite[Tag \texttt{09B8}]{StacksProject}. The system comparison
thus gives the assertion for complete deformations.
Taking limits of the idempotent matrices shows that $P=\varprojlim_n P_n$
and $U=\varprojlim_n U_n$ are bounded complexes of finitely generated
projectives over $\mathcal A$ and $\mathcal B$, respectively.
Put $H=\operatorname{End}_{\mathcal A}^\bullet(P)=\varprojlim_n H_n$.
Surjectivity on terms and on $H^0(H_n)$ gives $H^0(H)=\mathcal B$
and $H^i(H)=0$ for $i\ne0$ by the inverse limit cohomology sequence.

Suppose that $\mathcal A$ and $\mathcal B$ are $S$-flat.
Write $K=\tau_{\le0}H$ and put
$W=\operatorname{Hom}_{\mathcal A}^\bullet(P,\mathcal A)$,
a right dg $K$-module. The derived base change of the restricted
$K$-module $U$ to $K_n$ is $U_n$ via the quasi-isomorphisms
to $\mathcal B$ and $\mathcal B_n$. As
$W=\mathbf R\!\varprojlim_n W_n$, adjunction gives
\[
 \mathbf R\operatorname{Hom}_{K}(U,W)
  \simeq \mathbf R\!\varprojlim_n
     \mathbf R\operatorname{Hom}_{K_n}(U_n,W_n).
\]
The compatible isomorphisms lift to a morphism $U\to W$ by the
inverse limit cohomology sequence. This morphism is an
isomorphism, since $U$ and $W$ are the derived inverse limits
of the two termwise surjective systems and its composites
to $W_n$ are given by the chosen isomorphisms at the finite
stages. The map $K\to\mathcal B$ gives $\mathcal X$ as in (1).
Its right dual
$\mathcal Y$ is right perfect over $\mathcal B$ by $U\simeq W$.
The underlying $S$-complex of
$\mathcal Y\otimes_{\mathcal B}^{\mathbf L}\mathcal X$
belongs to the thick subcategory of $\mathcal D(S)$ generated
by $P$. For $P$ and $\mathcal A$ the map
$C'\to\mathbf R\!\varprojlim_n(S_n\otimes_S^{\mathbf L}C')$
is an isomorphism by completeness and termwise flatness.
This property is closed under thick operations. The evaluation
map is an isomorphism after reduction to each $S_n$ by (1),
hence it is an isomorphism over $S$. The other evaluation
isomorphism follows from right $\mathcal A$-perfection.
When $S$ is Noetherian, the required flatness of the inverse
limits follows from \cite[Theorem~5.11]{Yekutieli2018Completion}.
\end{proof}

\begin{theorem}\label{thm:bg-standard-derived}
Let $R=k$ be a field of characteristic zero and let
$\mathcal O=k[[t_1,\ldots,t_d]]$. Let $\mathcal A,\mathcal B$ be
corresponding marked formal deformations over $\mathcal O$ from
Proposition~\ref{prop:standard-derived-deformations}. Then in
the homotopy category of dg algebras
\[
        \operatorname{Deform}_{\mathcal B}
             =\Phi_{X,*}\operatorname{Deform}_{\mathcal A}.
\]
\end{theorem}

\begin{proof}
The algebras $\mathcal A,\mathcal B$ are $\mathcal O$-flat by
Proposition~\ref{prop:standard-derived-deformations}.
Their reductions are flat over local Artin quotients of
$\mathcal O$, hence free, as required in
\cite[Section~1.2]{BezrukavnikovGinzburg2007}.
Use the same $\mathcal R$ and $\Lambda$ for both algebras in
the construction of \cite[(3.1.4)--(3.1.6)]{BezrukavnikovGinzburg2007}.
Write $K(W)=(\mathcal R\otimes_{\mathcal O}\mathcal R)
\otimes_\Lambda^{\mathbf L}W$. Thus the functor for
$\mathcal D$ is represented by
$K(W)\otimes_{\mathcal O}^{\mathbf L}\mathcal D$.
For $\mathcal D=\mathcal A,\mathcal B$ put
$\mathcal D'=\mathcal R\otimes_{\mathcal O}\mathcal D$.
Extend $\mathcal X,\mathcal Y$ to $\mathcal R$ by derived
base change and denote the resulting kernels by
$\mathcal X',\mathcal Y'$. There are natural isomorphisms
\[
\begin{split}
 \mathcal X'\otimes_{\mathcal A'}^{\mathbf L}
       \Theta_{\mathcal A}(W)
       \otimes_{\mathcal A'}^{\mathbf L}\mathcal Y'
 &\simeq K(W)\otimes_{\mathcal O}^{\mathbf L}
       (\mathcal X\otimes_{\mathcal A}^{\mathbf L}\mathcal Y)\\
 &\simeq \Theta_{\mathcal B}(W).
\end{split}
\]
For the first isomorphism, choose resolutions to compute the
derived tensors and use associativity of tensor products. The
two $\mathcal R$-actions are carried by $K(W)$. For the second,
use the evaluation isomorphism
$\mathcal X\otimes_{\mathcal A}^{\mathbf L}\mathcal Y\simeq\mathcal B$.
These isomorphisms are natural in $W$ on the dg models.
Extension of scalars at both ends along $\mathcal R\to k$
commutes with this comparison by derived base change. The two
extended kernels become $X$ and $Y$, so the resulting functor
on bimodules over the special fibres is $\Phi_X$. For the
augmented $\Lambda$-module $k$ the identification is the
diagonal identification induced by $X\otimes_A^{\mathbf L}Y\simeq B$.
Taking derived dg endomorphisms and composing with the same
Koszul morphism, as in
\cite[(3.2.1)]{BezrukavnikovGinzburg2007}, gives the stated equality.
\end{proof}

\subsubsection*{Graded Koszul duality}

For a nonnegatively graded $k$-algebra $D$ with finite-dimensional
homogeneous components, write $\operatorname{Def}^{\mathrm{gr}}_D(S)$
for the isomorphism classes of graded $S$-algebras with finite free
homogeneous components and a graded marking over $k$. Here $S$ is
a commutative local Artin $k$-algebra with residue field $k$, placed
in degree zero. Isomorphisms respect the grading and the marking.
We use the convention that the quadratic dual reverses the sides
of the degree-one bimodule \cite{PolishchukPositselski2005}.

\begin{prp-nn}[Graded Koszul duality]
Let $k$ be a field of arbitrary characteristic. Let $D$ be a Koszul
algebra with finite-dimensional homogeneous components and with
$D_0=E=k^r$. Write $D^!$ for its quadratic dual. Then quadratic
duality induces a natural isomorphism
\[
 \operatorname{Def}^{\mathrm{gr}}_D
   \simeq \operatorname{Def}^{\mathrm{gr}}_{D^!}.
\]
\end{prp-nn}

\begin{proof}
Let $\mathcal D$ be a marked graded deformation over $S$. Lift the
vertex idempotents to $\mathcal D_0$. The off-diagonal Peirce
components vanish by Nakayama's lemma, and in each diagonal
component the lifted idempotent is an $S$-basis. This identifies
$\mathcal D_0$ with $E_S=S\otimes_k E$. The identification is
unique with the given marking because the primitive idempotents
of a product of local rings are determined by their reductions.
Put $V=\mathcal D_1$ and let $I$ be the kernel of
$V\otimes_{E_S}V\to\mathcal D_2$.
By Nakayama's lemma in each weight, $V$ generates $\mathcal D$
over $E_S$. The surjection to $\mathcal D_2$ splits in each
Peirce component, so $I$ commutes with base change. The
surjection $T_{E_S}(V)/(I)\to\mathcal D$ reduces to an isomorphism
over $k$, since $D$ is quadratic. In each weight the target is
$S$-free, so reduction preserves the kernel of the surjection.
Nakayama's lemma gives $\mathcal D=T_{E_S}(V)/(I)$.

Consider the weight-$w$ part of the reduced $E_S$-bar complex
\[
 B_n(\mathcal D)_w=(\mathcal D_{>0}^{\otimes_{E_S}n})_w .
\]
It is a bounded complex of finite free $S$-modules in bar degrees
at most $w$. Its reduction has homology only in bar degree $w$
by Koszulity of $D$. Splitting off two-term contractible
complexes at unit entries of the differentials leaves a
complex whose differentials vanish modulo the maximal ideal.
Its only nonzero term is therefore a free module in bar degree
$w$. In particular, $H_w(B_\bullet(\mathcal D)_w)$ is free and
commutes with arbitrary Artin base change, and the differentials
have split images.

Write $V^\vee=\operatorname{Hom}_S(V,S)$ with the reversed
bimodule actions. In the split Peirce decomposition, dualizing
a tensor reverses its factors. For $w\geq2$ the differential from
$B_w(\mathcal D)_w$ to $B_{w-1}(\mathcal D)_w$ consists of the
adjacent multiplications $V\otimes_{E_S}V\to\mathcal D_2$ into
separate summands, with bar signs. The image of its transpose
is the sum of the adjacent insertions of $I^\perp$. Since the
differential has split image, it follows that
\[
 \bigl(T_{E_S}(V^\vee)/(I^\perp)\bigr)_w
   \simeq
 \operatorname{Hom}_S\bigl(H_w(B_\bullet(\mathcal D)_w),S\bigr).
\]
The same formula in weights zero and one follows from $E_S$ and
$V$, respectively. Thus the quadratic dual defines a marked
graded deformation of $D^!$ and commutes with scalar extension.
It sends marked isomorphisms to marked isomorphisms with the
direction reversed, so it induces a map on isomorphism classes.
The same construction applies to $D^!$, which is again Koszul.
Biduality and $(I^\perp)^\perp=I$ recover $\mathcal D$ with its
marking and give the inverse natural map.
\end{proof}

In characteristic zero, the assertion also holds under the
separability hypothesis on the degree-zero algebra in
\cite[Section~3.5]{Keller2003Higher}. Indeed, take Adams-degree-zero
cochains in Keller's $B_\infty$ comparison. On the dual algebra,
the cohomological and Adams degrees of its homogeneous elements
are opposite. Thus a cochain of Adams degree zero has internal
cohomological degree zero, and its total degree is its arity.
The two shifted weight-zero complexes, with their Gerstenhaber
brackets, are linked by dg Lie quasi-isomorphisms. The
Maurer--Cartan equation is associativity of the graded
perturbation of the multiplication, and the gauge action is a
graded change of coordinates reducing to the identity.
Normalized cochains fix the unit and are quasi-isomorphic to
all cochains in Adams degree zero. The dg Lie deformation
criterion \cite{ManettiDeformationLectures} gives the assertion.

The grading restriction is essential even if both algebras are
finite-dimensional. Let $Q$ have arrows $a:1\to2$, $b:2\to3$
and $c:1\to3$. The radical-square-zero algebra $D=kQ/(ba)$,
with its path-length grading, is Koszul, and $D^!=kQ^{\mathrm{op}}$
is a path algebra. The normalized Hochschild complex for $D$
relative to the separable vertex algebra $k^3$ has a
one-dimensional degree-two term, given by sending $b\otimes a$
to a multiple of $c$. The incoming differential is zero because
products of two arrows vanish and every vertex-bilinear
degree-one cochain takes arrows to arrows. The degree-three
term vanishes because there is no composable triple. Thus
$\mathrm{HH}^2(D,D)\simeq k$, with generator represented by
$ba=\epsilon c$ over $k[\epsilon]/(\epsilon^2)$. The path algebra
$D^!$ has a bimodule resolution of length one, so
$\mathrm{HH}^2(D^!,D^!)=0$. The ordinary marked deformation
functors have different tangent spaces.

\subsection{Module deformations over a local Artin ring}
\label{subsec:artin-module-deformations}

Let $R$ be a commutative local Artin ring with residue field $k$,
and let $A$ be a finite free $R$-algebra. For an $R$-algebra $S$
write $A_S=S\otimes_R A$. The use of a complete local coefficient
ring with a specified residue field is standard in deformation
theory \cite{Schlessinger1968,BleherDihedral2009,BleherChinburgSoto2019}.
Bleher's comparison for block algebras uses a stable equivalence
over the coefficient ring \cite[Lemma~2.2.2]{BleherDihedral2009}.
Here we fix an $A$-module over $R$ and mark the special fibre
over the whole ring $R$.

We denote by $\mathcal C_R^{\mathrm{aug}}$ the category whose
objects are complete local commutative Noetherian $R$-algebras
$S$ with the same residue field $k$ as $R$, equipped with a local
$R$-algebra augmentation $\epsilon_S:S\to R$. The structure
map $R\to S$ induces the specified identity on $k$, and
$R\to S\xrightarrow{\epsilon_S}R$ is the identity. Morphisms
are continuous local $R$-algebra homomorphisms that induce the
identity on $k$ and commute with the augmentations. We write
$\mathcal C_{R,\mathrm{Art}}^{\mathrm{aug}}$ for the full
subcategory of Artin objects.

\begin{definition}\label{def:artin-module-deformations}
Let $V$ be an $A$-module that is finite free over $R$.
For $S\in\mathcal C_R^{\mathrm{aug}}$, a lift of $V$ over $S$
is a pair $(L,\phi)$ consisting of a finite $A_S$-module $L$
that is free over $S$ and an $A$-linear isomorphism
$\phi:R\otimes_{S,\epsilon_S}L\to V$. Two lifts are equivalent
if an $A_S$-module isomorphism between them commutes with the
markings after reduction along $\epsilon_S$. Base change
defines the functor $\operatorname{Def}_{A,V/R}$ on
$\mathcal C_R^{\mathrm{aug}}$.

A universal deformation ring of $V$ over $R$ is an object
$U_R(A,V)$ of $\mathcal C_R^{\mathrm{aug}}$ that represents
$\operatorname{Def}_{A,V/R}$. Thus base change of the universal
lift gives natural bijections
\[
 \operatorname{Hom}_{\mathcal C_R^{\mathrm{aug}}}
       \bigl(U_R(A,V),S\bigr)
  \xrightarrow{\ \sim\ }\operatorname{Def}_{A,V/R}(S)
 \qquad(S\in\mathcal C_R^{\mathrm{aug}}).
\]
\end{definition}

For $S\in\mathcal C_R^{\mathrm{aug}}$ put
$S_n=S/\mathfrak m_S^n$. For $n$ sufficiently large the
augmentation of $S$ factors through $S_n$, since the maximal
ideal of $R$ is nilpotent. We use these augmented quotients
when applying arguments on Artin objects to a complete test ring.

\begin{lemma}\label{lem:artin-deformation-splitting}
Put $V_k=k\otimes_R V$. Suppose $V_k$ is Gorenstein
projective over $A_k$. For every finite projective $A$-module
$P$, adjoining the lift $P_S$ gives a natural bijection
\[
 \operatorname{Def}_{A,V/R}(S)
     \simeq \operatorname{Def}_{A,V\oplus P/R}(S).
\]
\end{lemma}

\begin{proof}
We first record the lifting step used here and below. Let
$S'\to S$ be a small extension of local Artin $R$-algebras
with residue field $k$ and kernel $J$. Let $L'$ be $S'$-free
with an $A_{S'}$-action, and let $P'$ be finite projective over
$A_{S'}$. Write $L_0=k\otimes_{S'} L'$ and $P_0=k\otimes_{S'}P'$.
Applying $\operatorname{Hom}_{A_{S'}}(L',-)$ to
\[
 0\longrightarrow J\otimes_k P_0\longrightarrow P'
       \longrightarrow S\otimes_{S'}P'\longrightarrow 0
\]
shows that the obstruction to lifting a map from
$S\otimes_{S'}L'$ to $S\otimes_{S'}P'$ lies in
$J\otimes_k\operatorname{Ext}^1_{A_k}(L_0,P_0)$.
Indeed, resolve $L'$ by finite free $A_{S'}$-modules. Its
syzygies are $S'$-free, by splitting over $S'$ at each stage.
Reducing the resolution to $k$ computes the Ext group with
coefficients in $J\otimes_kP_0$. A map in the other direction,
from the projective, lifts by projectivity. These assertions
can be iterated for an arbitrary surjection of local Artin
rings with fixed residue field, by factoring it into small
extensions.

Now let $L$ be a marked lift of $V\oplus P$ over $S$.
If $S$ is Artin, the preceding step lifts the marked
projection to $P$, since
$\operatorname{Ext}^1_{A_k}(V_k\oplus P_k,P_k)=0$.
For a complete $S$, start with a lift of the projection over
one $S_n$ through which $\epsilon_S$ factors. The same vanishing
allows us to lift it successively over $S_{n+1}$, $S_{n+2}$,
and so on. Passing to the limit gives an $A_S$-linear map
$L\to P_S$, since both modules are finite free over the
complete ring $S$. The inclusion $P\to V\oplus P$ lifts to a
map $P_S\to L$ by projectivity.
Their composite reduces to the identity of $P$, so it is
invertible by Nakayama's lemma. Adjusting the inclusion by
this inverse splits off $P_S$. The kernel is $S$-free and
inherits a marking by $V$. If two direct sums of lifts with
$P_S$ are isomorphic respecting the marking, the component
of the isomorphism between the two $V$-lifts reduces to the
identity under the markings. This component is invertible
by Nakayama's lemma. The construction by adjoining $P_S$ is
natural, and so is its inverse on isomorphism classes.
\end{proof}

\begin{proposition}\label{prop:artin-module-deformation-universal}
If $V_k$ is Gorenstein projective and the scalar map
$k\to\underline{\operatorname{End}}_{A_k}(V_k)$ is an
isomorphism, then $\operatorname{Def}_{A,V/R}$ has a universal
deformation ring.
\end{proposition}

\begin{proof}
Let $D$ be the functor of $k$-marked free lifts of $V_k$ on
local Artin $R$-algebras with residue field $k$. For $L$ over
$S$ write $E_S(L)=\operatorname{End}_{A_S}(L)$ and let
$\mathcal P_S(L)$ be its ideal of maps through finite projectives.
We first prove
\begin{equation}\label{eq:module-def-end}
 E_S(L)=S\,1_L+\mathcal P_S(L).
\end{equation}
At $S=k$ this is the stable endomorphism hypothesis. For a small
extension $S'\twoheadrightarrow S$ with kernel $J$ and a lift $L'$
over $S'$, put $L=S\otimes_{S'}L'$. If
$L\xrightarrow{a}P\xrightarrow{b}L$ is a projective factorization,
lift $P$ to a projective $P'$ by lifting an idempotent. Then
\[
\begin{gathered}
 a\in\operatorname{im}\bigl(
 \operatorname{Hom}_{A_{S'}}(L',P')\to
 \operatorname{Hom}_{A_S}(L,P)\bigr),\\
 b\in\operatorname{im}\bigl(
 \operatorname{Hom}_{A_{S'}}(P',L')\to
 \operatorname{Hom}_{A_S}(P,L)\bigr).
\end{gathered}
\]
The first inclusion uses Lemma~\ref{lem:artin-deformation-splitting}
and $\operatorname{Ext}^1_{A_k}(V_k,P_k)=0$. The second uses
projectivity of $P'$. Hence $\mathcal P_{S'}(L')\to\mathcal P_S(L)$
is surjective. Moreover, freeness over $S'$ gives
\[
 \ker\bigl(E_{S'}(L')\to E_S(L)\bigr)
  =\operatorname{Hom}_{A_{S'}}(L',JL')
  \simeq J\otimes_k\operatorname{End}_{A_k}(V_k).
\]
By the same projective lifting argument, every residual projective
factorization lifts to $S'$. Multiplying its lift by an element of
$J$ represents the corresponding factorization in the kernel.
Using \eqref{eq:module-def-end} at $k$, we obtain
\[
 \ker\bigl(E_{S'}(L')\to E_S(L)\bigr)
       \subseteq S'\,1_{L'}+\mathcal P_{S'}(L').
\]
Lifting the two summands in \eqref{eq:module-def-end} over $S$
and using this kernel inclusion proves \eqref{eq:module-def-end}
over $S'$. This completes the induction. It also proves surjectivity
of the reduction on endomorphisms. By Nakayama's lemma and iteration
over small extensions, any surjection $S'\twoheadrightarrow S$
as above gives
\begin{equation}\label{eq:module-def-aut-lift}
 \operatorname{Aut}^{1}_{A_{S'}}(L')
 \twoheadrightarrow \operatorname{Aut}^{1}_{A_S}(S\otimes_{S'}L'),
\end{equation}
where the superscript $1$ means that the automorphism induces the
identity on $V_k$ under the marking.

Let $S_1\to S_0\leftarrow S_2$ be a diagram of local Artin
$R$-algebras with residue field $k$, with $S_2\to S_0$ surjective.
Given lifts $L_1,L_2$ and a marked identification of their
reductions with a common lift $L_0$, the glued module is
\[
 L_1\times_{L_0}L_2
        \quad\text{over}\quad S_1\times_{S_0}S_2 .
\]
It is free by choosing bases with the same reduction to $S_0$.
The two module actions give its action over
$A_{S_1\times_{S_0}S_2}$. Two identifications over $S_0$ differ by
an element of $\operatorname{Aut}^{1}_{A_{S_0}}(L_0)$.
By \eqref{eq:module-def-aut-lift} this difference can be lifted
on the $S_2$-side. Consequently
\[
 D(S_1\times_{S_0} S_2)
    \xrightarrow{\ \sim\ }D(S_1)\times_{D(S_0)}D(S_2),
 \qquad
 D(k[\varepsilon]/(\varepsilon^2))
      \simeq\operatorname{Ext}^1_{A_k}(V_k,V_k).
\]
The tangent space is finite-dimensional. Schlessinger's criterion
\cite[Theorem~2.11]{Schlessinger1968} gives a complete local
Noetherian $R$-algebra $U$ pro-representing $D$.

The fixed lift $V$ defines $[V]\in D(R)$, corresponding to
$\epsilon_U:U\to R$. For $S\in\mathcal C_{R,\mathrm{Art}}^{\mathrm{aug}}$,
the fibre over this point is
\[
 \operatorname{Def}_{A,V/R}(S)
 \ \xrightarrow{\ \sim\ }\
 D(S)\times_{D(R)}\{[V]\}
 \ \xleftarrow{\ \sim\ }\
 \{f:U\to S\mid \epsilon_S f=\epsilon_U\},
\]
where the last set consists of continuous local $R$-algebra maps
inducing the identity on $k$. The first map forgets the marking
over $R$ and retains the marking over $k$. A lift in the fibre
admits an $R$-marking. Two choices differ by an element of
$\operatorname{Aut}^{1}_{A}(V)$, and
\eqref{eq:module-def-aut-lift} along $\epsilon_S$ identifies the
resulting marked lifts.

Now let $S\in\mathcal C_R^{\mathrm{aug}}$ and consider the augmented
quotients $S_n=S/\mathfrak m_S^n$ for $n$ sufficiently large.
A compatible family of lift classes has representatives $L_n$
with chosen isomorphisms $S_n\otimes_{S_{n+1}}L_{n+1}\simeq L_n$
respecting the markings. Choosing compatible bases yields
\[
 L=\varprojlim_n L_n\simeq S^{\operatorname{rank}_R V}
 \quad\text{as an $S$-module},\qquad R\otimes_S L\simeq V
 \quad\text{as an $A$-module}.
\]
The module actions are compatible and define the action of $A_S$
on $L$. If two lifts have equivalent reductions at every stage,
their marked isomorphisms over $S_n$ can be chosen compatibly.
Indeed, the required correction at level $n$ is an automorphism
equal to the identity after reduction to $R$. It lifts to level
$n+1$ by \eqref{eq:module-def-aut-lift} and still induces the
identity over $R$. Thus
\[
 \operatorname{Def}_{A,V/R}(S)
 \simeq \varprojlim_n\operatorname{Def}_{A,V/R}(S_n)
 \simeq \varprojlim_n
       \operatorname{Hom}_{\mathcal C_R^{\mathrm{aug}}}(U,S_n)
 \simeq \operatorname{Hom}_{\mathcal C_R^{\mathrm{aug}}}(U,S).
\]
This proves the assertion on the complete category.
\end{proof}

\begin{theorem}\label{thm:artin-module-deformation-morita}
Let $A,B$ be finite free $R$-algebras. Let $M$ and $N$ be
$R$-central $A$-$B$ and $B$-$A$ bimodules, respectively,
finitely generated projective on both sides, with
\[
 M\otimes_BN\simeq A\oplus P,\qquad
 N\otimes_AM\simeq B\oplus Q,
\]
where $P,Q$ are projective over $A\otimes_R A^{\mathrm{op}}$
and $B\otimes_R B^{\mathrm{op}}$, respectively.
Let $V$ be as in Definition~\ref{def:artin-module-deformations},
with $V_k$ Gorenstein projective, and put $W=N\otimes_A V$.
Tensoring with $N_S$ gives a natural isomorphism
\[
 \operatorname{Def}_{A,V/R}
       \simeq\operatorname{Def}_{B,W/R}.
\]
If the scalar map $k\to\underline{\operatorname{End}}_{A_k}(V_k)$
is an isomorphism, then this identifies the two universal
deformation rings as augmented $R$-algebras.
\end{theorem}

\begin{proof}
The module $W$ is $R$-free since it is an $R$-module summand of
a finite sum of copies of $V$. The same argument over $S$
shows that tensoring a free lift gives a free lift, with the
marking obtained by base change. Base change of the bimodule
identities to $k$ gives a stable equivalence of Morita type.
Thus $W_k$ is Gorenstein projective by
Theorem~\ref{thm:gpd-base-injective-semt}, and the stable
endomorphism rings of $V_k$ and $W_k$ are identified.

For a lift $L$ on the $A$-side the composite tensor sends
$L$ to $L\oplus(P_S\otimes_{A_S}L)$. The second summand is
projective over $A_S$ because it is a summand of a finite sum of
copies of $A_S\otimes_S L$. Its reduction is $P\otimes_AV$,
with the marking induced by the marking of $L$. A marked
projective lift is unique up to marked isomorphism, since
one can lift an isomorphism on the special fibre by
projectivity and detect invertibility after reduction.
Consequently the composite agrees on classes with adjoining
the fixed lift of $P\otimes_A V$. This is a bijection by
Lemma~\ref{lem:artin-deformation-splitting}. Starting with $W$,
the other composite is likewise a bijection onto the
deformations of $W\oplus(Q\otimes_B W)$. Thus the map from
the deformations of $W$ to those of $M\otimes_BW$ is
injective, and its composite with the map from the
deformations of $V$ is bijective. It follows that both
maps in this composite are bijective. The tensor maps
commute with base change, which gives the required
natural isomorphism. Proposition~\ref{prop:artin-module-deformation-universal}
and uniqueness of pro-representing objects give the ring
comparison.
\end{proof}

For $R=k$, the universal ring is the usual ring in
\cite{BekkertGiraldoVelez2020}. The field-level comparison for
singular equivalences of Morita type with level inducing an
equivalence on stable Gorenstein projectives is established in
\cite{PanVelez2026}. The theorem above concerns level zero over
$R$ and retains the chosen lift over $R$ in the marking.

\begin{proposition}\label{prop:module-deformation-effective-base}
In the setting of Theorem~\ref{thm:artin-module-deformation-morita},
allow $V$ to be nonfree over $R$. Put $I=\operatorname{Ann}_R(V)$
and $\bar R=R/I$, and suppose $V$ is finite free over $\bar R$.
Then $W=N\otimes_A V$ is finite free over $\bar R$ and has
annihilator $I$. On augmented $R$-test rings define lifts by
requiring the underlying $S$-module to be isomorphic to
$S\otimes_R V$, still with a marking over $R$. The resulting
functor is naturally
\[
  S\longmapsto \operatorname{Def}_{A/IA,V/\bar R}(S/IS).
\]
Tensoring identifies this functor with the corresponding
functor for $W$. If the scalar map on the stable endomorphism
ring of $V_k$ is an isomorphism, the two functors over the
augmented $\bar R$-test category have isomorphic universal
rings.
\end{proposition}

\begin{proof}
The inclusion $\operatorname{Ann}_R(V)\subseteq
\operatorname{Ann}_R(W)$ follows from $R$-centrality. The
reverse inclusion follows because $V$ is a summand of
$M\otimes_B W$. Moreover, $W$ is an $\bar R$-module summand
of a finite sum of copies of $V$, so it is free over
$\bar R$. An $S$-module of the indicated type is annihilated
by $IS$ and is free over $S/IS$. Its action and marking
therefore factor through the corresponding quotient
algebras, and the objects and isomorphisms in the two
lift problems agree. Base change to $\bar R$ preserves the
bimodule identities and projectivity hypotheses of
Theorem~\ref{thm:artin-module-deformation-morita}, which now
gives the comparisons.
\end{proof}

\begin{example}\label{ex:artin-module-deformation-nonrepresentable}
Let $R=k[\pi]/(\pi^2)$ and $A=R[x]/(x^2)$. Let $V=k$ with
$xV=0$ and the natural $R$-action. The fixed-underlying-module
functor of Proposition~\ref{prop:module-deformation-effective-base}
is represented over $\bar R=k$ by $k[[u]]/(u^2)$, with
augmentation $u\mapsto0$. Its pullback along $S\mapsto S/\pi S$
has no universal augmented $R$-algebra, even when the $R$-test
rings are required to be free over $R$.

Indeed, put $S_n=R[t]/(t^n)$, augmented by $t\mapsto0$.
A lift is determined by the scalar $z$ giving the action of
$x$ on the rank-one $S_n/\pi S_n$-module. Thus
\[
 z\in t\,k[t]/(t^n),\qquad z^2=0.
\]
A module isomorphism preserves $z$. The lift $z=t$ at $n=2$
does not extend to $n=3$, since the square of any lift of $t$
has nonzero coefficient of $t^2$ in $k[t]/(t^3)$.

Suppose that an augmented $R$-algebra $U$ represents the
pullback, and let $J=\ker(U\to R)$. The endomorphism of $S_2$
given by $t\mapsto\pi t$ sends every deformation to the
trivial one. By naturality and the fact that there is a
unique deformation over $R$, for every augmented $f:U\to S_2$
the composite with this endomorphism annihilates $J$.
Because $f(J)\subset Rt$ and the annihilator of $\pi$ in
$R$ is $\pi R$, it follows that $f(J)\subset \pi Rt$.
The quotient $S_3\to S_2$ restricts to an isomorphism from
$R+\pi Rt$ in $S_3$ to $R+\pi Rt$ in $S_2$. Both are
subalgebras because the square of $\pi Rt$ is zero.
Hence every such $f$ lifts to an augmented map $U\to S_3$.
This contradicts the failure to lift $z=t$. In particular
pro-representability over the quotient base and invariance
on the original test category need not give
pro-representability on that original category.
\end{example}

\begin{example}\label{ex:artin-module-deformation-syzygy}
Take $R=\mathbb Z/8\mathbb Z$ and $A=R[x]/(x^3)$. Let
$V=A/(x-2)$, a free $R$-module of rank one. We give a stable
self-equivalence of Morita type which takes $V$ to a free
$R$-module $W$ of rank two, and compute both universal rings.

Put $A^e=R[u,v]/(u^3,v^3)$, where $u$ acts on the left and
$v$ on the right, and let
\[
 K=\ker(A^e\xrightarrow{\mu}A)=(u-v),
 \qquad h=u^2+uv+v^2.
\]
Here $\mu$ sends $u$ and $v$ to $x$. As a right $A$-module,
$A^e$ is free with basis $1,u,u^2$, and
\[
 (u-v)(a+bu+cu^2)=-va+(a-vb)u+(b-vc)u^2
 \qquad(a,b,c\in R[v]/(v^3)).
\]
It follows that $\operatorname{ann}_{A^e}(u-v)=A^eh$ and that
$\bar a\mapsto h a$ embeds $A=A^e/(u-v)$ into $A^e$. In
particular there are exact bimodule sequences
\[
 0\longrightarrow K\longrightarrow A^e
   \xrightarrow{\mu}A\longrightarrow0,\qquad
 0\longrightarrow A\xrightarrow{\cdot h}A^e
   \xrightarrow{\cdot(u-v)}K\longrightarrow0.
\]
The first splits on both sides over $A$, and $K$ is free of
rank two on each side. Tensor the first sequence on the right
with $K$. Since $A^e\otimes_A K\simeq (A^e)^{\oplus2}$ as
an $A$-bimodule, Schanuel's lemma applied to the two projective
presentations of $K$ gives
\[
 (K\otimes_A K)\oplus A^e
      \simeq A\oplus (A^e)^{\oplus2}.
\]
The enveloping algebra is Artinian, so Krull--Schmidt
cancellation gives $K\otimes_A K\simeq A\oplus A^e$.
Thus $(K,K)$ defines a stable equivalence of Morita type.
Reduction of $V$ to $k=\mathbb F_2$ is the nonprojective
simple over $k[x]/(x^3)$ with stable endomorphism ring $k$.
It is Gorenstein projective because $k[x]/(x^3)$ is symmetric.

Tensoring the first sequence with $V$ identifies
\[
 W=K\otimes_A V
    \simeq\ker\bigl(A\xrightarrow{x\mapsto2}R\bigr)
    =(x-2)A.
\]
For a test ring $S\in\mathcal C_R^{\mathrm{aug}}$, a marked
lift of $V$ has rank one and $x$ acts by a unique scalar
$c\in S$ with
\[
 \epsilon_S(c)=2,\qquad c^3=0.
\]
Putting $c=2+t$ yields the augmented universal ring
\[
 U_R(A,V)=
   (\mathbb Z/8\mathbb Z)[[t]]/(t^3+6t^2+4t),
 \qquad \epsilon_U(t)=0.
\]

We compute the deformation functor of $W$ independently.
The $R$-basis $w=x-2$, $xw$ of $W$ gives the relation
$x^2w=-2xw-4w$. In any marked lift $L$, choose $e\in L$
lifting $w$ under the marking. Then $(e,xe)$ is a basis
over $S$ by Nakayama's lemma. There are unique $c,d\in S$
such that
\[
 x^2e=-cxe-de,\qquad
 \epsilon_S(c)=2,\quad \epsilon_S(d)=4.
\]
The $S[x]$-module $L$ identifies with $S[x]/(x^2+cx+d)$.
The monic polynomial $x^2+cx+d$ generates the annihilator
of $L$ in $S[x]$. Thus its coefficients do not depend on the
choice of $e$. The condition that $x^3$ act by zero is
\[
 x^3=(x-c)(x^2+cx+d)+(c^2-d)x+cd,
 \qquad
 d=c^2,\quad c^3=0.
\]
Conversely any such $c$ gives a marked rank-two lift with
this cyclic presentation. Hence $U_R(A,W)$ has the same
presentation as $U_R(A,V)$.

To identify the isomorphism of Theorem~\ref{thm:artin-module-deformation-morita},
tensor the lift $A_S/(x-c)$ with $K_S$. The result is the
ideal $(x-c)A_S$, with basis $x-c$, $x(x-c)$ and action matrix
\[
 \begin{pmatrix}0&-c^2\\[1pt]1&-c\end{pmatrix}.
\]
Thus the tensor comparison preserves the parameter $t=c-2$.
The two modules have different lengths on the residue fibre,
so this is a nontrivial comparison on stable module classes.
\end{example}

\medskip
The single-residue-field convention above uses a local Artin
base. For a general nonzero commutative Artin ring
$R=\prod_iR_i$, the corresponding test category is the product
of the categories $\mathcal C_{R_i}^{\mathrm{aug}}$, with the
prescribed residue field on each factor. A lift is finite
projective over the resulting complete semilocal test ring,
which is equivalent to freeness on each factor. The central
idempotents of $R$ decompose the marking and the module
action. Thus the functor, its universal ring when it exists,
and the Morita-type comparison are obtained componentwise.

\section{Connected outer automorphism groups}
\label{sec:outer-groups}

Let $k$ be a field and let $A$ be a finite-dimensional $k$-algebra.
If $k$ is perfect, we put
\[
 \operatorname{Out}^0_k(A)=
 \bigl(\underline{\operatorname{Aut}}_k(A)_{\mathrm{red}}
       /\underline{\operatorname{Inn}}_k(A)\bigr)^0
\]
for the reduced connected algebraic group. We omit the subscript
$k$ if $k$ is algebraically closed. We also consider the identity
component with its full group-scheme structure. Rouquier proved
invariance of the connected group under derived equivalences and
under stable equivalences of Morita type between self-injective
algebras \cite[Theorems~4.6 and~4.15]{Rouquier2011Automorphisms}.
For a reduced pair between algebras without semisimple direct
factors we use projectivity and injectivity of its error bimodules
over their enveloping algebras
\cite[Lemma~4.1(3)]{CPX2012}.

\subsection{Splitting projective-injective summands in families}

\begin{lemma}\label{out:local-pi-splitting}
Let $\Lambda$ be a finite-dimensional algebra over an algebraically
closed field $k$, let $S$ be a finite-type $k$-scheme, and let $\mathcal X$ be a
family of $\Lambda$-modules which is locally free of finite rank
over $S$. Suppose that, at a closed point $s\in S$,
\[
 \mathcal X(s)\simeq V\oplus E,
\]
where $E$ is a projective-injective $\Lambda$-module. On an open
neighbourhood $U$ of $s$ there is a decomposition
\[
 \mathcal X|_U\simeq(E\otimes_k\mathcal O_U)\oplus\mathcal X'.
\]
\end{lemma}

\begin{proof}
Work on an affine neighbourhood with coordinate ring $R$, and
write $X$ for the corresponding $(\Lambda\otimes_kR)$-module.
Projectivity lifts the fibre inclusion $i_0:E\to X(s)$ to a map
$i:E\otimes_kR\to X$. Put $DE=\operatorname{Hom}_k(E,k)$.
Since $E$ is injective, $DE$ is a projective right $\Lambda$-module,
and tensor--Hom adjunction gives
\[
 \operatorname{Hom}_{\Lambda\otimes_kR}(X,E\otimes_kR)
 \simeq\operatorname{Hom}_R(DE\otimes_\Lambda X,R).
\]
The module $DE\otimes_\Lambda X$ is a direct summand of a finite
sum of copies of $X$, and hence is finite projective over $R$.
The displayed Hom therefore commutes with passage to the fibre.
Consequently the fibre projection $p_0:X(s)\to E$ lifts to
$p:X\to E\otimes_kR$. The endomorphism $pi$ specializes to
$1_E$ and is invertible after shrinking the neighbourhood.
Replace $p$ by $(pi)^{-1}p$. Then $pi=1$, which gives the required
splitting.
\end{proof}

\subsection{Stable conjugation and algebraic families}

Let $\mathcal L(A)$ be the category of finite-dimensional
$A$-bimodules which are projective as left and as right modules.
Its quotient by maps through $A^e$-projective bimodules is monoidal:
tensoring an enveloping-projective bimodule with an object of
$\mathcal L(A)$ on either side again gives an enveloping-projective
bimodule. Denote the group of tensor-invertible isomorphism classes
in this quotient by $\underline{\operatorname{Pic}}(A)$.

\begin{lemma}\label{out:stable-conjugation}
Suppose that $M,N$ satisfy \eqref{eq:MT} over a field. The maps
\[
 [X]\longmapsto[N\otimes_AX\otimes_AM],\qquad
 [Y]\longmapsto[M\otimes_BY\otimes_BN]
\]
induce inverse group isomorphisms between
$\underline{\operatorname{Pic}}(A)$ and
$\underline{\operatorname{Pic}}(B)$.
If $k$ is algebraically closed and $A$ has no semisimple direct
factor, then the natural map
$\operatorname{Pic}(A)\to\underline{\operatorname{Pic}}(A)$ is
injective.
\end{lemma}

\begin{proof}
The functor $\Phi(X)=N\otimes_AX\otimes_AM$ preserves
left-right projectivity and sends enveloping-projectives to
enveloping-projectives. Moreover,
\[
 \Phi(X)\otimes_B\Phi(Y)
 \simeq\Phi(X\otimes_AY)
 \oplus N\otimes_AX\otimes_AP\otimes_AY\otimes_AM,
\]
and the last summand is $B^e$-projective. The reverse construction
has the same property. Expanding the two composites by
\eqref{eq:MT} proves the first assertion.

A direct summand of the regular bimodule $A$ is a central direct
algebra factor. Such a summand is $A^e$-projective only if that
factor is separable, and hence semisimple. Thus $A$ has no
enveloping-projective summand under the stated hypothesis. Neither
does any invertible $A$-bimodule: tensoring with its inverse would
give such a summand of $A$. Stable isomorphism and Krull--Schmidt
cancellation now imply actual isomorphism for invertible
bimodules. This proves injectivity.
\end{proof}

For $\alpha\in\operatorname{Aut}_k(A)$, let $A_\alpha$ be the
regular left module with right action $x\cdot a=x\alpha(a)$.
Then $A_\alpha\otimes_AA_\beta\simeq A_{\alpha\beta}$, and
$A_\alpha\simeq A_\beta$ as bimodules precisely when $\alpha$ and
$\beta$ have the same outer class.

\begin{proposition}\label{out:pi-errors}
Let $k$ be algebraically closed. Let $A,B$ be finite-dimensional
$k$-algebras without semisimple direct factors, and let $M,N$
satisfy \eqref{eq:MT}. If the error bimodules $P,Q$ are injective
over $A^e,B^e$, respectively, stable conjugation induces an
isomorphism of algebraic groups
\[
 \operatorname{Out}^0(A)\simeq\operatorname{Out}^0(B).
\]
\end{proposition}

\begin{proof}
Put $G=(\underline{\operatorname{Aut}}_k(A)_{\mathrm{red}})^0$.
The universal automorphism on $G$ defines a family $\mathcal A$
with fibre $A_\alpha$ at $\alpha$. Set
\[
 \mathcal X=N\otimes_A\mathcal A\otimes_AM.
\]
This family is locally free over $G$. Indeed, twisting the
projective right module $N$ by the universal automorphism gives
a projective right $A\otimes_k\mathcal O_G$-module; tensoring
with $M$ makes $\mathcal X$ a direct summand of a finite sum of
copies of the constant underlying family of $M$.
At the identity, $\mathcal X(1)\simeq B\oplus Q$.
Lemma~\ref{out:local-pi-splitting}, applied to $\Lambda=B^e$,
gives an identity neighbourhood $U\subseteq G$ and a decomposition
\[
 \mathcal X|_U\simeq
 (Q\otimes_k\mathcal O_U)\oplus\mathcal W,
 \qquad\mathcal W(1)\simeq B.
\]

Shrink $U$ so that $\mathcal W$ is the regular left
$B\otimes_k\mathcal O_U$-module. To do this, lift the image of
$1\in B$ in the fibre isomorphism and restrict to the open set
where the resulting left-module map is invertible.
The right action then gives a regular family of algebra
endomorphisms of $B$. They are automorphisms on a smaller
identity neighbourhood, since their determinants do not vanish
at $1$. Passing to outer classes gives a regular map
\[
 f_U:U\longrightarrow\operatorname{Out}^0(B).
\]
The target is the identity component because $U$ is irreducible
and $f_U(1)=1$. A change of left-module trivialization changes
only the inner class. By construction, $f_U(\alpha)$ is
characterized by the stable bimodule class
$[N\otimes_AA_\alpha\otimes_AM]$.

Lemma~\ref{out:stable-conjugation} gives an abstract group
homomorphism
\[
 \lambda:G(k)\longrightarrow\underline{\operatorname{Pic}}(B).
\]
View $\operatorname{Out}^0(B)(k)$ as a subgroup of the target,
using the injectivity in that lemma. Its inverse image under
$\lambda$ contains $U(k)$. A connected algebraic group is
generated by any identity neighbourhood: for $g\in G(k)$,
irreducibility gives $U\cap gU\ne\varnothing$, and hence
$g\in U(k)U(k)^{-1}$. Thus $\lambda$ takes every point into
$\operatorname{Out}^0(B)(k)$ and defines a unique group
homomorphism $f$ on $k$-points extending $f_U$.

This homomorphism is regular. On every translate $gU$ it is
given by $f(gu)=f(g)f_U(u)$. These regular expressions agree on
closed points of overlaps, and therefore agree as morphisms
because the source is reduced. They define an algebraic-group
morphism $G\to\operatorname{Out}^0(B)$. For an inner
automorphism $\alpha$, the bimodule $A_\alpha$ is isomorphic
to $A$, so $f(\alpha)=1$. Hence $f$ descends through the
algebraic-group quotient to
\[
 \bar f:\operatorname{Out}^0(A)\longrightarrow
             \operatorname{Out}^0(B).
\]

The reverse construction, now using $P$, gives a regular
morphism $\bar g$ in the opposite direction. Stable conjugation
and the injectivity of the ordinary Picard groups imply that
both composites are the identity on closed points. Reducedness
then makes them the identity morphisms. Thus $\bar f$ and
$\bar g$ are inverse algebraic-group isomorphisms. In particular,
the proof supplies a regular inverse also in positive
characteristic.
\end{proof}

\begin{theorem}\label{out:general-semt}
Let $k$ be algebraically closed and let $A,B$ be finite-dimensional
$k$-algebras. A $k$-linear stable equivalence of Morita type
induces an isomorphism
\[
 \operatorname{Out}^0(A)\simeq\operatorname{Out}^0(B)
\]
of reduced algebraic groups. Neither self-injectivity nor a
Gorenstein hypothesis is required.
\end{theorem}

\begin{proof}
First suppose there are no semisimple direct factors. Delete
projective mixed-bimodule summands from $M,N$. Krull--Schmidt
cancellation shows that the reduced pair still satisfies
\eqref{eq:MT}; see also \cite[proof of Theorem~1.2]{CPX2012}.
By \cite[Lemma~4.1(3)]{CPX2012}, the new errors are
projective-injective over the enveloping algebras. Apply
Proposition~\ref{out:pi-errors}.

In general write $A=A_0\times S$ and $B=B_0\times T$, where
$S,T$ are semisimple and $A_0,B_0$ have no semisimple direct
factors. Let $e,f$ be the central idempotents of $A_0,B_0$.
The restricted bimodules $eMf,fNe$ are projective on both sides.
The composite through $T$ is enveloping-projective over $A_0$:
since $T$ is separable, its tensor product is a direct summand
of the corresponding outer tensor product over $k$. Similarly
the composite through $S$ is enveloping-projective over $B_0$.
Restricting \eqref{eq:MT} and applying Krull--Schmidt therefore
gives actual Morita-type identities for $eMf,fNe$. The preceding
case applies to $A_0,B_0$.

A connected automorphism group fixes the central block
idempotents, and a split semisimple algebra has trivial connected
outer automorphism group. Hence
$\operatorname{Out}^0(A)=\operatorname{Out}^0(A_0)$ and
$\operatorname{Out}^0(B)=\operatorname{Out}^0(B_0)$.
The isomorphism is characterized by stable conjugation with the
original pair. Indeed, even in the presence of semisimple factors,
$\operatorname{Out}^0(A)$ embeds in
$\underline{\operatorname{Pic}}(A)$: its automorphisms act
innerly on each semisimple block, and cancellation detects
their classes on the other blocks.
\end{proof}

\begin{corollary}\label{out:perfect-field}
Let $k$ be perfect and let $A,B$ be finite-dimensional $k$-algebras.
A $k$-linear stable equivalence of Morita type induces an
isomorphism of reduced connected algebraic $k$-groups
\[
 \operatorname{Out}^0_k(A)\simeq\operatorname{Out}^0_k(B).
\]
\end{corollary}

\begin{proof}
Extend scalars to an algebraic closure $\bar k$. Over a perfect
field the reduced automorphism group is smooth, and reduction
commutes with this scalar extension. The inner group
$A^\times/Z(A)^\times$ is smooth and geometrically connected.
Formation of the quotient and of its identity component commutes
with extension to $\bar k$. Thus the scalar extension of
$\operatorname{Out}^0_k(A)$ is
$\operatorname{Out}^0(A\otimes_k\bar k)$, and similarly for $B$.

Theorem~\ref{out:general-semt} supplies an isomorphism over
$\bar k$. Its characterization by stable conjugation with the
original $k$-defined bimodules makes it Galois equivariant on
$\bar k$-points. Since the groups are reduced, this is
equivariance of morphisms. Affine Galois descent applies to
the isomorphism and its inverse, and gives the assertion.
\end{proof}

\begin{remark}\label{out:scope}
The theorem concerns the algebraic group, including its
dimension, maximal torus rank, unipotent radical and reductive
quotient. The full outer automorphism group is not invariant
even under Morita equivalence: $k\times M_2(k)$ and $k\times k$
are Morita equivalent, whereas their outer groups are $1$ and
$C_2$, respectively. See also
\cite[Remark~4.3]{Rouquier2011Automorphisms}.
Over a nonclosed field, rational points of an algebraic-group
quotient are not silently identified with a quotient of the
groups of rational points.

In positive characteristic the Lie algebra of the reduced group
need not be all of $\mathrm{HH}^1(A)$. For instance, in
characteristic two, $A=k[x]/(x^2)$ has two-dimensional
$\mathrm{HH}^1(A)$, while its reduced outer group is
$\mathbb G_m$ and has one-dimensional Lie algebra. Thus neither
an abstract Lie algebra comparison nor an assertion about the
full nonreduced group scheme can replace the family argument.
The theorem also does not assert an Artin-base relative
group-scheme analogue.

For an almost $\nu$-stable derived equivalence over a field,
outer-group invariance already follows from Rouquier's derived
invariance theorem. The extension proved here applies when only
a classical stable equivalence of Morita type is given.
\end{remark}

\subsection{Formal outer automorphisms and integrability}

Let $k$ be algebraically closed. For a local Artin $k$-algebra $T$
with residue field $k$, put $A_T=A\otimes_k T$ and let
\[
 \widehat{\operatorname{Out}}_A(T)=
 \ker\bigl(\operatorname{Out}_T(A_T)\longrightarrow
                         \operatorname{Out}_k(A)\bigr),
\]
where $\operatorname{Out}_T(A_T)$ means $T$-algebra automorphisms
modulo conjugation by units of $A_T$. This is the formal completion
of the outer group scheme at the identity. Indeed, the inner group
scheme $A^\times/Z(A)^\times$ is smooth and its torsors over $T$
have sections, by smooth lifting from the algebraically closed
residue field. The same reasoning identifies inner $T$-points
with the image of $A_T^\times$. Representatives of the kernel
can be chosen congruent to the identity, since units lift over
nilpotent ideals.

\begin{theorem}\label{thm:formal-outer}
A stable equivalence of Morita type between finite-dimensional
$k$-algebras induces a natural group isomorphism
$\widehat{\operatorname{Out}}_A\simeq
 \widehat{\operatorname{Out}}_B$.
Its tangent map is the degree-one Hochschild comparison.
In particular the comparison preserves the subspaces of
classes integrable to each finite order, and the classes
integrable over $k[[t]]$.
\end{theorem}

\begin{proof}
First remove the semisimple factors and reduce the bimodule pair
as in Theorem~\ref{out:general-semt}. These operations do not
alter the formal comparison. An automorphism congruent to the
identity fixes the central block idempotents, since idempotents
lift uniquely in the center over a nilpotent ideal. For a matrix
algebra over $T$, the images of the matrix units determine a
decomposition into free rank-one summands of the defining module.
Choosing a generator for one of them and transporting it by the
matrix units gives a conjugating basis. Thus these factors have
trivial formal outer functor. Adjoining projective mixed bimodule
summands, or tensoring through a separable factor, changes the
tensor comparison only by projective summands. These remain
projective after scalar extension to $T$ and after twisting.

Let $\alpha\in\operatorname{Aut}_T(A_T)$ reduce to the identity.
The $T$-free bimodule
$X=N_T\otimes_{A_T}(A_T)_\alpha\otimes_{A_T}M_T$ reduces to
$B\oplus Q$. The maps splitting off $Q$ lift to maps with
$Q_T$ as source or target. For the inclusion use projectivity
of $Q$ over $B^e$. For the projection, the identity
\[
 \operatorname{Hom}_{B_T^e}(X,Q_T)
   \simeq \operatorname{Hom}_T(D_kQ\otimes_{B^e}X,T)
\]
and projectivity of $D_kQ$ show that reduction on this Hom is
surjective. The tensor product on the right is finite projective
over $T$ since it is a summand of a finite sum of copies of $X$.
Here $X$ is $T$-free by the one-sided projectivity argument in
Proposition~\ref{out:pi-errors}. The composite of the lifted
maps reduces to $1_Q$ and is invertible. The complement $W$
is $T$-free, and a lift of the left generator of its special
fibre gives an isomorphism $B_T\simeq W$ of left modules by
Nakayama and equality of ranks. The right action gives a twist
$\beta$ reducing to the identity, since its underlying linear
map is invertible modulo the maximal ideal.

The stable bimodule class of $X$ determines the outer class of
$\beta$. Indeed, a nonzero projective summand of a twisted regular
$B_T$-bimodule would reduce to a nonzero projective summand of
its special fibre, by Nakayama's lemma. After untwisting, this
would give a projective summand of the regular $B$-bimodule.
Thus Krull--Schmidt cancellation over $B_T^e$ identifies any
two stably isomorphic twisted regular bimodules, and an
isomorphism between them identifies their twists up to an inner
automorphism. Uniqueness also gives compatibility with maps of
test rings. The tensor identities over $T$ induce inverse monoidal
comparisons modulo projectives by the argument of
Lemma~\ref{out:stable-conjugation}. They therefore give inverse
homomorphisms of the formal outer groups.

For $T=k[\epsilon]/(\epsilon^2)$ the twist by
$1+\epsilon d$ represents the extension class of the derivation
$d$, with the convention of twisting the right action. Tensoring
gives the extension of $B\oplus Q$ by itself. The constant summand
$Q_T$ has zero extension class. The $B$-component is therefore
the compression of the tensor image, which identifies the tangent
map with $\phi_1$. The compression on positive cohomology agrees
with the original pair: when a projective mixed summand is
deleted, the two-sided tensor image changes by direct summands
which are projective at the unit and hence have zero image in
positive Ext. Changing the splittings of the regular summands
does not change compression on this image, by
Lemma~\ref{lem:hh-annihilation} and the projective-center
calculation in Theorem~\ref{thm:hh-cohomology-algebra}.

A class is integrable to order $r$ if it lies in the image of
reduction from $k[t]/(t^{r+1})$ to the dual numbers. A compatible
system of outer classes has compatible representatives. At each
stage, lift a unit that conjugates the reduction of the new
representative to the preceding one, and correct by that unit.
Taking the inverse limit gives an automorphism over $k[[t]]$.
Naturality proves both integrability assertions.
\end{proof}

Briggs and Rubio y Degrassi prove invariance of the Lie algebra
of integrable classes for self-injective algebras
\cite{BriggsRubioIntegrable}. The formal comparison applies
with general finite-dimensional algebras and retains the group
functor on local Artin test rings.

\begin{corollary}\label{thm:outer-group-scheme}
Over a perfect field a classical stable equivalence of Morita
type induces an isomorphism of the full identity components
of the outer automorphism group schemes, with their scheme
structures.
\end{corollary}

\begin{proof}
First work over an algebraically closed field and remove the
semisimple factors and reduce the pair as above. Those factors
have trivial connected outer group scheme: this follows from
the matrix-unit argument over local Artin test rings and the
corresponding statement over $k$.

Let $G=\underline{\operatorname{Aut}}(A)^0$, with its full
scheme structure, and let $H=\underline{\operatorname{Out}}(B)^0$.
Choose an affine neighbourhood $U$ of the identity in $G$.
The family argument of Proposition~\ref{out:pi-errors} applies
over its coordinate ring. Projectivity and the base-change
calculation in Lemma~\ref{out:local-pi-splitting} lift the
inclusion and projection of $Q$. The subsequent invertibility
conditions are open determinant conditions on $G$, so the
construction applies on the possibly nonreduced open
subscheme. Trivializing the complementary left module gives a
morphism $h_U:U\to H$. On local Artin test rings its value is
determined by the stable class of the transported twisted
regular bimodule. The cancellation argument of
Theorem~\ref{thm:formal-outer} applies at every closed point,
since a twisted regular bimodule has no nonzero projective
summand.

For $g\in G(k)$ let $h(g)$ be the value of the reduced
comparison. Define the morphism on $gU$ by
$h(gu)=h(g)h_U(u)$. The translates cover $G$, and the monoidal
tensor comparison gives the same characterization on local
Artin test rings for each translate. Two such morphisms agree on every Artin
test point in the overlap, so they agree as morphisms. To
see this detection of equality, pull back the ideal of the
diagonal of the separated target. At every closed point it
vanishes in all the quotients of the local ring by powers of
its maximal ideal. The Krull intersection theorem makes the
ideal zero, and closed points detect its support. Thus the
local maps glue. The same test on $G\times G$ proves the group
law. On the inner subgroup the map is trivial, since an inner
twist is the regular bimodule, and inner classes lift to units
over the test rings. The identity component of the automorphism
group scheme maps onto that of its outer quotient as an fppf
sheaf: the kernel is smooth and connected, and over $k$
the inverse image of the identity component is connected
by smoothness, surjectivity and connected fibres. Therefore
the map descends to the identity component of the outer
group scheme of $A$.

Repeat with $A,B$ interchanged. The two composites agree with
the identity on all local Artin test points, by stable
cancellation. Thus the scheme morphisms are inverse.
Over a perfect field extend scalars to an algebraic closure.
Formation of the inner group, its quotient and the identity
component commutes with this extension. The characterization
by stable conjugation on local Artin test rings makes both
morphisms Galois equivariant, and descent finishes the proof.
\end{proof}

\section{Algebraic \texorpdfstring{$K$}{K}-theory and relative singularity categories}
\label{sec:k-theory}

For a finite-dimensional algebra, the quotient by its radical
induces maps on ordinary $K$-groups. We first compare their kernels
by showing that the error bimodules act trivially on them. We then
use the $R$-split exact structure to compare the nonconnective
$K$-theory of the enhanced relative singularity categories.
For the $K$-theory obtained from the Waldhausen category of
Gorenstein projective modules with stable isomorphisms as weak
equivalences, we refer to \cite{Ren2022}.

\subsection{Radical kernels and ordinary \texorpdfstring{$K$}{K}-groups}

For a finite-dimensional algebra $A$ over a field $k$, write
$J_A=\operatorname{rad}(A)$ and $S_A=A/J_A$. We use $K_i(A)$ for
the Quillen $K$-group of the exact category $\mathcal P(A)$ of finitely
generated projective left $A$-modules. Put
\[
 \widetilde K_i(A)=
 \ker\bigl(q_{A*}:K_i(A)\longrightarrow K_i(S_A)\bigr),
 \qquad i\geq0,
\]
where $q_A:A\to S_A$ is the quotient homomorphism. The exact functor
inducing $q_{A*}$ is $S_A\otimes_A-$. The following result identifies
a part of ordinary algebraic $K$-theory on which the projective
bimodule errors in a stable equivalence of Morita type act trivially.

\begin{theorem}\label{thm:radical-kernel-invariance}
Let $A$ and $B$ be finite-dimensional $k$-algebras. Suppose that
bimodules ${}_AM_B$ and ${}_BN_A$, finitely generated and projective
on both sides, satisfy
\[
 M\otimes_BN\simeq A\oplus P,
 \qquad
 N\otimes_AM\simeq B\oplus Q,
\]
where $P$ is projective over $A^e=A\otimes_kA^{\mathrm{op}}$ and
$Q$ is projective over $B^e$. Then the maps induced by
$F=N\otimes_A-$ and $G=M\otimes_B-$ restrict to mutually inverse
isomorphisms
\[
 \widetilde K_i(A)\simeq\widetilde K_i(B)
 \qquad(i\geq0).
\]
No separability hypothesis on $S_A$ or $S_B$, and no adjunction
hypothesis on the tensor functors, is required.
\end{theorem}

\begin{proof}
We first check that $F_*$ preserves the indicated kernels.
Set $V=S_B\otimes_BN$, viewed as an $S_B$-$A$-bimodule. Choose
$d$ with $J_A^d=0$, and put
\[
 U_r=VJ_A^r/VJ_A^{r+1}\qquad(0\leq r<d).
\]
Each $U_r$ is an $S_B$-$S_A$-bimodule. All the left $S_B$-modules
$VJ_A^r$ and $U_r$ are finitely generated projective, since $S_B$
is semisimple. For every $X\in\mathcal P(A)$, tensoring the
radical filtration of $V$ with $X$ is exact and has successive
quotients
\[
 U_r\otimes_AX\simeq
 U_r\otimes_{S_A}(S_A\otimes_AX).
\]
Thus this filtration is a finite sequence of objectwise admissible
short exact sequences of exact functors
$\mathcal P(A)\to\mathcal P(S_B)$. The additivity theorem for
algebraic $K$-theory \cite[Chapter~V, Theorem~1.2 and Corollary~1.2.1]{WeibelKbook}
therefore gives
\begin{equation}\label{eq:k-radical-descent}
 q_{B*}F_*
 =\left(\sum_{r=0}^{d-1}
      (U_r\otimes_{S_A}-)_*\right)q_{A*}.
\end{equation}
In particular, $F_*(\widetilde K_i(A))\subseteq\widetilde K_i(B)$.
Interchanging $A$ and $B$ proves the corresponding assertion for $G$.

We next show that the endofunctor $P\otimes_A-$ acts trivially
on $\widetilde K_i(A)$. Choose a bimodule $P'$ and an integer $t$
such that
\[
 P\oplus P'\simeq(A\otimes_kA)^t
\]
as $A$-$A$-bimodules, where the first and second factors carry
the left and right actions, respectively. Taking right radical
filtrations and successive quotients shows that $PJ_A^r$ and
\[
 W_r=PJ_A^r/PJ_A^{r+1}
\]
are direct summands, as left $A$-modules, of
$(A\otimes_kJ_A^r)^t$ and
$(A\otimes_k(J_A^r/J_A^{r+1}))^t$, respectively. Hence both are
finitely generated projective left $A$-modules. Moreover, $W_r$
is an $A$-$S_A$-bimodule. Consequently all the functors occurring
in the filtration of $P\otimes_A-$ take values in
$\mathcal P(A)$, and $W_r\otimes_{S_A}-$ is an exact functor
$\mathcal P(S_A)\to\mathcal P(A)$. A second application of
additivity yields
\begin{equation}\label{eq:k-error-factorization}
 (P\otimes_A-)_*
 =\left(\sum_{r=0}^{d-1}
       (W_r\otimes_{S_A}-)_*\right)q_{A*}.
\end{equation}
It follows that $(P\otimes_A-)_*$ vanishes on
$\widetilde K_i(A)$. The same argument applies to $Q$.

Finally, the given bimodule decompositions and additivity give
\[
 G_*F_* =\operatorname{id}_{K_i(A)}+(P\otimes_A-)_*,
 \qquad
 F_*G_* =\operatorname{id}_{K_i(B)}+(Q\otimes_B-)_*.
\]
Restricting these identities to the radical kernels proves the assertion.
\end{proof}

\begin{corollary}\label{cor:radical-relative-k}
Under the assumptions of Theorem~\ref{thm:radical-kernel-invariance},
suppose in addition that both quotient homomorphisms $q_A$ and $q_B$
admit unital ring sections. Then
\[
 K_i(A,J_A)\simeq K_i(B,J_B)\qquad(i\geq0),
\]
where the relative groups are defined by the homotopy fibres of the
maps on connective algebraic $K$-theory spectra. Under the canonical
identifications with the radical kernels, these isomorphisms are
the restrictions of $F_*$ and $G_*$. In particular, the conclusion
holds if $S_A$ and $S_B$ are separable over $k$.
\end{corollary}

\begin{proof}
A ring section of $q_A$ induces a section on algebraic $K$-theory.
The long exact sequence of the homotopy fibre therefore identifies
$K_i(A,J_A)$ canonically with $\ker(q_{A*})$ for every $i\geq0$.
The same applies to $B$, and the assertion follows from
Theorem~\ref{thm:radical-kernel-invariance}. If the semisimple
quotients are separable, the required sections exist by the
Wedderburn--Malcev theorem.
\end{proof}

\begin{corollary}\label{cor:ordinary-k-simple-count}
Let $A$ and $B$ be split finite-dimensional $k$-algebras which are
stably equivalent of Morita type, and write $r_A$ and $r_B$ for
their numbers of isomorphism classes of simple modules. Then the
following conditions are equivalent:
\begin{enumerate}
\item $r_A=r_B$;
\item $K_0(A)\simeq K_0(B)$ as groups;
\item $K_i(A)\simeq K_i(B)$ as groups for every $i\geq0$.
\end{enumerate}
The isomorphisms in the last condition need not be induced by
the tensor functors defining the stable equivalence.
\end{corollary}

\begin{proof}
Since the algebras are split, their semisimple quotients are
products of full matrix algebras over $k$. The Wedderburn--Malcev
theorem, Morita invariance of algebraic $K$-theory, and the splitting
of the quotient maps give
\[
 K_i(A)\simeq\widetilde K_i(A)\oplus K_i(k)^{r_A},
 \qquad
 K_i(B)\simeq\widetilde K_i(B)\oplus K_i(k)^{r_B}.
\]
Thus (1) implies (3) by Theorem~\ref{thm:radical-kernel-invariance}.
The implication (3)$\Rightarrow$(2) is immediate, whereas
(2)$\Rightarrow$(1) follows from
$K_0(A)\simeq\mathbb Z^{r_A}$ and
$K_0(B)\simeq\mathbb Z^{r_B}$.
\end{proof}

\begin{theorem}\label{thm:radical-spectra}
Under the hypotheses of Theorem~\ref{thm:radical-kernel-invariance},
assume that the quotient maps $q_A:A\to S_A$ and $q_B:B\to S_B$
admit unital ring sections. The homotopy fibres of
$\mathbb K(A)\to\mathbb K(S_A)$ and
$\mathbb K(B)\to\mathbb K(S_B)$ are equivalent spectra.
In the stable homotopy category the comparison is uniquely
characterized by restriction of the defining tensor maps.
If $A$ and $B$ are split, their ordinary nonconnective
$K$-theory spectra are equivalent if and only if
$r_A=r_B$.
\end{theorem}

\begin{proof}
Apply additivity at the spectrum level to the two filtrations
in the proof of Theorem~\ref{thm:radical-kernel-invariance}.
In the stable homotopy category this gives factorizations
\[
 q_{B*}F_*=u q_{A*},\qquad (P\otimes_A-)_*=v q_{A*},
\]
with $u$ and $v$ induced by the sums of the layer functors in
\eqref{eq:k-radical-descent} and \eqref{eq:k-error-factorization}.
Interchanging the two algebras gives the factorizations for
$G$ and $Q$. After passage to the perfect derived enhancements,
the same additivity of exact functors gives these identities
for nonconnective $K$-theory \cite{BGT2013}. The ring sections split the
maps $q_{A*}$ and $q_{B*}$ as maps of spectra. Write $I_A,I_B$
for the two fibres, with split inclusions $i_A,i_B$.
The first factorization gives a unique morphism
$f:I_A\to I_B$ in the stable homotopy category with
$i_Bf=F_*i_A$. The same construction gives $g:I_B\to I_A$
with $i_Ag=G_*i_B$. By the second factorization and
$GF\simeq\operatorname{Id}\oplus(P\otimes_A-)$, one has
$i_Agf=i_A$. Since $i_A$ is split, $gf=1$. The symmetric
argument gives $fg=1$.

For a split algebra the complement of $I_A$ in $\mathbb K(A)$ is
$\mathbb K(k)^{r_A}$ by Morita invariance and the product
decomposition of $S_A$. This proves sufficiency of equal simple
counts. Necessity follows by taking $\pi_0$.
\end{proof}

\begin{remark}\label{rem:k-theory-scope}
The preceding results do not require the algebras to have no
separable direct factors. Excluding such factors does not, by
itself, establish equality of the numbers of simple modules.
Moreover, the equivalence in Corollary~\ref{cor:ordinary-k-simple-count}
uses degree zero; the argument does not give the converse if only
the positive $K$-groups are assumed to be isomorphic. Over an
arbitrary field, equal numbers of simple modules also do not
identify the $K$-theories of the possibly different division rings
occurring in the semisimple quotients. Finally, in the absence of
sections of $q_A$ and $q_B$, the radical kernels must not be
identified with the full relative $K$-groups without an additional
argument controlling the connecting homomorphisms.
\end{remark}

The absence of separable direct factors does not supply the
simple-count hypothesis in Corollary~\ref{cor:ordinary-k-simple-count}.
Over an algebraically closed field it excludes semisimple blocks;
already for self-injective algebras in this class, the remaining
counting assertion is the stable-equivalence Auslander--Reiten
problem discussed in \cite{LiuZhouZimmermann2012,LiXi2026Centralizer}.
This is the conjecture about nonprojective simple modules, not the
homological conjecture about vanishing self-extensions.
For Frobenius-finite algebras over an algebraically closed field,
with no semisimple direct factors, ordinary $K$-theory invariance
also follows from the derived lifting theorem
\cite[Theorem~1.1]{HuXi2018II}. That special case is an established
consequence of derived equivalence.

\subsection{Stable cyclic homology over a commutative base}
\label{subsec:relative-stable-cyclic}

Let $R$ be a commutative ring. In this and the next subsection
we use right modules. Write $\mathcal Q_R(D)$ for the $R$-linear
category of finite relative projectives
$\operatorname{add}\{V\otimes_R D\}$, with $V$ a finite
$R$-module, viewed as a dg category concentrated in degree zero.
Algebras and bimodule coefficients in these two subsections
are finite over $R$. All bimodules are central over $R$.
An induced $D$-$L$-bimodule means one of the form
$D\otimes_R V\otimes_R L$. A relatively projective bimodule
is a summand of a finite sum of such bimodules.

We use the derived mixed complex over $R$, denoted by
$\mathcal C_R$, as a dg Morita invariant. It can be evaluated
using an $R$-cofibrant dg replacement. It is additive on a
direct sum of functors and on objectwise split exact
sequences of functors between additive categories, by
additivity for the mixed complex of a split exact category
\cite{KellerCyclicExact,TabuadaGaloisAdditive}.
The ordinary cyclic, negative cyclic and periodic cyclic
realizations of a mixed complex are denoted by $HC$, $HC^-$
and $HP$. We use their derived functors on the derived
category of mixed complexes.

For $H$ one of these realizations in degree $n$, put
\[
 H^{\mathrm{st}}_{n,R}(D)=
 \bigcap_{\substack{L,\;P\\
      P\ \text{relatively projective }D\text{-}L}}
 \ker\!\left(
 H_n(\mathcal C_R(\mathcal Q_R(D)))
   \xrightarrow{\,P_*\,}
 H_n(\mathcal C_R(\mathcal Q_R(L)))\right).
\]
Here $L$ ranges over the finite central $R$-algebras, and
$P_*$ is induced by $-\otimes_D P$. The target functor is
well defined: its value on any finite module is a relative
projective when $P$ is an induced bimodule, and the same
holds on its summands. One can take representatives of
the tests in a fixed universe. This definition uses
relatively projective bimodules as transfer tests.

\begin{proposition}\label{prop:relative-stable-cyclic}
Suppose that $M$ and $N$ are sidewise finite projective and
that $M\otimes_BN$ and $N\otimes_AM$ are stably the regular
bimodules modulo $R$-relative projectives. The two tensor
maps induce inverse isomorphisms on $H^{\mathrm{st}}_{n,R}$
for $H=HC,HC^-,HP$. Over an algebraically closed field,
the construction for $HC$ agrees with stable cyclic
homology in \cite[Definition~9.7]{KoenigLiuZhou2012}.
\end{proposition}

\begin{proof}
The functors $-\otimes_A M$ and $-\otimes_B N$ preserve
relative projectives. If $P$ is a relatively projective
$B$-$L$-bimodule, then $M\otimes_B P$ is a relatively projective
$A$-$L$-bimodule, since $M$ is projective on the left. Thus
composition with a test on the target is a test on the
source, and the tensor maps preserve the intersections of
kernels. In a stabilized identity the errors belong to the
tests with target $L=A$ or $L=B$. Their maps vanish on the
indicated subgroups. Additivity and tensor composition now
show that the restricted maps are inverse.

Over an algebraically closed field, $\mathcal Q_k(D)$ is the
category of finite projective modules. The tests with $L=k$
are the transfers with finite projective left $D$-modules.
To identify the maps, a right $D$-module is tensored with
a left projective $P$. The resulting action on $P$ is the map
from $D$ to $\operatorname{End}_k(P)$. Passing through the matrix
trace on cyclic chains inserts a projective dual basis, which
is the formula in \cite[Section~9]{KoenigLiuZhou2012}.
Every projective $D$-$L$-bimodule is a direct sum
of external products $De\otimes_k fL$ of indecomposable
projectives. Indeed the images of primitive tensor idempotents
are split primitive idempotents in the quotient by the sum
of the two radical ideals. This sum is nilpotent, so these
idempotents give all the projective types. The transfer for
each external product factors through a test with target $k$.
Additivity therefore identifies the two intersections of
kernels. The field result and the role of the odd and even
degrees of ordinary cyclic homology are described in
\cite[Theorems~9.6 and~9.9 and Remark~9.10]{KoenigLiuZhou2012}.
\end{proof}

\subsection{Additive invariants of radical fibres}
\label{subsec:radical-additive}

We use the universal dg additive invariant $U_R$ of
\cite[Section~2]{TabuadaGaloisAdditive}. Its morphisms on
dg categories can be represented by Grothendieck classes
of bimodule correspondences. Equivalently for the following
arguments, apply any dg Morita-invariant functor $E$ with
values in an additive category that satisfies additivity
for objectwise split exact sequences of exact functors.
Indeed such a sequence gives an extension triangle on the
associated representable bimodules, and hence the sum
relation in the correspondence group. For the ordinary
case below write $\mathcal Q(A)$ for finite projective right
$A$-modules.

\begin{theorem}\label{thm:radical-additive}
Suppose $A,B,M,N$ satisfy the hypotheses of
Theorem~\ref{thm:radical-kernel-invariance}, with the same
notation for the semisimple quotients. For every dg additive
invariant with values in an abelian category,
\[
 \ker(E(\mathcal Q(A))\to E(\mathcal Q(S_A)))
 \ \simeq\
 \ker(E(\mathcal Q(B))\to E(\mathcal Q(S_B))).
\]
The isomorphisms are induced by the tensor maps on the
kernels. If the quotient maps admit algebra sections, the
tensor maps similarly identify the two complementary
summands for any idempotent-complete additive target.
In particular they identify the homotopy fibres when the
target is stable. At the level of $U_k$ this is an
isomorphism of summands in the idempotent completion of the
additive correspondence category.
\end{theorem}

\begin{proof}
Use right modules in the filtration argument of
Theorem~\ref{thm:radical-kernel-invariance}. The bimodule
$M\otimes_B S_B$ has the filtration by the powers of $J_A$
from the left. The layers are projective on the $S_B$
side and have left action through $S_A$. The filtration
induces objectwise split exact sequences on projectives
and hence $E(q_B)E(M)=uE(q_A)$ for a sum $u$ of layer
maps. The corresponding equation holds in the reverse
direction. The filtration of an error bimodule $P$ by
$J_A^iP$ has layers projective on the right over $A$,
since $P$ is a summand of a finite sum of $A\otimes_k A$.
Its action therefore factors through $E(q_A)$ by the same
additivity. The two composites of the tensor maps act
as the identities on the kernels.

If sections are given, let $i_A:F_A\to E(\mathcal Q(A))$
denote the inclusion of the complementary summand of the
section, with a similar notation for $B$. Then the factorization
through $E(q_A)$ implies $E(q_B)E(M)i_A=0$. Since $i_B$
is a split monomorphism and identifies the kernel of the
split quotient in the additive category, $E(M)i_A$ factors
uniquely through $i_B$. Both error maps annihilate these
summands, so the two restrictions are inverse. In a stable
target the complement is the fibre of the split quotient.
The restrictions depend only on the tensor maps and fibre
inclusions in the homotopy category.
\end{proof}

\begin{corollary}\label{cor:radical-mixed-complex}
Under the section hypothesis of Theorem~\ref{thm:radical-additive},
the fibres of $\mathcal C_k(A)\to\mathcal C_k(S_A)$ and
$\mathcal C_k(B)\to\mathcal C_k(S_B)$ are isomorphic in the
derived category of mixed complexes. Hence their Hochschild
and their ordinary, negative and periodic cyclic homologies
are isomorphic. This applies over a perfect field.
\end{corollary}

\begin{proof}
The mixed complex is dg additive and Morita invariant.
Apply Theorem~\ref{thm:radical-additive} to it. The section
condition over a perfect field is the Wedderburn--Malcev
splitting theorem. This fibre construction uses the quotient
by the radical. Proposition~\ref{prop:relative-stable-cyclic}
uses an intersection of transfer kernels.
\end{proof}

For instance, let $D=k[\varepsilon]/(\varepsilon^2)$ over an
algebraically closed field of characteristic $2$. Every finite
projective left $D$-module is free. The degree-zero transfer
of a free module is the trace of the regular action repeated
once per free generator, hence is zero. Thus $HC_0^{\mathrm{st}}(D)=D$,
whereas the degree-zero homology of the radical fibre is
$k\varepsilon$.

\begin{theorem}\label{thm:relative-nilpotent-additive}
Let $A,B$ be finite central algebras over an arbitrary
commutative ring $R$. Let $M,N$ be as in
Proposition~\ref{prop:relative-stable-cyclic}. Suppose there
are nilpotent ideals $I_D$ in $D=A,B$ with quotients
$S_D=D/I_D$ separable over $R$. Assume the following
filtrations split over $R$, in the sense that each successive
quotient sequence splits:
\[
 (I_A^j)_j,\quad (I_B^j)_j,\qquad
 (I_A^j(M\otimes_B S_B))_j,\quad
 (I_B^j(N\otimes_A S_A))_j.
\]
For every dg additive invariant, the corresponding
comparisons of Theorem~\ref{thm:radical-additive} hold with
$\mathcal Q_R$ in place of $\mathcal Q$ and these quotients
in place of the semisimple quotients. The conclusions about
complements and fibres use algebra sections as there.
\end{theorem}

\begin{proof}
Quotient extension of scalars preserves relative projectives.
In the two filtrations of $M\otimes_B S_B$ and
$N\otimes_A S_A$, the layers are relatively projective on
the target side. In fact, for a separable $R$-algebra
$S$, the separability idempotent splits
$Y\otimes_R S\to Y$ as a right $S$-module map for every
right $S$-module $Y$. Each filtration sequence stays
$R$-split after tensoring with $V\otimes_R A$ and hence
with any relative projective over $A$. This follows by
tensoring the given $R$-splitting with $V$ and cutting
it by the idempotent of the relative projective. Its
quotient is relatively projective over $S_B$, so the
resulting sequence is split over $S_B$. The corresponding
argument applies in the other direction.

For an induced error $D\otimes_R V\otimes_R D$ its left
filtration is $I_D^j\otimes_R V\otimes_R D$, as follows
from the assumed splitting of the filtration of $D$.
Its layers are induced on the right. For a direct
summand of the error, the defining idempotent preserves
the filtration and restricts the layer splittings and
relative projectivity. As above, tensoring these sequences
with a relative projective yields split sequences of
right $D$-modules. Therefore the layer maps factor the
error action through $\mathcal Q_R(S_D)$ in the additive
correspondence category. Apply the kernel and split
complement argument of Theorem~\ref{thm:radical-additive}.
\end{proof}

\subsection{The maps induced by the defining bimodules}

The isomorphisms of ordinary groups in
Corollary~\ref{cor:ordinary-k-simple-count} are abstract group
isomorphisms. The original tensor maps can fail to be isomorphisms
even for a reduced pair and an algebra without separable direct
factors.

\begin{example}\label{ex:k-tensor-syzygy}
Let $k$ be a field and put
\[
 A=k[x]/(x^3),\qquad E=A^e=k[x,y]/(x^3,y^3),\qquad
 L=\ker(E\xrightarrow{\mu}A)=(x-y).
\]
Then $(L,L)$ defines a stable equivalence of Morita type from $A$
to itself. The bimodule $L$ is reduced, and $A$ has no nonzero
separable algebra direct factor. Nevertheless the tensor functor
$F=L\otimes_A-$ induces multiplication by $2$ on $K_0(A)=\mathbb Z$.
If $\operatorname{char}k\ne2$, its map on $K_1(A)$ is not injective.
\end{example}

\begin{proof}
The sequence $0\to L\to E\xrightarrow{\mu}A\to0$ splits on each
side. Since $A$ is local and $E$ is free of rank $3$ on each side,
$L$ is free of rank $2$ on each side. Put $v=x^2+xy+y^2$.
The relations $(x-y)v=0$ and $x^3=y^3=0$ show directly that
\[
 0\longrightarrow A\xrightarrow{\,1\mapsto v\,}E
   \xrightarrow{\,x-y\,}L\longrightarrow0
\]
is an exact sequence of bimodules. For example, $v,xv,x^2v$
are linearly independent and span the kernel of multiplication
by $x-y$. Tensoring $0\to L\to E\to A\to0$ with $L$ gives
\[
 0\longrightarrow L\otimes_A L\longrightarrow E^2
   \longrightarrow L\longrightarrow0.
\]
Schanuel's lemma and Krull--Schmidt cancellation over $E$ give
\begin{equation}\label{eq:k-syzygy-square}
                  L\otimes_A L\simeq A\oplus E.
\end{equation}
This proves the Morita-type identities. Since $E$ is local and
$\dim_kL=6<9=\dim_kE$, the bimodule $L$ has no nonzero
$E$-projective direct summand. The algebra $A$ is local and
nonsemisimple, so it has no nonzero separable direct factor.

The equality $[F(A)]=[L]=2[A]$ proves the assertion about $K_0$.
For $K_1$, the commutative local ring $A$ has
$K_1(A)=A^\times$. Write $u=a+bx+cx^2$ with $a\ne0$.
On the free left $A$-module $E$, the map obtained by tensoring
multiplication by $u$ is multiplication by $u$ on the second
factor, and has determinant $a^3$. Additivity applied to
$0\to L\to E\to A\to0$ therefore gives, in multiplicative
notation,
\begin{equation}\label{eq:k-syzygy-units}
                         F_*(u)=a^3u^{-1}.
\end{equation}
In characteristic different from $2$, the nontrivial unit $-1$
is sent to $1$. Thus $F_*$ on $K_1$ is not injective.
\end{proof}

This example does not contradict an abstract comparison of the
ordinary groups: the two algebras are identical. It separates
such a comparison from an isomorphism induced by a prescribed
tensor functor. In contrast, the restrictions in
Theorem~\ref{thm:radical-kernel-invariance} are induced by the
original tensor functors and are mutually inverse.

\subsection{Enhanced relative singularity categories and algebraic \texorpdfstring{$K$}{K}-theory}
\label{subsec:enhanced-relative-singular-k}

In this subsection, $R$ is an arbitrary commutative ring and $A$ is an
$R$-algebra which is finitely generated as an $R$-module. All bimodules
are $R$-central. Let $\mathcal E_R(A)$ be the exact category of finitely
generated left $A$-modules, with the short exact sequences which split
after restriction to $R$ as its conflations. Finite modules are
closed under extensions and direct summands, so this is an exact
category without any Noetherian hypothesis. Its projective objects
are precisely
\[
 \mathcal P_R(A)=\operatorname{add}
 \{A\otimes_R U\mid U\in R\text{-}\mathrm{mod}\}.
\]
Indeed, the induced modules are projective for this exact structure,
and the multiplication map $A\otimes_R X\longrightarrow X$ is an
$R$-split deflation for every $X\in A\text{-}\mathrm{mod}$.

Write
\[
 \mathcal D_R(A)=D^b_\infty(\mathcal E_R(A))^{\natural},\qquad
 \operatorname{Perf}_R(A)
 =\operatorname{thick}_{\mathcal D_R(A)}\mathcal P_R(A),
\]
where $(-)^{\natural}$ denotes idempotent completion. We use the
bounded derived stable $\infty$-category of an exact category; equivalently,
one may use its stable hull \cite[Theorem~1.2 and equation~(1)]{Klemenc2022}.
Define
\[
 \mathcal S_R(A)=
 \bigl(\mathcal D_R(A)/\operatorname{Perf}_R(A)\bigr)^{\natural},
 \qquad
 \mathbb K^{\mathrm{sg}}_R(A)=\mathbb K(\mathcal S_R(A)).
\]
We write $K^{\mathrm{sg}}_{i,R}(A)=\pi_i\mathbb K_R^{\mathrm{sg}}(A)$
for every $i\in\mathbb Z$. All exact categories are replaced by
small skeletons when forming these constructions.
The quotient is formed in stable $\infty$-categories and
$\mathbb K$ denotes nonconnective algebraic $K$-theory. This definition
retains the enhancement and does not require $A$ to be flat over $R$.
The functorial enhancement is specified by the exact category
$\mathcal E_R(A)$, rather than chosen from its triangulated quotient.

\begin{lemma}\label{lem:relative-error-exact}
Let $P$ be a direct summand, as an $A$--$A$-bimodule, of a finite direct
sum of bimodules $A\otimes_R U\otimes_R A$, where
$U\in R\text{-}\mathrm{mod}$. Then
\[
 P\otimes_A-:\mathcal E_R(A)\longrightarrow\mathcal P_R(A)
\]
is an exact functor when $\mathcal P_R(A)$ carries its split exact
structure. No one-sided projectivity assumption on $P$ is needed.
\end{lemma}

\begin{proof}
For an induced bimodule the functor in question is
$A\otimes_R U\otimes_R-$. Its value on $X$ is an induced
$A$-module. If a sequence splits over $R$, tensoring its splitting
with $A\otimes_R U$ gives an $A$-linear splitting. Passing to finite
direct sums and direct summands proves the assertion. Explicitly,
if $e$ is the natural idempotent defining the summand and $s$ is a
splitting for the induced sequence, then the restriction of $ese$
is a splitting for the summand sequence.
\end{proof}

\begin{theorem}\label{thm:relative-singular-k-invariance}
Let $R$ be a commutative ring and let $A,B$ be finite central $R$-algebras. Suppose that
${}_AM_B$ and ${}_BN_A$ are finitely generated projective on each side,
and that there are bimodule isomorphisms
\[
 (M\otimes_BN)\oplus P_1\simeq A\oplus P_2,
 \qquad
 (N\otimes_AM)\oplus Q_1\simeq B\oplus Q_2,
\]
where each $P_j$ is a direct summand of a finite direct sum of
$A\otimes_R U\otimes_R A$, and each $Q_j$ is a direct summand of a
finite direct sum of $B\otimes_R V\otimes_R B$, for finitely generated
$R$-modules $U,V$. Then $N\otimes_A-$ and $M\otimes_B-$ induce mutually
inverse equivalences
\[
 \mathcal S_R(A)\simeq\mathcal S_R(B).
\]
Consequently they induce an equivalence of spectra
\[
 \mathbb K_R^{\mathrm{sg}}(A)\simeq\mathbb K_R^{\mathrm{sg}}(B)
\]
and isomorphisms
\[
 K^{\mathrm{sg}}_{i,R}(A)\simeq K^{\mathrm{sg}}_{i,R}(B)
 \qquad(i\in\mathbb Z).
\]
The assertion includes the case of actual tensor decompositions and
the case where the error bimodules are projective over $A^e$ and $B^e$.
\end{theorem}

\begin{proof}
Put $F=N\otimes_A-$ and $G=M\otimes_B-$. Since $N_A$ is finitely
generated projective, the underlying $R$-module functor of $F$ is a
direct summand of a finite direct sum of the restriction functor
$A\text{-}\mathrm{mod}\longrightarrow R\text{-}\mathrm{mod}$.
It therefore sends $R$-split conflations to $R$-split conflations.
The same argument applies to $G$.

For $U\in R\text{-}\mathrm{mod}$ there is an isomorphism
\[
 F(A\otimes_R U)\simeq N\otimes_R U.
\]
Since ${}_BN$ is finitely generated projective, its right-hand side
is a direct summand of a finite direct sum of $B\otimes_R U$.
Thus $F$ sends $\mathcal P_R(A)$ into $\mathcal P_R(B)$, and similarly
$G$ sends $\mathcal P_R(B)$ into $\mathcal P_R(A)$. Exact functoriality
of bounded derived stable $\infty$-categories, followed by idempotent
completion, gives exact functors between $\mathcal D_R(A)$ and
$\mathcal D_R(B)$ preserving the relative perfect subcategories.
Hence these functors descend to $\mathcal S_R(A)$ and
$\mathcal S_R(B)$.

By Lemma~\ref{lem:relative-error-exact}, each error functor factors
through the split exact category of relatively projective modules.
On bounded complexes its image therefore lies in
$\operatorname{Perf}_R(A)$ or $\operatorname{Perf}_R(B)$, respectively.
The given bimodule isomorphisms are natural isomorphisms of exact
functors
\[
 GF\oplus(P_1\otimes_A-)
 \simeq\operatorname{Id}\oplus(P_2\otimes_A-),
\]
and
\[
 FG\oplus(Q_1\otimes_B-)
 \simeq\operatorname{Id}\oplus(Q_2\otimes_B-).
\]
They extend functorially to the derived enhancements. All four error
functors become zero after passage to the enhanced quotients.
Consequently the induced composites are naturally equivalent to the
identity functors. This proves the asserted equivalence. Applying
nonconnective algebraic $K$-theory gives the equivalence of spectra
and the isomorphisms on their homotopy groups.
\end{proof}

\begin{remark}\label{rem:relative-singular-k-localization}
The definition gives a localization sequence
\[
 \mathbb K(\operatorname{Perf}_R(A))
 \longrightarrow \mathbb K(\mathcal D_R(A))
 \longrightarrow \mathbb K_R^{\mathrm{sg}}(A)
\]
by the localization theorem for nonconnective algebraic $K$-theory
\cite[Section~5 and Theorem~9.8]{BGT2013}. No Frobenius, Gorenstein,
adjointness, or base-flatness hypothesis is used. In particular, the
proof does not assume that the ordinary tensor functors give an
equivalence of the ordinary stable module categories.

When $R=k$ is a field, every short exact sequence is $k$-split and
$\mathcal P_k(A)=A\text{-}\mathrm{proj}$. Hence the construction is
the usual enhanced singularity category, with idempotent completion,
and the displayed sequence is
\[
 \mathbb K(A)\longrightarrow\mathbb G(A)
 \longrightarrow\mathbb K^{\mathrm{sg}}(A).
\]
For a general base ring the relative perfect objects need not be
ordinary perfect complexes, so this sequence must not be identified
with the ordinary $K$--$G$ localization sequence. Likewise, an
identification of $K^{\mathrm{sg}}_{0,R}(A)$ with a naive cokernel
requires a separate examination of idempotent completion and the
negative $K$-group boundary.
\end{remark}

\subsection{Specialization to fields and derived constructions}

\begin{corollary}\label{cor:field-singular-k}
Let $A,B$ be finite-dimensional algebras over a field $k$ which
are stably equivalent of Morita type. The defining tensor functors
induce inverse equivalences of their enhanced singularity categories
and an equivalence
\[
                  \mathbb K^{\mathrm{sg}}(A)
                    \simeq\mathbb K^{\mathrm{sg}}(B).
\]
In degree zero, this recovers the stable Grothendieck-group
invariant of the Cartan matrix:
\[
 K^{\mathrm{sg}}_0(A)\simeq\operatorname{coker}
       \bigl(K_0(A)\longrightarrow G_0(A)\bigr)
       \simeq\operatorname{coker}C_A.
\]
\end{corollary}

\begin{proof}
Apply Theorem~\ref{thm:relative-singular-k-invariance} with $R=k$.
The finite-dimensional singularity category is idempotent-complete;
see \cite[Section~2]{ChenLiuWang2023}. Its Grothendieck group is
therefore the cokernel of the map from perfect complexes to the
bounded derived category. Under the projective and simple bases
this map is the Cartan matrix, which proves the last assertion.
\end{proof}

The enhancement assertion over a field is consistent with the
existing tensor-induced dg comparison of
\cite[Proposition~3.1 and Lemma~3.2]{ChenLiuWang2023}: a tensor
functor defined by a bimodule projective on each side and inducing
a singular equivalence has a dg lift which is a quasi-equivalence. The field case is a consequence of this
established comparison and the bimodule identities. The proof of
Theorem~\ref{thm:relative-singular-k-invariance} supplies the relative
construction over an arbitrary commutative base.

\begin{corollary}\label{cor:derived-source-relative-k}
The bimodule pairs constructed in
Theorem~\ref{thm:rz-relative-derived-stable} and
Theorem~\ref{anu-thm-raw-semt} induce equivalences
\[
                  \mathbb K_R^{\mathrm{sg}}(A)
                    \simeq\mathbb K_R^{\mathrm{sg}}(B).
\]
For the Rouquier--Zimmermann construction, the stabilized
bimodule identities are sufficient. No realization on the entire
ordinary stable module categories is required for this conclusion.
\end{corollary}

\begin{proof}
The cited constructions provide bimodules projective on each side
and actual or stabilized tensor identities with enveloping-projective
errors. These are special cases of the hypotheses of
Theorem~\ref{thm:relative-singular-k-invariance}.
\end{proof}

The relative result thus applies directly to the bimodule
identities furnished by the derived constructions. Its proof uses
the $R$-split exact structure even when this exact category
is not Frobenius. The ordinary radical kernels, the
relative groups for the ideal $J_A$, and the enhanced $R$-relative
singularity groups are distinct constructions; each comparison
above retains the hypotheses required by its definition.

\section{Adjoint normalization}
\label{sec:adjoint-normalization}

In this section we construct biadjoint representatives of
stable bimodule pairs. Until Subsection~\ref{subsec:adjoint-complete-base},
we work over a commutative Artin ring $R$, unless otherwise
stated. In that subsection we pass to a Henselian Noetherian
local base. Throughout, the algebras are finite over $R$, and
the bimodules are finitely generated and $R$-central. Put
\[
 C_{A,B}=A\otimes_R B^{\mathrm{op}},
 \qquad A^e_R=C_{A,A}.
\]
A projective $A$-$B$-bimodule means a projective $C_{A,B}$-module.
For a finite $R$-algebra $C$, consider the two classes
\[
 \mathcal P^{\mathrm{ord}}(C)=C\text{-}\mathrm{proj},
 \qquad
 \mathcal P^{\mathrm{rel}}(C)
 =\operatorname{add}\{C\otimes_RU\mid U\in R\text{-}\mathrm{mod}\}.
\]
The second class consists of the projective objects in the exact structure
whose conflations split over $R$.  We fix one of these two classes
throughout and write $\mathcal P_{A,B}$ for the chosen class over $C_{A,B}$.
A bimodule is called \emph{reduced} with respect to this class if it has no
nonzero direct summand in $\mathcal P_{A,B}$.

For ordinary projective error bimodules, the existence of adjoint-type
representatives for stable equivalences of Morita type between Artin
algebras without separable direct factors is stated in
\cite[Lemma~3.11]{XZ2026}.  The corresponding statement for
finite-dimensional algebras over arbitrary fields without semisimple
direct factors appears in
\cite[Lemma~4.1]{CPX2012}.  We construct both adjunctions directly from
the decompositions of a given reduced bimodule pair.  The proof applies
to either projectivity class and does not require projectivity over $R$.

\begin{lemma}
\label{lem:projective-factorization-radical}
Let $C$ be a finite $R$-algebra and let $\mathcal S$ be an additive
class of finitely generated $C$-modules closed under direct summands.
Suppose that $X$ is a finitely generated $C$-module with no nonzero direct
summand in $\mathcal S$.  Denote by
$\mathcal I_{\mathcal S}(X,X)$ the ideal in
$\operatorname{End}_C(X)$ of maps factoring through objects of
$\mathcal S$.  Then
\[
 \mathcal I_{\mathcal S}(X,X)
 \subseteq \operatorname{rad}\operatorname{End}_C(X).
\]
Consequently, every endomorphism of $X$ which is invertible modulo
$\mathcal I_{\mathcal S}(X,X)$ is an automorphism of $X$.
\end{lemma}

\begin{proof}
Let $f\in\mathcal I_{\mathcal S}(X,X)$ and
$g\in\operatorname{End}_C(X)$.  Then $gf$ factors through an object
of $\mathcal S$.  The Fitting decomposition gives, for sufficiently large $n$,
\[
 X=\ker(gf)^n\oplus\operatorname{im}(gf)^n,
\]
and $gf$ restricts to an automorphism of
$Y=\operatorname{im}(gf)^n$.  Since this automorphism factors through
an object of $\mathcal S$, so does $\operatorname{id}_Y$.  Thus $Y$ is a direct summand
of an object of $\mathcal S$, so $Y\in\mathcal S$.
The assumption on $X$ gives $Y=0$, so $gf$ is nilpotent.  Hence
$1-gf$ is invertible for every $g$, which proves the radical inclusion.
The last assertion follows because units lift modulo any ideal
contained in the Jacobson radical.
\end{proof}

\begin{lemma}
\label{lem:projective-cross-tensors}
Let $A,B,C$ be finite central $R$-algebras, and let ${}_AU_B$ and
${}_BV_C$ be bimodules.
\begin{enumerate}
\item If $U\in\mathcal P_{A,B}$ and $V_C$ is projective,
then $U\otimes_BV\in\mathcal P_{A,C}$.
\item If $V\in\mathcal P_{B,C}$ and ${}_AU$ is projective,
then $U\otimes_BV\in\mathcal P_{A,C}$.
\end{enumerate}
\end{lemma}

\begin{proof}
In the relative case, $U$ is a summand of
$A\otimes_RW\otimes_RB$ for some finite $R$-module $W$.
Tensoring the splitting and using $V_C\in\operatorname{add}(C_C)$
gives
\[
 U\otimes_BV\in\operatorname{add}(A\otimes_RW\otimes_RV)
 \subseteq\operatorname{add}(A\otimes_RW\otimes_RC).
\]
This proves (1) in the relative case.  For the ordinary case, take $W$
to be finite free over $R$.  For (2), use a split presentation of $V$ by
$B\otimes_RW\otimes_RC$ and
${}_AU\in\operatorname{add}({}_AA)$.  Since all the maps being tensored
are split, no flatness assumption over $R$ is needed.
\end{proof}

An algebra has no \emph{$R$-separable direct factor} if
its regular bimodule has no nonzero projective direct summand.
Equivalently, it has no nonzero direct algebra factor $A_1$ that is
separable over $R$, where separability means that the multiplication
$A_1\otimes_RA_1\to A_1$ splits as a bimodule map.  Indeed, the
idempotents of $\operatorname{End}_{A^e_R}(A)=Z(A)$ are precisely
the central idempotents determining the direct algebra factors.

\begin{lemma}
\label{lem:unit-relative-projective-center}
Let $R$ be any commutative ring and $A$ any unital central
$R$-algebra.  For the regular bimodule $A$, the ideals of endomorphisms factoring
through ordinary and through $R$-relative projective $A^e_R$-modules
coincide.  In particular, these two classes have the same direct
summands of $A$, and the relative projective center is the ordinary
projective center $Z_{\mathrm{pr}}(A)$.
\end{lemma}

\begin{proof}
The multiplication $\mu:A\otimes_RA\to A$ is an $R$-split
epimorphism, with $R$-linear section $a\mapsto1\otimes a$.
If $A\to U\xrightarrow{h}A$ is a factorization through an
$R$-relative projective module, lift $h$ across $\mu$.
The resulting factorization is through the ordinary free
$A^e_R$-module $A\otimes_RA$.  The opposite inclusion of ideals
is immediate.  The second assertion follows by applying this equality
to the idempotent associated with a direct summand of $A$.
\end{proof}

\begin{lemma}
\label{lem:remove-projective-bimodule-summands}
Let $A$ and $B$ have no $R$-separable direct factors.  Suppose that
${}_AM_B$ and ${}_BN_A$ are projective on both sides and
satisfy
\begin{equation}
\label{eq:adjoint-morita-identities}
 M\otimes_BN\cong A\oplus P,
 \qquad N\otimes_AM\cong B\oplus Q,
\end{equation}
where $P\in\mathcal P_{A,A}$ and $Q\in\mathcal P_{B,B}$.  Write
\[
 M=M_0\oplus U,\qquad N=N_0\oplus V,
\]
where $M_0,N_0$ are reduced with respect to the chosen class, and
$U\in\mathcal P_{A,B}$, $V\in\mathcal P_{B,A}$.
Then $M_0,N_0$ are projective on both sides and satisfy
\[
 M_0\otimes_BN_0\cong A\oplus P_0,
 \qquad N_0\otimes_AM_0\cong B\oplus Q_0
\]
with $P_0\in\mathcal P_{A,A}$ and $Q_0\in\mathcal P_{B,B}$.
\end{lemma}

\begin{proof}
All four summands are projective on each side.  Expanding
$(M_0\oplus U)\otimes_B(N_0\oplus V)$ and applying
Lemma~\ref{lem:projective-cross-tensors}, we see that every term except
$M_0\otimes_BN_0$ belongs to $\mathcal P_{A,A}$.  Thus
\[
 (M_0\otimes_BN_0)\oplus T\cong A\oplus P
\]
for $T\in\mathcal P_{A,A}$.  By
Lemma~\ref{lem:unit-relative-projective-center}, $A$ has no nonzero
summands in $\mathcal P_{A,A}$.  The Krull--Schmidt theorem therefore
identifies the summands of $M_0\otimes_BN_0$ outside this class with
the summands of $A$.
The remaining summands form $P_0$.  The same argument over $B^e_R$
gives the second decomposition.
\end{proof}

We can now construct the two adjunctions for a reduced bimodule pair.

\begin{theorem}
\label{thm:adjoint-normalization}
Let $A,B$ be finite central algebras over a commutative Artin ring
$R$, and fix one of the two projectivity classes defined above.
Suppose that ${}_AM_B$ and ${}_BN_A$ are reduced with
respect to the chosen class, projective as ordinary modules on both
sides, and satisfy \eqref{eq:adjoint-morita-identities}, with
$P\in\mathcal P_{A,A}$ and $Q\in\mathcal P_{B,B}$.  Then
\[
 M\otimes_B-\ \dashv\ N\otimes_A-
 \qquad\text{and}\qquad
 N\otimes_A-\ \dashv\ M\otimes_B-
\]
on the ordinary module categories.  In particular, there are bimodule
isomorphisms
\begin{align*}
 {}_BN_A
 &\cong \operatorname{Hom}_A(M,A)
 \cong \operatorname{Hom}_{B^{\mathrm{op}}}(M,B),\\
 {}_AM_B
 &\cong \operatorname{Hom}_B(N,B)
 \cong \operatorname{Hom}_{A^{\mathrm{op}}}(N,A).
\end{align*}
\end{theorem}

\begin{proof}
Fix the decompositions in \eqref{eq:adjoint-morita-identities},
with inclusions and projections
\[
 i:A\longrightarrow M\otimes_BN,
 \quad p:M\otimes_BN\longrightarrow A,
 \qquad
 v:B\longrightarrow N\otimes_AM,
 \quad q:N\otimes_AM\longrightarrow B.
\]
Thus $pi=1_A$, $qv=1_B$, $\ker p\cong P$, and
$\operatorname{coker}v\cong Q$.  Associativity and unit
identifications for tensor products will be understood.

Consider the $A$-$B$-bimodule endomorphism
\[
 t_M=(p\otimes_A1_M)(1_M\otimes_Bv):
 M\longrightarrow M\otimes_BN\otimes_AM\longrightarrow M.
\]
The first map $1_M\otimes_Bv$ is a split monomorphism with
complement $M\otimes_BQ$, and the second map
$p\otimes_A1_M$ is a split epimorphism with kernel
$P\otimes_AM$.  Both complements belong to $\mathcal P_{A,B}$
by Lemma~\ref{lem:projective-cross-tensors}.  Each map is therefore
invertible in the quotient of $C_{A,B}$-modules by maps through
objects of $\mathcal P_{A,B}$.  Consequently, the image of $t_M$ in
\[
 \operatorname{End}_{C_{A,B}}(M)/\mathcal I_{\mathcal P_{A,B}}(M,M)
\]
is a unit.  As $M$ is reduced,
Lemma~\ref{lem:projective-factorization-radical} shows that $t_M$
is a bimodule automorphism.

Define
\begin{equation}
\label{eq:normalized-counit}
 p'=p(t_M^{-1}\otimes_B1_N):M\otimes_BN\longrightarrow A.
\end{equation}
The map $p'$ is again a split epimorphism and its kernel is
isomorphic to $P$.  Interchanging the maps acting on separate tensor
factors gives
\begin{align*}
 (p'\otimes_A1_M)(1_M\otimes_Bv)
 &=(p\otimes_A1_M)
   (t_M^{-1}\otimes_B1_N\otimes_A1_M)(1_M\otimes_Bv)\\
 &=(p\otimes_A1_M)(1_M\otimes_Bv)t_M^{-1}
   =1_M.
\end{align*}
This is the first triangle identity.

For the second triangle, put
\[
 s_N=(1_N\otimes_Ap')(v\otimes_B1_N):
 N\longrightarrow N\otimes_AM\otimes_BN\longrightarrow N.
\]
We first show that $s_N$ is idempotent.  Write
$v(1_B)=\sum_{j=1}^r n_j\otimes m_j$.  The first triangle identity gives
\begin{equation}
\label{eq:first-triangle-elements}
 \sum_{j=1}^r p'(m\otimes n_j)m_j=m
 \qquad(m\in M),
\end{equation}
while
$s_N(n)=\sum_i n_i p'(m_i\otimes n)$.  By the bimodule linearity
of $p'$ and \eqref{eq:first-triangle-elements},
\begin{align*}
 s_N^2(n)
 &=\sum_{i,j} n_i p'(m_i\otimes n_j)p'(m_j\otimes n)\\
 &=\sum_i n_i p'\left(
       \left(\sum_j p'(m_i\otimes n_j)m_j\right)\otimes n
       \right)\\
 &=\sum_i n_i p'(m_i\otimes n)=s_N(n).
\end{align*}

The map $v\otimes_B1_N$ is a split monomorphism with complement
$Q\otimes_BN\in\mathcal P_{B,A}$.  The map
$1_N\otimes_Ap'$ is a split epimorphism with kernel isomorphic to
$N\otimes_AP$, again in $\mathcal P_{B,A}$.  Hence $s_N$ is
invertible modulo maps through objects of $\mathcal P_{B,A}$.
Since $N$ is reduced, Lemma~\ref{lem:projective-factorization-radical}
implies that $s_N$ is an automorphism.  An invertible idempotent is
the identity, so $s_N=1_N$.

The bimodule maps $v$ and $p'$ therefore satisfy both triangle
identities.  Tensoring them with modules gives the unit
and counit of
\[
 F=M\otimes_B-\ \dashv\ G=N\otimes_A-.
\]
The unit at a $B$-module $X$ is
$v\otimes_B1_X:X\to GF(X)$, and the counit at an $A$-module $Y$
is $p'\otimes_A1_Y:FG(Y)\to Y$.

Apply the same construction with $(A,B,M,N)$ replaced by
$(B,A,N,M)$, starting with $i$ and $q$.  It gives $G\dashv F$.
The two constructions do not require any compatibility between the
chosen decompositions.

Finally, ${}_AM$ and ${}_BN$ are finitely generated projective.
The usual tensor--Hom adjunction and uniqueness of right adjoints
therefore give
\[
 N\cong\operatorname{Hom}_A(M,A),
 \qquad M\cong\operatorname{Hom}_B(N,B)
\]
as bimodules.  Projective biduality, applied respectively to
${}_BN$ and ${}_AM$, yields the other two displayed isomorphisms.
\end{proof}

\begin{corollary}
\label{cor:adjoint-type-without-separable-factors}
Let $A,B$ be finite central algebras over a commutative Artin ring
$R$, with no $R$-separable direct factors.  Suppose bimodules $M,N$
are projective on both sides and satisfy
\eqref{eq:adjoint-morita-identities}, with
$P\in\mathcal P_{A,A}$ and $Q\in\mathcal P_{B,B}$ for the chosen
projectivity class.  Then deleting the bimodule direct summands
in this class gives a pair satisfying the same type of
decompositions and inducing tensor functors which are adjoint in both
orders.  For ordinary projective bimodules, the resulting pair defines
a stable equivalence of adjoint type between $A$ and $B$.
\end{corollary}

\begin{proof}
Apply Lemma~\ref{lem:remove-projective-bimodule-summands} and
Theorem~\ref{thm:adjoint-normalization}.
\end{proof}

The same construction applies when the bimodule identities are given
only in the corresponding stable categories.

\begin{corollary}
\label{cor:adjoint-relative-stable-identities}
Fix either projectivity class.  Suppose that $M,N$ are
projective on both sides and that
\[
 M\otimes_BN\simeq A,
 \qquad N\otimes_AM\simeq B
\]
in the respective bimodule categories modulo maps through
$\mathcal P_{A,A}$ and $\mathcal P_{B,B}$, with separable direct factors
of $A,B$ allowed.  After removing from $M,N$
the summands in $\mathcal P_{A,B}$ and $\mathcal P_{B,A}$, the resulting
bimodules $M_0,N_0$ still satisfy these stable bimodule identities
and induce tensor functors which are adjoint in both orders on the
ordinary module categories.
They induce the same functors as $M,N$ on the $R$-relative stable
module categories, and these induced functors are inverse
equivalences.

The reduced pair has zero components on every $R$-separable
direct factor.  If either algebra has such a nonzero factor, the
corresponding full regular bimodule cannot be an actual direct
summand of the normalized composite.

If, in addition, $A,B$ have no $R$-separable direct factors, then
the stable identities can be represented by actual decompositions
\[
 M_0\otimes_BN_0\cong A\oplus P_0,
 \qquad N_0\otimes_AM_0\cong B\oplus Q_0,
\]
with $P_0\in\mathcal P_{A,A}$ and $Q_0\in\mathcal P_{B,B}$.
\end{corollary}

\begin{proof}
Deleting the indicated summands preserves the stable identities by
Lemma~\ref{lem:projective-cross-tensors}.  Replace $M,N$ by their
reduced parts.  Choose bimodule maps
\[
 p:M\otimes_BN\longrightarrow A,
 \qquad v:B\longrightarrow N\otimes_AM
\]
representing the given stable isomorphisms.  Tensoring preserves
factorizations through the chosen classes, by the same lemma.
Thus $(p\otimes_A1_M)(1_M\otimes_Bv)$ is invertible in the quotient
category of cross bimodules, and hence is an automorphism by
Lemma~\ref{lem:projective-factorization-radical}.  Normalize $p$ as
in \eqref{eq:normalized-counit}.  This gives the first triangle identity.
The calculation using
\eqref{eq:first-triangle-elements} proves that the second triangle
endomorphism is idempotent; it remains invertible in the quotient,
so reducedness implies that it is the identity.  Only the stable
invertibility of $p,v$ is used; their representatives need not be split.
Thus $M\otimes_B-\dashv N\otimes_A-$.
The symmetric construction proves the reverse adjunction.

For the assertion about separable factors, if $eA$ is an
$R$-separable direct factor, then
$eM\cong eA\otimes_AM\in\mathcal P_{A,B}$ by
Lemma~\ref{lem:projective-cross-tensors}.  Since it is a direct
summand of the reduced bimodule $M$, it is zero.  The same argument
applies to each component of $M,N$ supported on a separable factor.

Recall that the $R$-relative stable category of $A$-modules is the quotient
by maps through $\operatorname{add}\{A\otimes_RU\mid
U\in R\text{-}\mathrm{mod}\}$.  If a cross bimodule belongs to
$\mathcal P_{A,B}$, tensoring it with any $B$-module gives an
$R$-relative projective $A$-module.  Also, one-sided projectivity
ensures that $M\otimes_B-$ preserves relative projective modules.
It follows that the stable bimodule identities induce inverse
equivalences of the relative stable module categories.  Deleting the
indicated bimodule summands does not change these induced functors.

For the final assertion, decompose $M\otimes_BN=X_0\oplus P_0$
with $X_0$ reduced and $P_0\in\mathcal P_{A,A}$.
Under the additional assumption, $A$ is reduced by
Lemma~\ref{lem:unit-relative-projective-center}.  Thus $X_0$ and $A$ are reduced
objects isomorphic modulo the chosen class.  Representatives of
inverse isomorphisms have composites that are units by
Lemma~\ref{lem:projective-factorization-radical}, so the
representatives are actual isomorphisms.  Hence $X_0\cong A$.
The other decomposition follows in the same way.
\end{proof}

\begin{remark}
\label{rem:ordinary-relative-coincide-base-projective}
If a finite $C$-module $X$ is projective over $R$, then it is
$R$-relative projective over $C$ if and only if it is projective over
$C$ in the ordinary sense.  Indeed, an $R$-split epimorphism $R^n\to X$
induces an $R$-split $C$-epimorphism $C^n\to X$, which relative
projectivity splits over $C$.  Consequently, if $A,B$ are projective
over $R$, the two choices above coincide on all direct summands of
the one-sided projective bimodules $M,N$ and of their products.
They therefore coincide on the error modules in the displayed
decompositions as well.  The factorization ideals between
$R$-projective modules also agree: lift a map from a relative
projective module to the target across an $R$-split ordinary free
$C$-module epimorphism.  In this case the ordinary and relative
versions of the normalization statement agree.
\end{remark}

\begin{remark}
\label{rem:adjunction-versus-ordinary-stable-equivalence}
The adjunctions in Theorem~\ref{thm:adjoint-normalization} hold on
ordinary module categories, and the tensor functors are
exact by one-sided projectivity.  Over a nonfield base, they need not
induce inverse equivalences of the entire ordinary stable module
categories.  Indeed,
\[
 (A\otimes_RA)\otimes_AX\cong A\otimes_RX
\]
need not be projective as an $A$-module for an arbitrary $X$.
The proof uses the quotients of the module categories over
$C_{A,B}$ and $C_{B,A}$ by maps through the chosen projectivity classes
to establish the invertibility of the two triangle endomorphisms.
The stable bimodule identities give the relative stable equivalences
in Corollary~\ref{cor:adjoint-relative-stable-identities}.  An equivalence
of ordinary stable module categories is not required.
\end{remark}

\subsection{Henselian Noetherian local bases}
\label{subsec:adjoint-complete-base}

The normalization argument uses the radical inclusion in
Lemma~\ref{lem:projective-factorization-radical}, while the removal of
error summands also uses the Krull--Schmidt theorem.  Over a Henselian
Noetherian local base, both follow from semiperfectness of the relevant
endomorphism rings.  Completeness of the base is not needed.

\begin{proposition}
\label{prop:adjoint-complete-base-extension}
Let $(R,\mathfrak m)$ be a commutative Henselian Noetherian local
ring, and let $A,B$ be finite central $R$-algebras.  Fix either the
ordinary or the $R$-relative projectivity class from
Section~\ref{sec:adjoint-normalization}.
Then Theorem~\ref{thm:adjoint-normalization} and
Corollaries~\ref{cor:adjoint-type-without-separable-factors}
and~\ref{cor:adjoint-relative-stable-identities} remain valid over
$R$.

In particular, a reduced pair of finitely generated bimodules which are
projective on both sides and satisfy the Morita-type decompositions
induces tensor functors which are adjoint in both orders.  A pair
satisfying only the stable bimodule identities also has such a reduced
representative, with the same stable identities and induced
$R$-relative stable functors.  If $A,B$ have no $R$-separable direct
factors, the regular bimodules remain actual direct summands of the
composites.  No projectivity of $A,B$ over $R$ is required.
\end{proposition}

\begin{proof}
Let $C$ be a finite central $R$-algebra and $X$ a finitely generated
$C$-module.  Then $X$ is finite over $R$.  Since $R$ is Noetherian,
\[
 E=\operatorname{End}_C(X)
 \subseteq\operatorname{Hom}_R(X,X)
\]
is a finite $R$-module.  We first prove that $E$ is semiperfect.

Put $k=R/\mathfrak m$.  For $z\in\mathfrak mE$ and $y\in E$,
left and right multiplication by $1-yz$ induce the identity modulo
$\mathfrak m$.  Nakayama's lemma shows that both multiplication maps
are surjective, so $1-yz$ has a left and a right inverse.  Hence
\[
 \mathfrak mE\subseteq\operatorname{rad}E.
\]
Moreover, $E/\mathfrak mE$ is a finite-dimensional $k$-algebra.

We claim that idempotents lift modulo $\mathfrak mE$.  Let $a\in E$
lift an idempotent $\bar a\in E/\mathfrak mE$.  Choose generators
$e_1,\ldots,e_n$ of the $R$-module $E$ whose residue classes form a
$k$-basis of $E/\mathfrak mE$, and choose a matrix $T\in M_n(R)$
representing left multiplication by $a$ on these generators.  If
$f(Z)=\det(ZI-T)$, the determinant trick gives $f(a)=0$.
Indeed, the natural epimorphism $R^n\to E$ intertwines $T$ with
left multiplication by $a$, and Cayley--Hamilton then gives this
identity by evaluation at $1_E$.
The reduction of $T$ represents the idempotent operator given by
left multiplication by $\bar a$.  Consequently,
\[
 \bar f(Z)=Z^s(Z-1)^q
 \qquad(s+q=n).
\]
If $s=0$ or $q=0$, then $\bar a=1$ or $\bar a=0$, respectively,
and the claim is immediate.  Otherwise, the factorization criterion
for Henselian local rings gives monic polynomials $g,h\in R[Z]$ with
\[
 f=gh,\qquad \bar g=Z^s,\qquad \bar h=(Z-1)^q;
\]
see \cite[Lemma~10.153.3]{StacksHenselian2026}.
Their resultant is a unit, so $(g)+(h)=R[Z]$.  Choose $e(Z)\in R[Z]$
with
\[
 e\equiv0\pmod g,\qquad e\equiv1\pmod h.
\]
Then $e^2-e\in(f)$, and $e(a)$ is an idempotent of $E$.
The congruences $\bar e(0)=0$ and $\bar e(1)=1$ show that
$\bar e(\bar a)=\bar a$, so $e(a)$ lifts $\bar a$.
This argument applies to the possibly noncommutative ring $E$,
since the coefficients of the polynomials lie in the central
image of $R$.

The Artinian ring $E/\mathfrak mE$ has nilpotent radical, modulo
which idempotents lift.  Combining this fact with the preceding
lifting and the inclusion $\mathfrak mE\subseteq\operatorname{rad}E$
shows that $E$ is semiperfect.  Splitting primitive idempotents gives
a finite decomposition of $X$ into modules with local endomorphism
rings.  Thus the category of finitely generated $C$-modules is
Krull--Schmidt.

We now establish the radical inclusion needed for normalization.
Let $\mathcal S$ be an additive class of finite $C$-modules closed
under direct summands, and suppose that $X$ has no nonzero summand in
$\mathcal S$.  Write
\[
 X=\bigoplus_{i=1}^r X_i^{n_i},
 \qquad
 E/\operatorname{rad}E\cong\prod_{i=1}^r M_{n_i}(D_i),
\]
where the $X_i$ are pairwise nonisomorphic indecomposables and
\[
 D_i=\operatorname{End}_C(X_i)/
       \operatorname{rad}\operatorname{End}_C(X_i).
\]
For a map factoring through $\mathcal S$, each component between
isomorphic indecomposable summands also factors through
$\mathcal S$.  Under an identification with $X_i$, such a component
cannot be invertible: otherwise $X_i$ would be a direct summand of
an object of $\mathcal S$.  Since $\operatorname{End}_C(X_i)$ is
local, the component lies in its radical.  Components between
nonisomorphic indecomposables vanish in the displayed semisimple
quotient.  Therefore
\[
 \mathcal I_{\mathcal S}(X,X)
 \subseteq\operatorname{rad}\operatorname{End}_C(X).
\]
Units lift modulo an ideal contained in the Jacobson radical; no
nilpotence assumption on this ideal is needed.

The cross-tensor lemma and the equality of the ordinary and relative
projective centers use split maps and remain valid over $R$.
The triangle calculations in Theorem~\ref{thm:adjoint-normalization}
and Corollary~\ref{cor:adjoint-relative-stable-identities} now apply
using the radical inclusion just proved.  The Krull--Schmidt property
allows the indicated error summands to be removed.  When $A,B$ have
no $R$-separable direct factors, their regular bimodules are reduced
by Lemma~\ref{lem:unit-relative-projective-center}; the same
decomposition argument therefore retains the actual regular
bimodule summands.  This proves the assertions.
\end{proof}

\begin{corollary}
\label{cor:complete-base-individual-dimensions}
Under the hypotheses of Proposition~\ref{prop:adjoint-complete-base-extension},
suppose that the normalized bimodules $M,N$ satisfy the actual
decompositions \eqref{eq:adjoint-morita-identities}.  Then, for every
$A$-module $X$ and every $B$-module $Y$, including modules that are
not finitely generated, and for each
\[
 d\in\{\operatorname{pd},\operatorname{id},
       \operatorname{Gpd},\operatorname{Gid},\operatorname{domdim}\},
\]
one has
\[
 d_A(X)=d_B(N\otimes_AX),
 \qquad d_B(Y)=d_A(M\otimes_BY).
\]
The corresponding equalities hold for right modules, with the tensor
factors reversed.
\end{corollary}

\begin{proof}
The bimodule units and counits define both adjunctions on the
categories of all modules, and one-sided projectivity makes both
tensor functors exact.  The actual decompositions give natural
split inclusions
\[
 X\mid M\otimes_BN\otimes_AX,
 \qquad Y\mid N\otimes_AM\otimes_BY.
\]
Theorem~\ref{thm:adjoint-abstract-dimensions} applies.  The same
argument over the opposite algebras proves the right-module
assertions.  Stable bimodule identities alone do not supply the
actual direct-summand inclusions used here.
\end{proof}

The preceding results apply in the complete local setting of
Abe--Hoshino \cite{AbeHoshino2007} whenever the required one-sided
projectivity and stable bimodule identities hold.  Applying them to
a given derived equivalence requires a separate construction of
these bimodules.

The radical condition can also be stated separately from the
Henselian hypothesis.

\begin{proposition}
\label{prop:adjoint-radical-criterion}
Let $R$ be a commutative ring, let $A,B$ be unital central
$R$-algebras, and fix either projectivity class.  Suppose that
${}_AM_B$ and ${}_BN_A$ are finitely generated projective on each
side, and that
\[
 M\otimes_BN\simeq A,\qquad N\otimes_AM\simeq B
\]
in the corresponding stable bimodule categories.  Assume that
\begin{align*}
 \mathcal I_{\mathcal P_{A,B}}(M,M)
 &\subseteq\operatorname{rad}\operatorname{End}_{C_{A,B}}(M),\\
 \mathcal I_{\mathcal P_{B,A}}(N,N)
 &\subseteq\operatorname{rad}\operatorname{End}_{C_{B,A}}(N).
\end{align*}
Then $M\otimes_B-$ and $N\otimes_A-$ are adjoint in both orders on
the categories of all ordinary modules.  In particular, this conclusion
holds if $M,N$ are reduced and their displayed endomorphism rings
are semiperfect.
\end{proposition}

\begin{proof}
Choose representatives $p:M\otimes_BN\to A$ and
$v:B\to N\otimes_AM$ of the stable isomorphisms.  The cross-tensor
lemma shows that
\[
 t_M=(p\otimes_A1_M)(1_M\otimes_Bv)
\]
is invertible modulo $\mathcal I_{\mathcal P_{A,B}}(M,M)$.
The assumed radical inclusion makes $t_M$ an automorphism.
Replacing $p$ by $p(t_M^{-1}\otimes_B1_N)$ gives the first triangle
identity.  The calculation in
\eqref{eq:first-triangle-elements} shows that the second triangle
endomorphism is idempotent.  It is also invertible modulo the
second factorization ideal, and hence is an actual automorphism
by the other radical inclusion.  It is therefore the identity.
Interchanging $M,N$ gives the reverse adjunction.

For the last assertion, semiperfectness gives finite decompositions
of $M,N$ into indecomposables with local endomorphism rings.
The matrix argument in the proof of
Proposition~\ref{prop:adjoint-complete-base-extension} then gives
the two radical inclusions from reducedness.
\end{proof}

This criterion concerns the specified bimodule pair.  It does not
assert the existence of reduced representatives over an arbitrary
base, or replace stable bimodule identities by actual decompositions.

The next example shows why reducedness alone does not give the
radical condition over an arbitrary Noetherian local base.

\begin{example}
\label{ex:nonhenselian-stable-unit}
Let $k=\mathbb Q$ and put
\[
 R=k[x]_{(x)},\qquad
 A=R[t]/\bigl(t^2(t-1)-x\bigr).
\]
The algebra $A$ is free of rank three over $R$.  Eliminating $x$
identifies $A$ with a localization of $k[t]$, so $A$ is a domain.
Since $A$ is finite over the local ring $R$, its maximal ideals are
the two ideals lying over $(x)$, namely $(t)$ and $(t-1)$.  Thus
\[
 \operatorname{rad}A=(t)\cap(t-1)=(t(t-1)).
\]
In particular, $R$ is not Henselian: the finite commutative
$R$-algebra $A$ is a domain with two maximal ideals, whereas a
finite commutative algebra over a Henselian local ring is a product
of local rings \cite[Lemma~10.153.3]{StacksHenselian2026}.

Set $C=A^e_R$.  It is the corresponding localization of
\[
 k[t,u]/\bigl((t-u)g(t,u)\bigr),
 \qquad g(t,u)=t^2+tu+u^2-t-u.
\]
The regular bimodule is $A=C/(t-u)$, and
\[
 \operatorname{End}_C(A)=A,\qquad
 \operatorname{Hom}_C(A,C)=\operatorname{ann}_C(t-u)=(g).
\]
The annihilator equality follows by cancellation of $t-u$ in
$k[t,u]$ before localization.  Composing a map $A\to C$ with a map
$C\to A$ therefore gives precisely the ideal generated by
$g(t,t)$.  Since every finite projective $C$-module is a summand of
a finite free one, the full projective-factorization ideal is
\[
 \mathcal I_{C\text{-}\mathrm{proj}}(A,A)
 =(g(t,t))=(3t^2-2t)=(t).
\]
Here $3t-2$ is a unit in $A$, since its residues at the two maximal
ideals are $-2$ and $1$.

The domain $\operatorname{End}_C(A)=A$ has no nontrivial
idempotents, so $A$ is indecomposable as a bimodule.  It is not
$C$-projective, since its identity does not belong to $(t)$.
Consequently, $A$ is reduced, but
\[
 \mathcal I_{C\text{-}\mathrm{proj}}(A,A)=(t)
 \not\subseteq\operatorname{rad}\operatorname{End}_C(A)
 =(t(t-1)).
\]
By Lemma~\ref{lem:unit-relative-projective-center}, the same
factorization ideal and the same reducedness assertion hold for
the $R$-relative class.  The algebra $A$ has no nonzero
$R$-separable direct factor, because it is a domain and its regular
bimodule is not projective.

Multiplication by $1-t$ is a unit in the stable endomorphism ring
$A/(t)\cong k$, but it is not an automorphism of $A$, since
$1-t\in(t-1)$.  Thus an endomorphism can be stably invertible
without being an automorphism, even when the algebra is $R$-free
and the bimodule is free on both sides.
This shows that the radical argument cannot be extended to
arbitrary Noetherian local bases by reducedness alone.  It is not
a counterexample to the existence of a biadjoint normalization:
the identity pair $(A,A)$ is already biadjoint.
\end{example}

\subsection{Fields and the role of separable factors}
\label{sec:adjoint-field}

\subsubsection{The field case}

Separability of the semisimple quotients is distinct from the absence
of separable direct factors.  The former is not needed for adjoint
normalization.  Chen--Pan--Xi
\cite[Lemma~4.1, p.~650]{CPX2012} prove that, over an arbitrary field,
if $A,B$ have no semisimple direct factors and a
Morita-type pair $M,N$ has no projective bimodule summands, then
\[
 {}_BN_A\simeq\operatorname{Hom}_A(M,A)
       \simeq\operatorname{Hom}_{B^{\mathrm{op}}}(M,B),
 \qquad
 M\otimes_B-\dashv N\otimes_A-\dashv M\otimes_B-.
\]
They also prove the corresponding dual isomorphisms for $M$ and
projective-injectivity of the error bimodules over the enveloping
algebras.  Neither semisimple quotient is assumed separable.

Their proof removes the assumption on the separability of the
semisimple quotients from the arguments of Dugas--Mart\'inez-Villa
\cite[Theorem~2.7 and Corollary~3.1]{DMV2007}.
The criterion in \cite[Theorem~2.8]{DMV2007} uses separability to
deduce projectivity of a bimodule from projectivity of all its tensor
images.  It is used to construct Morita-type decompositions from
weaker functorial hypotheses.  Thus the inseparable-field example
following that criterion does not preclude normalization of a pair
which already satisfies the Morita-type decompositions.

\begin{corollary}
\label{thm:adjoint-arbitrary-field}
Let $A,B$ be finite-dimensional algebras over a field $k$, with
no separable direct factors.  Suppose that $M,N$ define a stable
equivalence of Morita type.  Then deleting their projective bimodule
direct summands gives a pair $M_0,N_0$ whose tensor functors are
adjoint in both orders and induce
the same functors on the ordinary stable categories.
\end{corollary}

\begin{proof}
Separability of a finite-dimensional algebra means projectivity
of its diagonal bimodule.  Since
$\operatorname{End}_{A^e}(A)=Z(A)$, diagonal direct summands
correspond to direct algebra factors.  Apply
Lemma~\ref{lem:remove-projective-bimodule-summands} and
Theorem~\ref{thm:adjoint-normalization} with $R=k$.
Over a field, a projective bimodule tensored with any module
is projective on the remaining side; hence the removed summands
induce zero functors on the ordinary stable categories.
\end{proof}

The existence assertion is included in
\cite[Lemma~3.11]{XZ2026}.  The proof above gives explicit units and
counits, without scalar extension or absolute indecomposability.
Over a perfect field, absence of separable
direct factors is equivalent to absence of semisimple direct
factors.  Over an imperfect field the former is weaker: a finite
purely inseparable field extension of degree greater than one
is semisimple as a ring but is not separable over the ground field.

\subsubsection{Reduction and separable direct factors}

An unreduced Morita-type pair need not induce adjoint tensor functors.

\begin{example}
\label{ex:adjoint-added-projective}
Let $A=k[t]/(t^2)$, put $B=A$, and set
\[
 M=A\oplus A^e,\qquad N=A.
\]
The bimodules are projective on both sides and satisfy
$M\otimes_AN\simeq A\oplus A^e\simeq N\otimes_AM$.
The tensor functor defined by $N$ is the identity, whose left
and right adjoints are both the identity up to natural isomorphism.
The functor defined by $M$ is not isomorphic to the identity, since it sends $A$ to
$A\oplus A^e\not\simeq A$.  Thus the specified tensor functors
are not adjoint to one another.  Deleting the projective bimodule
summand gives the adjoint pair $(A,A)$.
\end{example}

The condition on separable direct factors cannot be omitted, even
if the bimodules realizing the equivalence may be replaced.

\begin{example}
\label{ex:adjoint-kronecker-extra-block}
Let $A=kQ$ for the Kronecker quiver $1\rightrightarrows2$, and put
$B=A\times k$.  Define
\[
 {}_AM_B=({}_AA_A,{}_AA_k),\qquad
 {}_BN_A=({}_AA_A,{}_kA_A)^{\mathsf t}.
\]
The entries are the components for the two factors of $B$;
the $A$-actions are regular and the $k$-actions are scalar.
These bimodules are projective on each side, and
\begin{equation}
 M\otimes_BN\simeq A\oplus(A\otimes_kA).
 \label{eq:adjoint-extra-block-first}
\end{equation}
The error is free of rank one over $A^e$.
Choose a vector-space complement $V$ to $k1_A$ in $A$.
Recording the four bimodule components in matrix form, we obtain
\begin{equation}
 \begin{split}
 N\otimes_AM
 &\simeq
 \begin{pmatrix}{}_AA_A&{}_AA_k\\{}_kA_A&{}_kA_k\end{pmatrix}
 \simeq B\oplus Q,\\
 Q&=
 \begin{pmatrix}0&{}_AA_k\\{}_kA_A&{}_kV_k\end{pmatrix}.
 \end{split}
 \label{eq:adjoint-extra-block-second}
\end{equation}
Since $B^e\simeq A^e\times A\times A^{\mathrm{op}}\times k$,
the four components of $Q$ are projective over the corresponding
factors.  Thus $M,N$ define a classical stable equivalence of Morita type.

There is, however, no Morita-type pair between these algebras whose
tensor functors are adjoint in both orders.  With the path convention in which
left modules represent the displayed quiver, the indecomposable
projective $A$-modules have dimension vectors $(1,2),(0,1)$,
and the indecomposable injectives have vectors $(1,0),(2,1)$.
Thus $A$ has no nonzero projective-injective module, whereas the
$B$-module $S=(0,k)$ is nonzero and projective-injective.
If any Morita-type pair gave biadjoint exact functors
$F:B\text{-}\mathrm{mod}\to A\text{-}\mathrm{mod}$ and
$G:A\text{-}\mathrm{mod}\to B\text{-}\mathrm{mod}$, then $F$
would preserve projectives and injectives, forcing $F(S)=0$.
The Morita-type decomposition makes the identity functor a direct
summand of $GF$, so $S$ would be a summand of $GF(S)=0$.
The argument applies over every field, including algebraically closed
fields.
\end{example}

\subsubsection{Adjunction versus ordinary stable equivalence}

Xi--Zhang \cite[Lemma~3.11, p.~13]{XZ2026} state the existence of
adjoint-type representatives for Artin algebras without separable
direct summands.
Their Definition~3.7 is stated over a field, and their proof invokes
extensions of the earlier field arguments.
Theorem~\ref{thm:adjoint-normalization} proves the assertion directly
under the Morita-type bimodule conditions used here.  The adjunctions
hold on ordinary module categories, without projectivity over the
base or an equivalence of the ordinary stable module categories.

For example, take $R=k[\epsilon]/(\epsilon^2)$,
$A=R[t]/(t^2)$, $P=A\otimes_RA$, $M=A\oplus P$, and $N=A$.
The algebra is free over $R$, the bimodules are projective on
both sides, and both error bimodules are the projective
$A^e_R$-module $P$.  For $X=A/\epsilon A$, however, we have
\[
 P\otimes_AX\simeq X^{\oplus2},
 \qquad M\otimes_AX\simeq X^{\oplus3}.
\]
Here $X$ is indecomposable, since its endomorphism ring is
$k[t]/(t^2)$, and nonprojective, since $A$ is local and
$\epsilon X=0$.  Krull--Schmidt implies that $X^{\oplus3}$ is
not isomorphic to $X$ in the ordinary stable category.
Thus the given tensor functors are not inverse on the ordinary stable
category.  They are inverse on the $R$-relative stable category, and
their reduced pair is the identity pair.

\section{Self-injective dimensions and preservation of injectives}
\label{sec:bel-correction}

For Artin algebras, Tang proved that a stable equivalence of
Morita type preserves each of the left and right self-injective
dimensions \cite[Theorem~1]{bel-Tang2026}. We use retractions on
extension groups to prove the dimension equalities for unital
rings with the indicated bimodule splittings. The Iwanaga--Gorenstein
consequence concerns two-sided Noetherian rings. We also use
Tang's triangular matrix example \cite[Section~3]{bel-Tang2026}
to examine the assertion concerning injective modules in
Beligiannis's proof.

\subsection{Natural retractions on extension groups}

We use the usual order of composition and compute injective dimensions
in the categories of all modules, allowing the value $\infty$.

\begin{lemma}
\label{lem:bel-ext-retract}
Let $\mathcal A$ and $\mathcal B$ be abelian categories, and let
$F\colon\mathcal A\to\mathcal B$ and
$G\colon\mathcal B\to\mathcal A$ be exact additive functors.  Suppose
there are natural transformations
\[
 \eta\colon\operatorname{Id}_{\mathcal A}\longrightarrow GF,
 \qquad
 \pi\colon GF\longrightarrow\operatorname{Id}_{\mathcal A},
 \qquad \pi\eta=\operatorname{Id}.
\]
For every $X,Y\in\mathcal A$ and every $n\geq0$, applying $F$ to
extensions gives a split injection
\[
 F_*\colon \operatorname{Ext}^n_{\mathcal A}(X,Y)
 \longrightarrow
 \operatorname{Ext}^n_{\mathcal B}(FX,FY).
\]
The extension groups are understood in the Yoneda sense; for module
categories, they agree with the usual derived functors.
\end{lemma}

\begin{proof}
For an extension class
$\xi\in\operatorname{Ext}^n_{\mathcal B}(FX,FY)$, define
\[
 r^n_{X,Y}(\xi)
   :=(\pi_Y)_*(\eta_X)^*G_*(\xi).
\]
Here $G_*$ is followed by pullback along $\eta_X$ and pushout along
$\pi_Y$.  Since $F$ and $G$ are exact, the induced maps are well defined.  Naturality
of $\eta$, applied to a representative extension of $X$ by $Y$, gives
\[
 (\eta_X)^*(GF)_*(\alpha)=(\eta_Y)_*\alpha.
\]
Consequently
\[
 r^n_{X,Y}F_*(\alpha)
   =(\pi_Y)_*(\eta_Y)_*\alpha
   =(\pi_Y\eta_Y)_*\alpha
   =\alpha.
\]
For $n=0$, the same formula is
$r^0_{X,Y}(h)=\pi_YG(h)\eta_X$.  The construction is natural in both
variables.
\end{proof}

\begin{remark}
\label{rem:bel-cochain}
For module categories, if $F$ preserves projectives, the preceding
retraction can be described on cochain complexes.  For a projective
resolution $L_\bullet\to X$, its degree-$i$ component is
\[
 \operatorname{Hom}_{B}(F L_i,FY)
 \longrightarrow\operatorname{Hom}_{A}(L_i,Y),
 \qquad h\longmapsto\pi_YG(h)\eta_{L_i}.
\]
Naturality of $\eta$ shows that this is a cochain map, and
$\pi\eta=1$ shows that it retracts the map induced by $F$.  This is
the comparison used in \cite[proof of Theorem~1]{bel-Tang2026}.
\end{remark}

\begin{corollary}
\label{cor:bel-id-reflect}
Let $F:A\text{-}\mathrm{Mod}\to B\text{-}\mathrm{Mod}$ and
$G:B\text{-}\mathrm{Mod}\to A\text{-}\mathrm{Mod}$ satisfy the hypotheses
of Lemma~\ref{lem:bel-ext-retract}.  Then
\[
 \operatorname{id}_A X\leq \operatorname{id}_B(FX)
\]
for every left $A$-module $X$.  In particular, $F$ reflects injectivity.
\end{corollary}

\begin{proof}
If $\operatorname{id}_B(FX)\leq d$, the lemma gives
$\operatorname{Ext}^{d+1}_A(U,X)=0$ for every left $A$-module $U$.
The usual criterion for injective dimension gives the assertion.
\end{proof}

\subsection{The regular modules}

\begin{theorem}
\label{thm:bel-selfinj}
Let $A$ and $B$ be unital rings.  Suppose that
${}_AM_B$ and ${}_BN_A$ are finitely generated projective on each side
and that there are bimodule decompositions
\[
 M\otimes_BN\simeq A\oplus P,
 \qquad N\otimes_AM\simeq B\oplus Q.
\]
Then, without any projectivity assumption on $P,Q$, one has
\[
 \operatorname{id}_A A=\operatorname{id}_B B,
 \qquad
 \operatorname{id}_{A^{\mathrm{op}}} A
   =\operatorname{id}_{B^{\mathrm{op}}} B.
\]
Both equalities hold in $\mathbb N_0\cup\{\infty\}$.  No Noetherian
hypothesis, adjunction, or underlying commutative base is required.
\end{theorem}

\begin{proof}
Put $F=N\otimes_A-$ and $G=M\otimes_B-$.  These functors are exact
because $N_A$ and $M_B$ are projective.  The bimodule
decompositions yield natural splittings
\[
 GF(X)\simeq X\oplus(P\otimes_AX),
 \qquad FG(Y)\simeq Y\oplus(Q\otimes_BY).
\]
Suppose $d=\operatorname{id}_B B<\infty$.  The left $B$-module
$F(A)\simeq{}_BN$ is a summand of $B^t$ for a finite $t$, and hence
$\operatorname{id}_B F(A)\leq d$.  By
Corollary~\ref{cor:bel-id-reflect},
\[
 \operatorname{id}_A A
 \leq\operatorname{id}_B F(A)
 \leq\operatorname{id}_B B.
\]
Interchanging $A,B$ and $F,G$ gives the reverse inequality whenever
either dimension is finite.  Hence the left self-injective dimensions
are equal, including the case in which both are infinite.

For right modules, use $-\otimes_AM$ and $-\otimes_BN$.  Their exactness
follows from left projectivity, and their values on the regular modules are
the finitely generated projective right modules $M_B$ and $N_A$.
The two composites again contain the identity functors as natural direct
summands.  The same argument proves the right-hand equality.
\end{proof}

For stable equivalences of Morita type between Artin algebras, this
recovers Tang's theorem \cite[Theorem~1]{bel-Tang2026}.

\begin{corollary}
\label{cor:bel-artin-selfinj}
\label{cor:bel-noetherian-selfinj}
Let $A,B$ be two-sided Noetherian rings.  Suppose that bimodules
$ {}_AM_B,{}_BN_A$, finitely generated projective on both sides, satisfy
\[
 M\otimes_BN\simeq A\oplus P,
 \qquad N\otimes_AM\simeq B\oplus Q
\]
for arbitrary complement bimodules $P,Q$.  Then the left and right
self-injective dimensions are preserved separately.  In particular,
$A$ is Iwanaga--Gorenstein if and only if $B$ is.  This applies to
finite algebras over a commutative Noetherian ring, and hence to
Artin algebras, whenever the displayed bimodule decompositions hold.
\end{corollary}

\begin{proof}
Apply Theorem~\ref{thm:bel-selfinj}.  For a two-sided Noetherian ring,
simultaneous finiteness of its left and right self-injective dimensions
is the Iwanaga--Gorenstein property.  An algebra finite over a
commutative Noetherian ring is Noetherian on both sides.
\end{proof}

The proof uses neither projectivity over a commutative base nor an equivalence of
ordinary stable module categories.  It also requires no adjunction or
preservation of injective modules.  Thus, over a nonfield base, the
equalities follow from the bimodule splittings themselves.  The earlier result of Liu
and Xi assumes adjunction in both orders
\cite[Proposition~4.1, p.~579]{bel-LiuXi2007}; their Question~1 on
p.~584 asks for one-sided finiteness without that extra assumption.

In the Artin and Henselian settings of
Theorem~\ref{thm:adjoint-normalization} and
Proposition~\ref{prop:adjoint-complete-base-extension}, respectively,
and in the absence of $R$-separable direct factors, deleting the
projective bimodule summands gives a pair whose tensor functors are
adjoint in both orders.
Preservation of injectives then holds for this reduced pair.  The
retraction argument applies to the original functors, allows arbitrary
complement bimodules, and imposes no condition on direct algebra factors.

\begin{proposition}
\label{prop:bel-projective-dimensions}
Under the hypotheses of Theorem~\ref{thm:bel-selfinj}, the original
tensor functors preserve the projective dimensions of all left and
right modules.  More precisely,
\[
\begin{aligned}
 \operatorname{pd}_B(N\otimes_AX)&=\operatorname{pd}_A X,\\
 \operatorname{pd}_A(M\otimes_BY)&=\operatorname{pd}_B Y,\\
 \operatorname{pd}_{B^{\mathrm{op}}}(U\otimes_AM)
   &=\operatorname{pd}_{A^{\mathrm{op}}}U,\\
 \operatorname{pd}_{A^{\mathrm{op}}}(V\otimes_BN)
   &=\operatorname{pd}_{B^{\mathrm{op}}}V .
\end{aligned}
\]
All equalities include infinite values.  Neither adjunction nor
enveloping-algebra projectivity of the complement bimodules is needed.
\end{proposition}

\begin{proof}
The exact functors $F=N\otimes_A-$ and $G=M\otimes_B-$ preserve
projectives.  Since $X$ is a direct summand of $GF(X)$, one has
\[
 \operatorname{pd}_A X
 \leq \operatorname{pd}_A GF(X)
 \leq \operatorname{pd}_B F(X)
 \leq \operatorname{pd}_A X .
\]
This proves the first equality, and the second follows by symmetry.
The same argument applies to the right-module tensor functors.
\end{proof}

\subsection{A tensor functor that does not preserve injectives}

The following example is due to Tang \cite[Section~3]{bel-Tang2026}.

\begin{example}
\label{ex:bel-triangular}
Let $K$ be a field and let $A=T_2(K)$ be the algebra of upper triangular
$2\times2$ matrices.  Write $e_{ij}$ for its matrix units.  Set
\[
 E=A\otimes_K A,
 \qquad M=A,
 \qquad N=A\oplus E,
\]
where $a\cdot(u\otimes v)\cdot b=au\otimes vb$.  The map
\[
 A\otimes_K A^{\mathrm{op}}\longrightarrow E,
 \qquad u\otimes v^{\mathrm{op}}\longmapsto u\otimes v
\]
is an isomorphism of left modules over $A\otimes_KA^{\mathrm{op}}$.
Thus $E$ is a projective $A$-bimodule.  It is also free of rank $3$ on
each side.  Both $M$ and $N$ are therefore finitely generated projective
on each side, and
\[
 M\otimes_A N\simeq A\oplus E,
 \qquad N\otimes_A M\simeq A\oplus E.
\]
These are the defining identities of a stable self-equivalence of
Morita type.

For every left $A$-module $X$ one has
\[
 N\otimes_A X\simeq X\oplus(A\otimes_KX).
\]
The second summand is a free left $A$-module.  Hence this tensor functor
induces the identity on the ordinary stable category of all modules.

The left ideal $Ae_{11}=Ke_{11}$ is not injective.  Indeed, right
multiplication by $e_{12}$ gives a monomorphism
\[
 j\colon Ae_{11}\longrightarrow Ae_{22},
 \qquad e_{11}\longmapsto e_{12}.
\]
For an $A$-linear $r\colon Ae_{22}\to Ae_{11}$,
\[
 r(e_{12})=r(e_{12}e_{22})=e_{12}r(e_{22})=0,
\]
because $e_{12}Ae_{11}=0$.  Thus $rj$ cannot be the identity.

Let $I=D_K(e_{22}A)=\operatorname{Hom}_K(e_{22}A,K)$.  It is a
one-dimensional injective left $A$-module: $e_{22}A$ is a projective
right $A$-module and
\[
 \operatorname{Hom}_A(-,I)
 \simeq D_K(e_{22}A\otimes_A-)
\]
is exact.  Nevertheless,
\[
 N\otimes_A I\simeq I\oplus(A\otimes_KI)\simeq I\oplus A
\]
is not injective, since it contains the noninjective summand $Ae_{11}$.

The injective dimensions can also be computed explicitly.  The pairing
\[
 e_{11}A\times Ae_{22}\longrightarrow K,
 \qquad (y,x)\longmapsto
       \text{the coefficient of $e_{12}$ in $yx$}
\]
identifies $Ae_{22}$ with $D_K(e_{11}A)$, so $Ae_{22}$ is injective.
Moreover,
\[
 0\longrightarrow Ae_{11}\xrightarrow{j}Ae_{22}
   \longrightarrow I\longrightarrow0
\]
is exact.  Hence
\[
 \operatorname{id}_A I=0,
 \qquad \operatorname{id}_A A=1,
 \qquad \operatorname{id}_A(N\otimes_A I)=1.
\]
Thus the functor does not preserve the injective dimension of every
module, although the self-injective dimensions of the algebras agree.
The same example gives
\[
 \operatorname{Gid}_A I=0,\qquad
 \operatorname{Gid}_A(N\otimes_AI)=1,
\]
because Gorenstein injective dimension agrees with injective dimension
when the latter is finite; see \cite{Holm2004}.

It also distinguishes the dominant dimensions of individual modules.
The monomorphism $j$ is essential, and the exact sequence above gives
the minimal injective resolution
\[
 0\longrightarrow A\longrightarrow (Ae_{22})^2
   \longrightarrow I\longrightarrow0 .
\]
Here $Ae_{22}$ is projective-injective and $I$ is not projective.
Thus $\operatorname{domdim}_A A=1$, whereas
$\operatorname{domdim}_A(Ae_{22})=\infty$.  Since
\[
 N\otimes_A Ae_{22}\simeq Ae_{22}\oplus A^2,
\]
the latter module has dominant dimension $1$.  Thus individual
injective, Gorenstein injective, and dominant dimensions need not be
preserved by the original pair, even over a field.
Passing to opposite algebras gives the analogous right-module example.
\end{example}

\subsection{The gap in Beligiannis's argument}

Beligiannis \cite[Theorem~8.12, p.~192]{bel-Beligiannis2005} states
preservation of the stable Cohen--Macaulay category, virtual
Gorensteinness, and Gorensteinness, in parts (i), (ii), and (iii),
respectively.  His definition immediately preceding that theorem
includes an equivalence of ordinary stable categories induced by the
specified tensor functor.  Over an arbitrary Artinian base, this is
stronger than the bimodule splitting conditions alone.

In the proof of part (i), on the same page, the left tensor functor
\[
 M\otimes_\Gamma-\colon
 \operatorname{Mod}\text{-}\Gamma^{\mathrm{op}}
 \longrightarrow
 \operatorname{Mod}\text{-}\Lambda^{\mathrm{op}}
\]
is asserted to preserve injectives.  Example~\ref{ex:bel-triangular},
with the names of the two defining bimodules interchanged, disproves
this assertion under the stated hypotheses.  The example satisfies
even the ordinary stable-equivalence requirement, so this failure is
independent of any ambiguity about the nonfield base.

Thus the injectivity assertion used in the proof fails, whereas
Theorem~8.12(iii) remains valid.  Corollary~\ref{cor:bel-artin-selfinj}
proves this part without injective preservation and gives the stronger
separate equalities of self-injective dimensions.  The other assertions
in parts (i) and (ii) require separate arguments.

For comparison, Example~\ref{ex:bel-triangular} also disproves the
preservation direction of Sun--Zhao
\cite[Lemma~3.1(1), pp.~2592--2593]{bel-SunZhao2025} and the
individual-module equality in their Proposition~3.2(1), p.~2593.
The opposite example addresses the right-module assertions.
The equalities for the regular modules and the resulting Gorenstein
conclusion remain valid.

\begin{remark}
\label{rem:bel-projections}
The extension-of-maps argument in
\cite[Lemma~3.1(1), p.~2593]{bel-SunZhao2025} uses an additional
compatibility between two projections; see also
\cite[Section~3]{bel-Tang2026}.
Let $\pi\colon GF\to\operatorname{Id}$ and
$\rho\colon FG\to\operatorname{Id}$ be the two projections furnished
by the bimodule splittings.  They do not imply
\[
 F(\pi_I)=\rho_{F(I)}.
\]
In Example~\ref{ex:bel-triangular}, $G=\operatorname{Id}$ and
$F(Z)=Z\oplus(A\otimes_KZ)$.  On the summand $A\otimes_KF(I)$ of
$F(F(I))$, the map $\rho_{F(I)}$ is zero, whereas $F(\pi_I)$ is
$1_A\otimes\pi_I$ and has nonzero image $A\otimes_KI$.
The retraction in Lemma~\ref{lem:bel-ext-retract} does not require
these projections to agree.
\end{remark}

\section{Individual Gorenstein dimensions}
\label{sec:individual-gorenstein-dimensions}

Let $A$ be a unital ring. We denote by $\mathrm{GProj}(A)$
the category of Gorenstein projective left $A$-modules and by
$\mathrm{Gproj}(A)$ its full subcategory of finitely generated
modules. We work with arbitrary modules unless a category of
finitely generated modules is specified. We take dimensions in
$\mathbb N_0\cup\{\infty\}$, putting $\operatorname{Gpd}_A(0)=0$.
Then adjoining a projective summand leaves the dimension unchanged.
With this convention, $\operatorname{Gpd}_A$ agrees with the
Auslander--Bridger dimension $\mathrm{G\!\! -\! dim}_A$ on finitely
generated modules over a two-sided Noetherian ring.

Recall that $\mathrm{GProj}(A)$ is resolving and closed under direct
summands, that Gorenstein projective dimension does not increase on
direct summands, and that, for every $n\geq 0$,
\begin{equation}
  \operatorname{Gpd}_A X\leq n
  \quad\Longleftrightarrow\quad
  \Omega_A^nX\in\mathrm{GProj}(A).
  \label{eq:gpd-syzygy-test}
\end{equation}
These facts hold without assuming that $A$ is Gorenstein; see
\cite[Theorem~2.5, Propositions~2.16, 2.18--2.19, and
Theorem~2.20]{Holm2004}.  Over a two-sided Noetherian ring, finitely
generated Gorenstein projectives admit totally acyclic complexes of
finitely generated projectives, and the same criterion holds in the
finite-module category.  Isomorphic objects in either ordinary stable
module category have the same Gorenstein projective dimension.

\subsection{Dimensions from exact biadjoints}
\label{subsec:adjoint-individual-dimensions}

We first isolate the homological argument from the hypotheses used
to construct an adjoint pair.  All dimensions in this subsection are
computed in the categories of all modules.  In particular, injective
resolutions are not required to have finitely generated terms.
We allow infinite values, put
$\operatorname{pd}(0)=\operatorname{id}(0)
=\operatorname{Gpd}(0)=\operatorname{Gid}(0)=0$, and use
$\operatorname{domdim}(0)=\infty$.

The Gorenstein-dimension part of the next theorem also follows from
\cite[Theorem~3.2 and Remark~3.3]{ChenRen2022}: the identity
retractions make each Frobenius functor faithful. We include the
short resolution argument to treat the five dimensions together.
The adjoint-type Artin antecedents are
\cite[Lemma~4.2 and Proposition~4.5]{bel-LiuXi2007}.

\begin{theorem}
\label{thm:adjoint-abstract-dimensions}
Let $A,B$ be unital rings, and let
\[
 F:A\text{-}\mathrm{Mod}\longrightarrow B\text{-}\mathrm{Mod},
 \qquad
 G:B\text{-}\mathrm{Mod}\longrightarrow A\text{-}\mathrm{Mod}
\]
be exact additive functors which are adjoint in both orders.  Suppose
that the identity functors are natural direct summands of their
composites:
\begin{equation}
 \mathrm{Id}_{A\text{-}\mathrm{Mod}}\mid GF,\qquad
 \mathrm{Id}_{B\text{-}\mathrm{Mod}}\mid FG.
 \label{eq:adjoint-identity-retracts}
\end{equation}
Then, for every left $A$-module $X$, every left $B$-module $Y$, and
every
\[
 d\in\{\operatorname{pd},\operatorname{id},
          \operatorname{Gpd},\operatorname{Gid},
          \operatorname{domdim}\},
\]
one has
\[
 d_B(FX)=d_A(X),\qquad d_A(GY)=d_B(Y).
\]
No Noetherian or finite generation hypothesis is required.
\end{theorem}

\begin{proof}
An exact functor with an exact right adjoint preserves projectives,
and one with an exact left adjoint preserves injectives.  Thus both
$F$ and $G$ preserve projective and injective modules.

Let $T^\bullet$ be a totally acyclic complex of projective $A$-modules.
Then $F(T^\bullet)$ is an acyclic complex of projective $B$-modules.
For every projective $B$-module $Q$, the adjunction gives
\[
 \operatorname{Hom}_B(F(T^\bullet),Q)
 \simeq\operatorname{Hom}_A(T^\bullet,GQ),
\]
which is acyclic since $GQ$ is projective.  Hence $F$ preserves
Gorenstein projective modules.  If $I^\bullet$ is a totally acyclic
complex of injective $A$-modules and $J$ is injective over $B$, the
other adjunction gives
\[
 \operatorname{Hom}_B(J,F(I^\bullet))
 \simeq\operatorname{Hom}_A(GJ,I^\bullet).
\]
The right-hand complex is acyclic since $GJ$ is injective.  Thus $F$
also preserves Gorenstein injective modules.  The same arguments
apply to $G$.

Applying these functors to finite resolutions or coresolutions gives
$d_B(FX)\leq d_A(X)$ and $d_A(GY)\leq d_B(Y)$ for the first four
dimensions.  These dimensions do not increase on direct summands;
for the Gorenstein dimensions over arbitrary rings, see
\cite{Holm2004}.  Hence \eqref{eq:adjoint-identity-retracts} gives
\[
 d_A(X)\leq d_A(GFX)\leq d_B(FX)\leq d_A(X),
\]
and similarly for $Y$.  This proves the first four equalities,
including infinite values.

For dominant dimension, take a minimal injective resolution of $X$.
If its first $n$ terms are projective, applying $F$ gives an injective
resolution of $FX$ whose first $n$ terms are projective.  Over any
ring, a minimal injective resolution is a termwise direct summand of
an arbitrary injective resolution: this follows by successively
splitting the injective envelopes of the successive cokernels.
Consequently,
\[
 \operatorname{domdim}_B(FX)\geq\operatorname{domdim}_A(X).
\]
The same argument applies to $G$.  Finite direct sums of injective
envelopes are injective envelopes, so dominant dimension of a finite
direct sum is the minimum of the dominant dimensions of its summands.
Using $X\mid GFX$, we obtain
\[
 \operatorname{domdim}_A(X)
 \geq\operatorname{domdim}_A(GFX)
 \geq\operatorname{domdim}_B(FX)
 \geq\operatorname{domdim}_A(X).
\]
This proves the remaining equality.
\end{proof}

\begin{remark}
\label{rem:adjoint-noetherian-dimensions}
If $A,B$ are two-sided Noetherian and the functors preserve finite
modules, the Gorenstein projective dimension equalities also give
equalities of Auslander--Bridger $G$-dimensions on finite modules.
Injective, Gorenstein injective, and dominant dimensions are still
computed in $A\text{-}\mathrm{Mod}$ and $B\text{-}\mathrm{Mod}$:
over a Noetherian ring, an injective envelope of a finite module need
not be finite.
\end{remark}

\begin{corollary}
\label{cor:adjoint-individual-dimensions}
Let $A,B$ be finite algebras over a commutative Artin ring, or over a
commutative Henselian Noetherian local ring.  Let $M_0,N_0$ be a
reduced pair to which Theorem~\ref{thm:adjoint-normalization} or
Proposition~\ref{prop:adjoint-complete-base-extension} applies.
Assume that the regular bimodules are actual direct summands:
\[
 M_0\otimes_BN_0\simeq A\oplus P_0,\qquad
 N_0\otimes_AM_0\simeq B\oplus Q_0.
\]
Put $F=N_0\otimes_A-$ and $G=M_0\otimes_B-$.
For every left $A$-module $X$, every left $B$-module $Y$, and every
\[
 d\in\{\operatorname{pd},\operatorname{id},
          \operatorname{Gpd},\operatorname{Gid},
          \operatorname{domdim}\},
\]
one has
\begin{equation}
 d_B(FX)=d_A(X),\qquad d_A(GY)=d_B(Y).
 \label{eq:adjoint-individual-dimensions}
\end{equation}
For every right $A$-module $U$ and right $B$-module $V$, one has
\begin{equation}
 \begin{aligned}
 d_{B^{\mathrm{op}}}(U\otimes_A M_0)&=d_{A^{\mathrm{op}}}(U),\\
 d_{A^{\mathrm{op}}}(V\otimes_B N_0)&=d_{B^{\mathrm{op}}}(V).
 \end{aligned}
 \label{eq:adjoint-right-dimensions}
\end{equation}
In particular, the reduced functors preserve and reflect Gorenstein
projective and Gorenstein injective modules.  On finite modules they
preserve the left and right Auslander--Bridger $G$-dimensions.
No projectivity or injectivity over the base is required, and the
functors need not give inverse equivalences of ordinary stable
module categories.
\end{corollary}

\begin{proof}
The normalization results give bimodule dual isomorphisms, hence
adjunctions in both orders on the categories of all modules.
One-sided projectivity makes the tensor functors exact.  The displayed
bimodule decompositions give
\eqref{eq:adjoint-identity-retracts}, so
Theorem~\ref{thm:adjoint-abstract-dimensions} proves
\eqref{eq:adjoint-individual-dimensions}.  The two bimodule dual
isomorphisms also give adjunctions in both orders between the
right-module tensor functors.  Applying the theorem to the opposite
rings proves \eqref{eq:adjoint-right-dimensions}.
The rings are two-sided Noetherian and the tensor functors preserve
finite modules, which gives the assertion about $G$-dimension.
\end{proof}

\begin{corollary}
\label{cor:adjoint-algebra-dimensions}
Under the hypotheses of Corollary~\ref{cor:adjoint-individual-dimensions},
one has
\[
 \begin{aligned}
 \operatorname{domdim}({}_AA)&=\operatorname{domdim}({}_BB),&
 \operatorname{domdim}(A_A)&=\operatorname{domdim}(B_B),\\
 \operatorname{id}_A A&=\operatorname{id}_B B,&
 \operatorname{id}_{A^{\mathrm{op}}}A
   &=\operatorname{id}_{B^{\mathrm{op}}}B.
 \end{aligned}
\]
\end{corollary}

\begin{proof}
The left $B$-module $N_0$ is a finite projective generator.  Indeed,
$B$ is a summand of $N_0\otimes_A M_0$, which belongs to
$\operatorname{add}({}_B N_0)$ because $ {}_A M_0$ is finite projective.
Thus $\operatorname{add}({}_B N_0)=\operatorname{add}({}_B B)$.
Injective dimension and dominant dimension each have the same value
on two modules with the same additive closure.  Apply
\eqref{eq:adjoint-individual-dimensions} to $X=A$, so that $F(A)=N_0$.
The same argument for right projective generators proves the
equalities of right self-injective and right dominant dimensions.
\end{proof}

\begin{remark}
\label{rem:adjoint-normalized-versus-original}
The five individual equalities above concern the reduced pair with
actual identity summands.  Adding projective bimodule summands can
destroy preservation of individual injective, Gorenstein injective,
and dominant dimensions, even over a field; see
Example~\ref{ex:bel-triangular}.  Projective dimension is already
preserved by the original pair under the bimodule splittings alone,
by Proposition~\ref{prop:bel-projective-dimensions}.
For the original pair, the criteria below give preservation of
individual Gorenstein projective dimensions.  In particular,
Theorem~\ref{thm:gpd-base-injective-semt} treats base injectivity, and
Corollary~\ref{cor:gpd-gorenstein-noetherian-base} treats finite
projective algebras over a commutative Noetherian base of finite
self-injective dimension.  The separate equalities of self-injective
dimensions for the rings also follow from
Theorem~\ref{thm:bel-selfinj}, without normalization or adjunction.
\end{remark}

\subsection{A criterion in the ordinary stable module category}

\begin{proposition}
\label{prop:gpd-stable-criterion}
Let $A$ and $B$ be two-sided Noetherian rings.  Suppose that
\[
  \Phi:A\text{-}\underline{\mathrm{mod}}
       \longrightarrow B\text{-}\underline{\mathrm{mod}}
\]
is an equivalence satisfying the following two conditions:
\begin{enumerate}
\item there is a natural isomorphism
      $\Phi\Omega_A\simeq\Omega_B\Phi$;
\item the restriction of $\Phi$ is an equivalence
      $\underline{\mathrm{Gproj}}(A)
       \longrightarrow\underline{\mathrm{Gproj}}(B)$.
\end{enumerate}
Then, for every $X\in A\text{-}\mathrm{mod}$, one has
\[
  \operatorname{Gpd}_B\Phi(X)=\operatorname{Gpd}_A X.
\]
Here $\Phi(X)$ denotes any module representing its stable isomorphism
class.  The equality includes the value $\infty$.
\end{proposition}

\begin{proof}
We first show that $X$ is Gorenstein projective if and only if $\Phi(X)$
is Gorenstein projective.  The forward implication follows from (2).
For the converse, suppose that $\Phi(X)$ is Gorenstein projective.
Essential surjectivity of the restricted functor gives
$U\in\mathrm{Gproj}(A)$ and an isomorphism
$\Phi(X)\simeq\Phi(U)$ in the stable category.  Full faithfulness of
$\Phi$ on the entire stable module category gives $X\simeq U$ there.
Thus $X$ is Gorenstein projective, since this property is unchanged by
adding or removing projective direct summands.

For every $n\geq 0$, condition~(1) and
\eqref{eq:gpd-syzygy-test} now give
\[
\begin{aligned}
 \operatorname{Gpd}_A X\leq n
 &\Longleftrightarrow \Omega_A^nX\in\mathrm{Gproj}(A)\\
 &\Longleftrightarrow \Phi(\Omega_A^nX)\in\mathrm{Gproj}(B)\\
 &\Longleftrightarrow \Omega_B^n\Phi(X)\in\mathrm{Gproj}(B)\\
 &\Longleftrightarrow \operatorname{Gpd}_B\Phi(X)\leq n.
\end{aligned}
\]
The equality follows for both finite and infinite dimensions.
\end{proof}

\begin{remark}
\label{rem:gpd-requires-reflection}
The proof uses full faithfulness on the entire ordinary stable module
category to reflect Gorenstein projectivity of arbitrary syzygies.
An equivalence between the stable categories of Gorenstein projective
modules alone does not provide this reflection property.
\end{remark}

\subsection{Exact functors with exact stable inverses}

\begin{lemma}
\label{lem:stable-full-faithful-reflects-split-epi}
Let $A,B$ be unital rings and
$G:B\text{-}\mathrm{Mod}\to A\text{-}\mathrm{Mod}$ an additive
functor preserving projectives.  If its induced functor on ordinary
stable module categories is fully faithful, then $G$ reflects split
epimorphisms.
\end{lemma}

\begin{proof}
Let $p:E\twoheadrightarrow Z$ be an epimorphism, and let
$s:GZ\to GE$ satisfy $G(p)s=1_{GZ}$.  Fullness supplies a map $t:Z\to E$
whose stable class is sent to the stable class of $s$.  Faithfulness
then implies that $pt-1_Z$ factors through a projective module $Q$.
Write $pt-1_Z=ba$ with $a:Z\to Q$ and $b:Q\to Z$.  Lift $b$ through
$p$ to a map $\widetilde b:Q\to E$.  Then
$t-\widetilde b a$ is a section of $p$.
\end{proof}

\begin{theorem}
\label{thm:gpd-exact-stable-pair}
Let $A$ and $B$ be unital rings, and let
\[
 F:A\text{-}\mathrm{Mod}\to B\text{-}\mathrm{Mod},\qquad
 G:B\text{-}\mathrm{Mod}\to A\text{-}\mathrm{Mod}
\]
be exact functors preserving projectives.  Suppose that their induced
functors on the ordinary stable module categories are mutually inverse
equivalences.  Then $F$ and $G$ preserve Gorenstein projective modules,
and
\[
 \operatorname{Gpd}_B F(X)=\operatorname{Gpd}_A X,
 \qquad
 \operatorname{Gpd}_A G(Y)=\operatorname{Gpd}_B Y
\]
for all $X\in A\text{-}\mathrm{Mod}$ and $Y\in B\text{-}\mathrm{Mod}$.
No adjunction or preservation of injective modules is required.
\end{theorem}

\begin{proof}
Let $T^\bullet$ be a totally acyclic complex of
projective $A$-modules.  Exactness and preservation of projectives imply
that $F(T^\bullet)$ is an acyclic complex of projective $B$-modules.
For each cycle $Z$ of $T^\bullet$, the module $GF(Z)$ is stably
isomorphic to $Z$ and is therefore Gorenstein projective.

Fix a projective $B$-module $Q$, and consider an extension
\[
 0\longrightarrow Q\longrightarrow E\xrightarrow{p}F(Z)
 \longrightarrow0.
\]
After applying $G$, this extension splits because $GQ$ is projective
and $GF(Z)$ is Gorenstein projective.  Lemma
\ref{lem:stable-full-faithful-reflects-split-epi} shows that $p$ splits.
It follows that
$\operatorname{Ext}_B^1(F(Z),Q)=0$ for every cycle $Z$ and every
projective $Q$.  Applying this vanishing to the short exact sequences
of cycles of $F(T^\bullet)$ shows that
$\operatorname{Hom}_B(F(T^\bullet),Q)$ is acyclic.  Hence $F$ preserves
Gorenstein projective modules.  Interchanging $F$ and $G$ gives the
same conclusion for $G$.

An exact functor preserving Gorenstein projectives does not increase
Gorenstein projective dimension.  Since $GF(X)$ is stably isomorphic
to $X$, we obtain
\[
 \operatorname{Gpd}_A X
 =\operatorname{Gpd}_A GF(X)
 \leq\operatorname{Gpd}_B F(X)
 \leq\operatorname{Gpd}_A X.
\]
The inequalities hold in $\mathbb N_0\cup\{\infty\}$, so the equality
also holds at infinite dimension.  The formula for $G$ follows by
symmetry.
\end{proof}

\begin{remark}
The theorem requires two exact functors preserving projectives and
inducing inverse ordinary stable equivalences.  For a single exact
functor inducing a stable equivalence, one must also verify that the
stable inverse has such an exact realization.
\end{remark}

\subsection{Bimodule decompositions over a commutative base}

We next work directly with the bimodule decompositions.  Let $R$ be a
commutative ring, let $A,B$ be $R$-algebras, and suppose that the
$R$-central bimodules $ {}_AM_B$ and $ {}_BN_A$ are finitely generated
projective on each side and satisfy
\begin{equation}
 M\otimes_BN\simeq A\oplus P,
 \qquad
 N\otimes_AM\simeq B\oplus Q.
 \label{eq:gpd-morita-decompositions}
\end{equation}
Put
\[
 F=N\otimes_A-,\quad G=M\otimes_B-,\quad
 E_A=P\otimes_A-,\quad E_B=Q\otimes_B-.
\]
The complement bimodules $P,Q$ are projective on each side, being
direct summands of the displayed tensor products.  Thus all four
functors are exact and preserve projectives.  The next criterion
does not require $P,Q$ to be projective over the enveloping algebras.
We do not assume that $F,G$ induce inverse equivalences of ordinary
stable module categories.

\begin{proposition}
\label{prop:gpd-error-functor-criterion}
Assume \eqref{eq:gpd-morita-decompositions} and suppose that $E_A,E_B$
preserve Gorenstein projective modules in the categories of all
modules.  Then $F,G$ preserve and reflect Gorenstein projectivity, and
\begin{equation}
 \operatorname{Gpd}_B(N\otimes_AX)=\operatorname{Gpd}_A X,\qquad
 \operatorname{Gpd}_A(M\otimes_BY)=\operatorname{Gpd}_B Y
 \label{eq:gpd-bimodule-equalities}
\end{equation}
for every left $A$-module $X$ and left $B$-module $Y$.
\end{proposition}

\begin{proof}
The decompositions give natural isomorphisms
$GF\simeq\mathrm{Id}\oplus E_A$ and
$FG\simeq\mathrm{Id}\oplus E_B$.  They also show that both $F$ and $G$
reflect split epimorphisms.  Indeed, if $G(p)$ has a section, then
$FG(p)$ has a section $s$.  If $p:E\to Z$ and $\iota,\pi$ are the
natural inclusion and projection for the identity summand of $FG$,
then $\pi_Es\iota_Z$ is a section of $p$.

Let $Z$ be a cycle of a totally acyclic complex of projective
$A$-modules.  Then $GF(Z)=Z\oplus E_A(Z)$ is Gorenstein projective.
The extension argument of Theorem~\ref{thm:gpd-exact-stable-pair}
applies: an extension of $F(Z)$ by any projective $B$-module splits
after applying $G$, and hence already splits.  Thus $F$ preserves
Gorenstein projectives.  The same argument applies to $G$.

Since $X$ is a direct summand of $GF(X)$, one has
\[
 \operatorname{Gpd}_A X
 \leq\operatorname{Gpd}_A GF(X)
 \leq\operatorname{Gpd}_B F(X)
 \leq\operatorname{Gpd}_A X.
\]
This proves the first equality, including infinite values, and the
second follows by symmetry.  Dimension zero gives reflection of
Gorenstein projectivity.
\end{proof}

\begin{lemma}
\label{lem:gpd-error-target-criterion}
Let $A$ be a unital ring and let $P$ be an $A$-bimodule which is
projective on both sides.  Suppose that
$\operatorname{Hom}_A(P,U)$ has finite injective dimension as a left
$A$-module for every projective left $A$-module $U$, where
$(a f)(p)=f(pa)$.  Then $P\otimes_A-$ preserves all Gorenstein
projective modules.
\end{lemma}

\begin{proof}
Let $T^\bullet$ be a totally acyclic complex of projective
$A$-modules.  Then $P\otimes_AT^\bullet$ is acyclic with projective
terms.  For any projective left $A$-module $U$, tensor--Hom
adjunction gives
\[
 \operatorname{Hom}_A(P\otimes_AT^\bullet,U)
 \simeq
 \operatorname{Hom}_A(T^\bullet,\operatorname{Hom}_A(P,U)).
\]
An acyclic complex of projectives remains Hom-acyclic against a
module of finite injective dimension.  To see this, use a finite
injective resolution of the target: Hom from each term of the
projective complex is exact, and Hom from an acyclic complex to an
injective module is acyclic.  Induction on the length of the
injective resolution proves the assertion.  Hence the displayed
complex is acyclic for every projective $U$, which is precisely
total acyclicity of $P\otimes_AT^\bullet$.
\end{proof}

\begin{lemma}
\label{lem:gpd-base-injective-error}
Let $R$ be a commutative Noetherian ring and let $A$ be an $R$-algebra
whose underlying $R$-module is injective.  Let $P$ be a finitely
generated projective $A^e$-module which is also projective on both
sides.  Then $P\otimes_A-$ preserves all Gorenstein projective modules.
\end{lemma}

\begin{proof}
Let $U$ be any projective left $A$-module.  Its underlying $R$-module
is a direct summand of a direct sum of copies of $A_R$.  Over a
Noetherian ring, direct sums of injectives are injective: by Baer's
criterion, a map from an ideal into such a sum has image in a finite
subsum because the ideal is finitely generated.  Thus $U_R$ is
injective.

Since $P$ is a direct summand of $(A\otimes_R A)^r$ for a finite $r$,
the module $\operatorname{Hom}_A(P,U)$ is a direct summand of
$\operatorname{Hom}_R(A,U)^r$.  The latter is injective over $A$,
as follows from the coinduction adjunction
\[
 \operatorname{Hom}_A(-,\operatorname{Hom}_R(A,U))
 \simeq \operatorname{Hom}_R(\operatorname{Res}(-),U).
\]
The action is $(a f)(x)=f(xa)$, and the right-hand functor is exact.
Lemma~\ref{lem:gpd-error-target-criterion} now applies.
\end{proof}

\begin{theorem}
\label{thm:gpd-base-injective-semt}
Assume \eqref{eq:gpd-morita-decompositions}.  Suppose that $R$ is
Noetherian, that $A,B$ are injective as $R$-modules, and that $P,Q$
are projective over $A^e,B^e$, respectively.  Then
\eqref{eq:gpd-bimodule-equalities} holds for all left modules.
For all right $A$-modules $U$ and right $B$-modules $V$, one also has
\begin{equation}
 \begin{aligned}
 \operatorname{Gpd}_{B^{\mathrm{op}}}(U\otimes_AM)
   &=\operatorname{Gpd}_{A^{\mathrm{op}}}U,\\
 \operatorname{Gpd}_{A^{\mathrm{op}}}(V\otimes_BN)
   &=\operatorname{Gpd}_{B^{\mathrm{op}}}V.
 \end{aligned}
 \label{eq:gpd-right-bimodule-equalities}
\end{equation}
In particular, the conclusions hold if $R$ is a field, or if $R$ is
a commutative self-injective Artin ring and $A,B$ are finite
projective over $R$.  All assertions include infinite dimension and
require no adjoint-type hypothesis.
\end{theorem}

\begin{proof}
The complements $P,Q$ are finite on each side, and hence finite over
their enveloping algebras.  Apply
Lemma~\ref{lem:gpd-base-injective-error} to $P,Q$, followed by
Proposition~\ref{prop:gpd-error-functor-criterion}.  Applying the same
argument to the opposite algebras proves
\eqref{eq:gpd-right-bimodule-equalities}.  Over a field all modules
are injective; over a self-injective Artin ring finite projective
modules are injective.
\end{proof}

\begin{corollary}
\label{cor:gpd-gorenstein-noetherian-base}
Let $R$ be a commutative Noetherian ring with
$\operatorname{id}_R R<\infty$, and let $A,B$ be finite projective
$R$-algebras.  Suppose that $R$-central bimodules
$ {}_AM_B,{}_BN_A$, finitely generated projective on each side, satisfy
\eqref{eq:gpd-morita-decompositions}, with $P,Q$ projective over
$A^e,B^e$, respectively.  Then the original tensor pair preserves
Gorenstein projective dimensions of all left and right modules:
\eqref{eq:gpd-bimodule-equalities} and
\eqref{eq:gpd-right-bimodule-equalities} hold, including infinite
values.  On finitely generated modules, the same equalities hold
for the left and right Auslander--Bridger $G$-dimensions.
No adjunction or preservation of injective modules is required.
\end{corollary}

\begin{proof}
Put $d=\operatorname{id}_R R$.  Every projective $R$-module has
injective dimension at most $d$.  Indeed, direct sums of injectives
are injective over the Noetherian ring $R$, as in the preceding
lemma.  Taking a direct sum of copies of a length-$d$ injective
resolution of $R$, and then using closure under direct summands,
gives the assertion.

Let $U$ be an arbitrary projective left $A$-module.  Since $A_R$ is
projective, $U_R$ is projective and
$\operatorname{id}_R U\leq d$.  Apply the exact functor
$\operatorname{Hom}_R(A,-)$ to a length-$d$ injective resolution of
$U_R$.  Exactness uses projectivity of $A_R$, and each resulting
term is injective over $A$ by the coinduction adjunction.  Therefore
\[
 \operatorname{id}_A\operatorname{Hom}_R(A,U)\leq d .
\]
The finite $A^e$-projectivity of $P$ implies that
$\operatorname{Hom}_A(P,U)$ is a direct summand of a finite direct
sum of copies of $\operatorname{Hom}_R(A,U)$.  Hence
\[
 \operatorname{id}_A\operatorname{Hom}_A(P,U)\leq d
 \qquad\text{for every projective left $A$-module $U$}.
\]
Lemma~\ref{lem:gpd-error-target-criterion} shows that
$P\otimes_A-$ preserves all Gorenstein projectives.  Repeating the
argument for $B,Q$ and applying
Proposition~\ref{prop:gpd-error-functor-criterion} proves the
left-module equalities.  The opposite algebras are also finite
projective over $R$, so the same proof gives the right-module
equalities.

Finally, $A,B$ are two-sided Noetherian because they are finite over
$R$, and the tensor functors preserve finite modules.  On such
modules Gorenstein projective dimension agrees with
Auslander--Bridger $G$-dimension, with the convention at zero fixed
above.
\end{proof}

\begin{remark}
\label{rem:gpd-ordinary-versus-relative}
Over a field, $E_A,E_B$ are projective-valued, so the tensor pair also
induces inverse equivalences of ordinary stable module categories.
Over a general Noetherian base, that conclusion is not part of
Theorem~\ref{thm:gpd-base-injective-semt} or
Corollary~\ref{cor:gpd-gorenstein-noetherian-base}.  The proofs use
the actual bimodule retracts and preservation of Gorenstein
projectives by the error functors.  They do not identify the
ordinary and relative stable module categories.
\end{remark}

\begin{remark}
\label{rem:gpd-finite-target-test}
For finite-module statements over two-sided Noetherian rings, it is
enough that $\operatorname{Hom}_A(P,A)$ and
$\operatorname{Hom}_B(Q,B)$ have finite injective dimension, with
the analogous conditions over the opposite rings for right modules.
Indeed, use totally acyclic complexes of finite projectives; their
Hom complexes commute with direct sums, so testing the transformed
complex against the regular module suffices to test against every
projective module.  For arbitrary modules we instead use the
projective-target condition in
Lemma~\ref{lem:gpd-error-target-criterion}.  The proofs above verify
that condition for every projective target, without a finite
generation assumption on the target.
\end{remark}

\subsection{Almost \texorpdfstring{$\nu$}{nu}-stable derived equivalences over Artin rings}

\begin{theorem}
\label{thm:gpd-almost-nu}
Let $A$ and $B$ be Artin algebras over an arbitrary commutative Artin
ring.  Let
\[
 F:D^b(A\text{-}\mathrm{mod})\longrightarrow D^b(B\text{-}\mathrm{mod})
\]
be an almost $\nu$-stable derived equivalence, with the usual
nonnegative normalization, and let $\overline F$ be its stable
functor.  Then, for every finitely generated left $A$-module $X$,
\[
 \operatorname{Gpd}_B\overline F(X)
   =\operatorname{Gpd}_A X
 \quad\text{in }\mathbb N_0\cup\{\infty\}.
\]
Equivalently, the same equality holds for Auslander--Bridger
$G$-dimension.  No projectivity over the base ring and no tensor
realization of $\overline F$ are needed.
\end{theorem}

\begin{proof}
By \cite[Theorem~3.7]{HuXi2010}, $\overline F$ is an equivalence of
the entire ordinary stable module categories of $A$ and $B$.
By \cite[Corollary~5.6]{HuPan2017}, its restriction is a triangle
equivalence between their stable categories of finitely generated
Gorenstein projective modules.  Finally,
\cite[Corollary~4.12]{HuPan2017} gives
$\overline F\Omega_A\simeq\Omega_B\overline F$.
Proposition~\ref{prop:gpd-stable-criterion} proves the assertion.
\end{proof}

Passing to opposite algebras gives the following right-module version.

\begin{corollary}
\label{cor:gpd-almost-nu-opposite}
Under the hypotheses of Theorem~\ref{thm:gpd-almost-nu}, let $G$ be a
quasi-inverse of $F$, and write $D_A,D_B$ for the standard Artin
dualities, extended to derived categories.  The derived equivalence
\[
 H=D_A G D_B:
 D^b(B^{\mathrm{op}}\text{-}\mathrm{mod})
 \longrightarrow D^b(A^{\mathrm{op}}\text{-}\mathrm{mod})
\]
is almost $\nu$-stable.  Its stable functor satisfies
\[
 \operatorname{Gpd}_{A^{\mathrm{op}}}\overline H(V)
   =\operatorname{Gpd}_{B^{\mathrm{op}}}V
\]
for every finitely generated right $B$-module $V$.  Equivalently,
any quasi-inverse of $\overline H$ gives a right-module
correspondence from $A$ to $B$ preserving individual $G$-dimensions.
\end{corollary}

\begin{proof}
Let $T^\bullet=G(B)$ and $\overline T^\bullet=F(A)$ be the radical
tilting complexes used in the definition of almost $\nu$-stability.
They are concentrated in nonpositive and nonnegative degrees,
respectively.  Compatibility of a derived equivalence with the
Nakayama functor on perfect complexes, as in
\cite[Lemma~2.3 and the discussion preceding it]{HuXi2010}, identifies
the tilting complexes for $H$ and its inverse with
\[
 (\overline T^\bullet)^*
   =\operatorname{Hom}_B(\overline T^\bullet,B),
 \qquad
 (T^\bullet)^*=\operatorname{Hom}_A(T^\bullet,A).
\]
Their degrees have the required signs.  If $P$ is projective, then
$P^*=D\nu P$ and $\nu_{A^{\mathrm{op}}}(P^*)=DP$.
Consequently, equality of the additive closures of $P$ and $\nu P$
implies equality of those of $P^*$ and
$\nu_{A^{\mathrm{op}}}(P^*)$.  Applying this observation to the
nonzero-degree terms of the two tilting complexes shows that $H$ is
almost $\nu$-stable.  The preceding theorem applied to the opposite
algebras gives the equality.
\end{proof}

\begin{remark}
The right-module correspondence in
Corollary~\ref{cor:gpd-almost-nu-opposite} is defined by the stable
functor of $D_A G D_B$ and its quasi-inverse.  This construction
requires neither a tensor realization of $\overline F$ nor standardness
of the given derived equivalence.
\end{remark}

\section{Rigidity of ordinary stable tensor equivalences}
\label{sec:ordinary-stable-rigidity}

Let $R$ be a commutative ring. We use the ordinary stable
categories, obtained by factoring out morphisms through projective
modules. We shall give conditions under which an $R$-central
tensor functor inducing a stable equivalence is a Morita
equivalence.

\subsection{Rigidity of a tensor realization}

\begin{lemma}
\label{lem:rigidity-special-fibre-hom}
Let $(R,\mathfrak m)$ be a commutative local ring which is not a
field, and let $A$ be a finite projective $R$-algebra.  If $X$ and $Y$
are finitely generated $A$-modules annihilated by $\mathfrak m$, then
every map $X\to Y$ factoring through a projective $A$-module is zero.
Consequently
\[
 \underline{\operatorname{Hom}}_A(X,Y)
   =\operatorname{Hom}_A(X,Y).
\]
\end{lemma}

\begin{proof}
Every finitely generated projective $A$-module $P$ is projective, hence
free, over the local ring $R$.  The image of any map $X\to P$ lies in
\[
 \operatorname{Ann}_P(\mathfrak m)
   =\operatorname{Soc}(R)P\subseteq\mathfrak mP.
\]
The inclusion holds because $R$ is local and not a field: an element
of $\operatorname{Soc}(R)$ outside $\mathfrak m$ would be a unit
annihilating $\mathfrak m$, forcing $\mathfrak m=0$.  Every map
$P\to Y$ kills $\mathfrak mP$, so the composite is zero.
\end{proof}

\begin{theorem}
\label{thm:rigidity-nonfield-artin-base}
Let $(R,\mathfrak m)$ be a commutative local ring which is not a
field.  Let $A$ and $B$ be nonzero finite projective $R$-algebras,
and let ${}_BN_A$ be an $R$-central bimodule such that ${}_BN$ is
finitely generated projective.  Suppose that the tensor functor
\[
 F=N\otimes_A-:A\text{-}\mathrm{mod}
       \longrightarrow B\text{-}\mathrm{mod}
\]
induces an equivalence of ordinary stable module categories.
Then ${}_BN$ is a projective generator, the right-action map
\[
 \rho:A^{\mathrm{op}}\longrightarrow\operatorname{End}_B(N)
\]
is an isomorphism, and $N\otimes_A-$ is a Morita equivalence of
module categories.  In particular, $N_A$ is finitely generated projective.
\end{theorem}

\begin{proof}
Set $k=R/\mathfrak m$, $A_0=A/\mathfrak mA$,
$B_0=B/\mathfrak mB$, and $N_0=N/\mathfrak mN$.
Since ${}_BN$ is projective, $F$ preserves projectives and hence
induces a stable functor.  Exactness of $F$ is not required here.  Centrality over $R$ implies that $F$ takes modules
annihilated by $\mathfrak m$ to modules annihilated by $\mathfrak m$,
and its restriction to these modules identifies with
\[
 N_0\otimes_{A_0}-:A_0\text{-}\mathrm{mod}
                       \longrightarrow B_0\text{-}\mathrm{mod}.
\]
By Lemma~\ref{lem:rigidity-special-fibre-hom} for both $A$ and $B$,
the restriction to this subcategory is fully faithful.  Apply it to the left regular module $A_0$.  Since
$F(A_0)=N_0$, it gives the isomorphism induced by the right action
\begin{equation}
 A_0^{\mathrm{op}}\xrightarrow{\sim}
       \operatorname{End}_{B_0}(N_0).
 \label{eq:rigidity-residue-endomorphisms}
\end{equation}

Put $C=\operatorname{End}_B(N)$.  Choose a finite projective
complement of ${}_BN$ in ${}_BB^s$.  Then $C$ is a corner of
$\operatorname{End}_B(B^s)$, and hence is a direct summand as an
$R$-module of a matrix algebra over $B^{\mathrm{op}}$.  Thus $C$ is
finite projective, and therefore finite free, over $R$.  The same
corner description shows that
\[
 C/\mathfrak mC\simeq\operatorname{End}_{B_0}(N_0).
\]
Under this identification, the reduction of $\rho$ is precisely
\eqref{eq:rigidity-residue-endomorphisms}.  Its source and target are finite free $R$-modules of the same rank,
since their reductions have the same $k$-dimension.  The reduction
of its representing square matrix is invertible, so its determinant
is a unit in $R$.  Hence $\rho$ is an isomorphism.

It remains to show that $N$ is a generator.  Let $S$ be a simple
left $B$-module.  The module $S$ is finite over $R$, since it is a
quotient of $B$.  Centrality implies that $\mathfrak mS$ is a
$B$-submodule of $S$, and Nakayama's lemma excludes
$\mathfrak mS=S$.  Hence $\mathfrak mS=0$.
Such a nonzero $S$ cannot be projective over
$B$: a projective $B$-module is free over $R$, whereas a nonzero
free $R$-module is not annihilated by the nonzero ideal
$\mathfrak m$.

Stable essential surjectivity supplies $X\in A\text{-}\mathrm{mod}$
and a stable isomorphism $S\simeq F(X)$.  Choose representatives
$i:S\to F(X)$ and $p:F(X)\to S$ of mutually inverse stable maps.
The endomorphism $pi-1_S$ factors through a projective $B$-module,
and is zero by Lemma~\ref{lem:rigidity-special-fibre-hom}.
Therefore $pi=1_S$: the module $S$ is an actual direct summand of
$F(X)$.

Choose an epimorphism $A^r\twoheadrightarrow X$ with $r$ finite.
Right exactness of tensoring yields an epimorphism
$N^r\twoheadrightarrow F(X)$, and composition with $p$ gives
$N^r\twoheadrightarrow S$.  Thus
$\operatorname{Hom}_B(N,S)\ne0$ for every simple $B$-module $S$.

Let
\[
 I=\operatorname{Tr}_B(N)
   =\sum_{f\in\operatorname{Hom}_B(N,B)}f(N).
\]
If $B/I\ne0$, choose a simple quotient $S$ of $B/I$ and a nonzero
map $N\to S$.  Projectivity of $N$ lifts this map to $B/I$ and
then to $B$.  The latter lift has image in $I$, so its composite
with $B\to B/I$ is zero, a contradiction.  Therefore $I=B$.
Expressing $1_B$ as a finite sum of elements $f(n)$ gives a split
epimorphism from a finite direct sum of copies of $N$ onto $B$.
Hence ${}_BN$ is a projective generator.

Morita theory now applies to the projective generator ${}_BN$ and
the isomorphism
$A^{\mathrm{op}}\simeq\operatorname{End}_B(N)$, proving the assertion.
\end{proof}

\begin{remark}
Theorem~\ref{thm:rigidity-nonfield-artin-base} applies to a specified
$R$-central tensor functor on the ordinary stable categories of all
finitely generated modules.  An equivalence on relative stable
categories or on lattice categories does not satisfy this hypothesis.
Related results are proved in
\cite[Theorem~8.8 and Corollary~8.9]{DugasMartinezVilla2008} under
different assumptions.  Their domain corollary concerns torsionfree
algebras over a non-field domain and uses the vanishing of the socle,
whereas the proof above applies to every commutative local
non-field base, including bases with nonzero socle.
\end{remark}

\begin{corollary}\label{cor:rigidity-global}
Let $R$ be a commutative ring such that $R_{\mathfrak m}$ is not
a field at any maximal ideal $\mathfrak m$. Let $A,B$ be finite
projective $R$-algebras, and $N$ a central $B$--$A$-bimodule
finitely generated projective on the left. If $N\otimes_A-$
induces an equivalence of ordinary stable categories of
finitely generated modules, then it is a Morita equivalence.
For the proof, full faithfulness on the stable morphisms
between the residue regular modules and themselves, and stable
essential surjectivity on simple $B$-modules, suffice.
\end{corollary}

\begin{proof}
Fix a maximal ideal $\mathfrak m$. For modules annihilated by
$\mathfrak m$, a map factoring through an ordinary projective
is zero after localization by
Lemma~\ref{lem:rigidity-special-fibre-hom}. Such localization
leaves the source and target unchanged, so the map is zero.
Consequently stable full faithfulness on $A/\mathfrak m A$
identifies the right-action isomorphism on residue endomorphism
rings. The corner argument in
Theorem~\ref{thm:rigidity-nonfield-artin-base} identifies these
endomorphism rings with the residue fibres of
$\rho:A^{\mathrm{op}}\to\operatorname{End}_B(N)$.
Thus at every maximal ideal $\rho$ is an isomorphism of finite
free modules, by the determinant argument after localization.
It is therefore an isomorphism over $R$.

Every simple $B$-module $S$ is finite over $R$ and is annihilated
by some maximal ideal $\mathfrak m$: choose one in its support
and use Nakayama at $\mathfrak m$ to exclude $\mathfrak m S=S$.
The stable essential surjectivity argument of
Theorem~\ref{thm:rigidity-nonfield-artin-base} now makes $S$
a direct summand of some $N\otimes_A X$, because the projective
factorization ideal on $\operatorname{End}_B(S)$ is zero.
It follows that $\operatorname{Hom}_B(N,S)\ne0$ for all such $S$.
The trace-ideal argument proves that $N$ is a generator,
and Morita theory finishes the proof.
\end{proof}

For example, the base condition holds for any nonfield
integral domain. The argument permits bases with nonzero socle
in their localizations.

\subsection{Bimodule decompositions over a non-field base}

\begin{example}
\label{ex:rigidity-raw-semt-product}
Let $R=\mathbb Z/4\mathbb Z$, $A=R$, and $B=R\times R$.  Give
\[
 {}_AM_B={}_RB_B,\qquad {}_BN_A={}_BB_R
\]
the scalar $R$-actions induced by the diagonal map $R\to B$.
Both bimodules are projective on each side.  Moreover,
\[
 M\otimes_BN\simeq B\simeq A\oplus P,
 \qquad P=R\in A^e\text{-}\mathrm{proj}.
\]
If $e_1,e_2$ are the two central primitive idempotents of $B$, then
\[
 N\otimes_AM\simeq B\otimes_RB\simeq B\oplus Q,
 \qquad
 Q=(Be_1\otimes_Re_2B)\oplus(Be_2\otimes_Re_1B).
\]
The diagonal summands give $B$, and the two off-diagonal summands
give the projective $B^e$-module $Q$.  Thus $M,N$ satisfy the bimodule conditions for a stable equivalence
of Morita type.

Nevertheless, the tensor functors on modules are
\[
 X\longmapsto(X,X),\qquad (Y_1,Y_2)\longmapsto Y_1\oplus Y_2.
\]
For $S=R/(2)$, Lemma~\ref{lem:rigidity-special-fibre-hom} gives
\[
 \underline{\operatorname{End}}_A(S)=\mathbb F_2,
 \qquad
 \underline{\operatorname{End}}_B((S,S))
       =\mathbb F_2\times\mathbb F_2.
\]
The diagonal functor induces the map $a\mapsto(a,a)$
on these rings. Its image has dimension one over $\mathbb F_2$.  Thus this functor is not
an equivalence of ordinary stable module categories.  Also $A$
and $B$ are not Morita equivalent: they have one and two simple
modules, respectively.

In contrast, every $A$-module is relatively $R$-projective.  Every $B$-module
$(Y_1,Y_2)$ is a direct summand of
$B\otimes_R(Y_1\oplus Y_2)$.  Thus both $R$-relative stable
categories are zero, and the displayed tensor pair induces an
equivalence between them.
\end{example}

\subsection{An almost \texorpdfstring{$\nu$}{nu}-stable functor without an ordinary tensor realization}

\begin{example}
\label{ex:rigidity-almost-nu-z8}
Let $R=\mathbb Z/8\mathbb Z$ and $A=R\times R$.  Under the product
identification of derived categories, consider
\[
 F=\operatorname{Id}\times[-1]:D^b(A\text{-}\mathrm{mod})
                          \longrightarrow D^b(A\text{-}\mathrm{mod}).
\]
The tilting complexes associated with $F$ and its inverse are
\[
 T=(R,0)\oplus(0,R)[1],\qquad
 \overline T=(R,0)\oplus(0,R)[-1].
\]
They have the required nonpositive and nonnegative supports.
The ring $R$ is self-injective and its Nakayama functor fixes its
indecomposable projective.  Hence the nonzero-degree terms satisfy
the almost $\nu$-stable conditions of \cite{HuXi2010}.
By \cite[Example~4.7(b)]{HuPan2017}, the associated stable functor is
\[
 \overline F=\operatorname{Id}\times\Omega_R.
\]

Put $S=R/(2)$.  The projective cover $R\twoheadrightarrow S$ has
kernel $2R\simeq R/(4)$, and therefore
\[
 \overline F(0,S)\simeq(0,R/(4)).
\]
For any finite $R$-central $A$-bimodule $N$, the module
$N\otimes_A(0,S)$ is annihilated by $2$.  It consequently has the
form $(S^a,S^b)$ for some $a,b\geq0$.  Such a module cannot be
stably isomorphic to $(0,R/(4))$: the indecomposable nonprojective
$R$-modules $R/(2)$ and $R/(4)$ are nonisomorphic, and
Krull--Schmidt uniqueness remains valid after adjoining projective
summands.  Thus $\overline F$ has no realization by tensoring with
an $R$-central bimodule projective on the target side, even at the
level of isomorphism classes of stable objects.

Thus the specified stable functor has no such tensor realization,
although the identity bimodule gives a Morita-type self-equivalence
of $A$.
\end{example}

\begin{remark}
A Morita realization of a stable functor need not determine its
underlying derived equivalence.  Replace $\mathbb Z/8\mathbb Z$ in the preceding construction by
$R=\mathbb Z/4\mathbb Z$.  Every finite $R$-module is a direct sum
of copies of $R$ and $S=R/(2)$, and
$R\text{-}\underline{\mathrm{mod}}$ identifies with the category
of finite-dimensional $\mathbb F_2$-vector spaces.  Under this
identification $\Omega_R$ is naturally isomorphic to the identity:
it fixes $S$ and acts identically on its matrix endomorphisms.
Thus the stable functor of $\operatorname{Id}\times[-1]$ is
realized by the identity Morita bimodule.  The derived functor
itself is not $t$-exact, and is not a common shift of a Morita
equivalence, since its two factors have different shifts.
\end{remark}

\section{Relative Frobenius categories and derived equivalences}
\label{sec:relative-frobenius-derived}

We recall the $R$-split Frobenius structure needed for the
construction of Rouquier and Zimmermann. We then study the
Gorenstein condition of Abe and Hoshino for endomorphism
algebras of tilting complexes.

\subsection{The relative Frobenius structure}

Put $D_RA=\operatorname{Hom}_R(A,R)$.

\begin{proposition}\label{prop:relative-frobenius-criterion}
Assume that $A$ is finitely generated projective over $R$ and
$D_RA$ is a projective generator as a left $A$-module.  Then
\[
 \mathcal P_R(A)=\mathcal I_R(A).
\]
Consequently the $R$-split exact category of finitely generated
left $A$-modules is Frobenius, and
$A\text{-}\underline{\mathrm{mod}}_R$ has its canonical stable
triangulated structure.  The analogous right-module statement
holds when $(D_RA)_A$ is a projective generator.
\end{proposition}

\begin{proof}
Since $A_R$ is finitely generated projective, evaluation gives a
natural left $A$-module isomorphism
\[
 D_RA\otimes_RV\simeq\operatorname{Hom}_R(A,V)
 \qquad (V\in R\text{-}\mathrm{mod}).
\]
Since ${}_A(D_RA)$ is projective, every module on the right belongs
to $\mathcal P_R(A)$.  Conversely,
$A\in\operatorname{add}({}_A D_RA)$ implies
$A\otimes_RV\in\operatorname{add}(D_RA\otimes_RV)$, hence
$A\otimes_RV\in\mathcal I_R(A)$.  Taking direct summands gives
equality of the two classes.  Since
the exact category has enough projectives and enough injectives,
it is Frobenius.  Relative syzygy and cosyzygy therefore induce
mutually inverse shift functors on its stable category.
\end{proof}

\subsection{The Rouquier--Zimmermann construction}

The following theorem gives the derived-to-relative construction
used over Noetherian bases.  Its bimodule assertion is due to
Rouquier and Zimmermann.

\begin{theorem}
\label{thm:rz-relative-derived-stable}
Suppose $A$ and $B$ are finite projective $R$-algebras and
$D_RA$ and $D_RB$ are projective generators on both sides.
Then every standard $R$-linear derived equivalence between $A$ and $B$
yields bimodules ${}_AM_B$ and ${}_BN_A$, finitely generated projective on both
sides, satisfying \eqref{eq:relative-stabilized-mt}.  In
particular, $N\otimes_A-$ and $M\otimes_B-$ induce mutually
inverse triangle equivalences
\[
 A\text{-}\underline{\mathrm{mod}}_R
   \simeq B\text{-}\underline{\mathrm{mod}}_R
\]
on the full relative stable categories of finitely generated
modules, with their canonical triangulated structures from
Proposition~\ref{prop:relative-frobenius-criterion}.
\end{theorem}

\begin{proof}
The bimodule assertion follows from
\cite[Section~2.6, Proposition~2.14 and Corollary~2.15]{RouquierZimmermann2003}.
That section assumes finite projectivity over the base and uses
the term Gorenstein over the base for an algebra whose
base dual is projective.  The present hypotheses imply these
conditions.  For inverse two-sided tilting complexes, the cited
results construct the required one-sided projective bimodules
and the identities with ordinary projective bimodules added to
both sides.  In particular, the error terms are projective over
the ordinary enveloping algebras $A^e$ and $B^e$ defined over $R$.

Proposition~\ref{prop:relative-morita-bimodules} now gives the
additive equivalence on all finitely generated modules in the
relative quotients.  Proposition~\ref{prop:relative-frobenius-criterion}
shows that both exact categories are Frobenius.  The tensor
functors are exact for the $R$-split structures and preserve
projective-injective objects.  Their induced stable functors
are therefore mutually inverse triangle equivalences.
\end{proof}

\begin{corollary}\label{cor:rz-lattices-and-symmetric-algebras}
Under the hypotheses of
Theorem~\ref{thm:rz-relative-derived-stable}, let $\mathcal L_R(A)$
denote the full subcategory of finite $A$-modules which are
projective over $R$, and define $\mathcal L_R(B)$ similarly.
The constructed tensor functors restrict to mutually inverse
triangle equivalences
\[
 \mathcal L_R(A)/[A\text{-}\mathrm{proj}]
   \simeq \mathcal L_R(B)/[B\text{-}\mathrm{proj}].
\]
In particular, the theorem and this corollary apply to finite
projective symmetric $R$-algebras $A,B$, that is, algebras with
$D_RA\simeq A$ and $D_RB\simeq B$ as bimodules.
\end{corollary}

\begin{proof}
The one-sided projectivity of the bimodules implies that the
tensor functors preserve $R$-projective modules.  For
$X\in\mathcal L_R(A)$, the canonical $R$-split epimorphism
$A\otimes_RX\to X$ has an $A$-projective middle term and an
$R$-projective kernel.  The canonical $R$-split monomorphism
$X\to D_RA\otimes_RX$ has an $A$-projective middle term and an
$R$-projective cokernel.  The base-dual hypothesis therefore
makes $\mathcal L_R(A)$ a Frobenius exact category whose
projective-injective objects are the finite projective
$A$-modules.  The same holds for $B$.

If $P$ is projective over $A^e$ and $X$ is projective over $R$,
then $P\otimes_AX$ is projective over $A$.  Hence the stabilized
identities give inverse stable equivalences on these
subcategories.  The tensor functors are exact and preserve
projective-injective objects, so these are triangle equivalences.
The last assertion follows from the displayed bimodule
isomorphisms for symmetric algebras.
\end{proof}

\begin{corollary}\label{cor:rz-henselian-adjoint-type}
Assume the hypotheses of
Theorem~\ref{thm:rz-relative-derived-stable}.  Suppose further
that $R$ is Henselian local and that $A,B$ have no
$R$-separable direct factors.  Then the constructed bimodules
have reduced representatives $M_0,N_0$ which are projective on
both sides and satisfy
\[
 M_0\otimes_BN_0\simeq A\oplus P,
 \qquad N_0\otimes_AM_0\simeq B\oplus Q,
 \qquad P\in A^e\text{-}\mathrm{proj},\quad
        Q\in B^e\text{-}\mathrm{proj}.
\]
Their tensor functors are adjoint in both orders and induce the
same equivalences of the $R$-relative stable categories as the
original pair.
\end{corollary}

\begin{proof}
Apply Proposition~\ref{prop:adjoint-complete-base-extension} and
Corollary~\ref{cor:adjoint-relative-stable-identities} to the
ordinary projectivity class and the stabilized identities
from Theorem~\ref{thm:rz-relative-derived-stable}.
\end{proof}

\begin{remark}
Rouquier and Zimmermann use the stabilized identities
\eqref{eq:relative-stabilized-mt} in the definition preceding
their Corollary~2.15.  These identities suffice for
Theorem~\ref{thm:rz-relative-derived-stable} over an arbitrary
commutative Noetherian base satisfying its hypotheses.
Corollary~\ref{cor:rz-henselian-adjoint-type} gives the
direct-summand form \eqref{eq:relative-raw-mt} under the stated
additional assumptions.  Neither assertion requires an
equivalence of the entire ordinary stable module categories.
\end{remark}

\subsection{The question of Abe and Hoshino}
\label{subsec:abe-hoshino-answer}

In this subsection we use right modules. For a commutative Noetherian
ring $R$, a finite $R$-algebra $C$ is Gorenstein of codimension $c$
in the sense of \cite[Definition~2.1]{AbeHoshino2009} if
$R_{\mathfrak p}$ is Gorenstein for every
$\mathfrak p\in\operatorname{Supp}_R C$, the groups
$\operatorname{Ext}^i_R(C,R)$ vanish for $i\ne c$, and
$\operatorname{Ext}^c_R(C,R)$ is a projective generator as a right
$C$-module. After Example~4.6, Abe and Hoshino ask whether the
Ext vanishing for a derived equivalent algebra is sufficient for
it to be Gorenstein of the same codimension. We give an affirmative
answer.

\begin{theorem}\label{thm:abe-hoshino-answer}
Let $A$ be a Gorenstein $R$-algebra of codimension $c$ in the sense
above. Let $P\in K^b(\mathrm{proj}\text{-}A)$ be a tilting complex and
put $B=\operatorname{End}_{K^b(\mathrm{proj}\text{-}A)}(P)$, with the
induced $R$-algebra structure. If
\[
 \operatorname{Ext}^i_R(B,R)=0\qquad(i\ne c),
\]
then $B$ is a Gorenstein $R$-algebra of codimension $c$.
\end{theorem}

We first establish a self-injectivity result that will be used after
reduction modulo a system of parameters. The case of
finite-dimensional algebras over an arbitrary field is known
\cite[Corollary~5.3]{XiZhangSelfInjective2025}.

\begin{lemma}\label{lem:artin-self-injective-base}
Let $S$ be a commutative Noetherian ring and let $C$ be a finite
$S$-algebra. If $C$ is self-injective, then the opposite endomorphism
algebra of every tilting complex in $K^b(C\text{-}\mathrm{proj})$
is self-injective.
\end{lemma}

\begin{proof}
Put $S_0=S/\operatorname{Ann}_S C$. Since $C$ is finite over $S$,
it is Noetherian. By the characterization of quasi-Frobenius rings
\cite[Theorem~15.1]{Lam1999ModulesRings}, self-injectivity implies
that $C$ is Artinian. Thus the left regular module ${}_CC$ has a
finite composition series.

Let $V$ be one of its simple factors, and choose a maximal ideal
$\mathfrak m$ in $\operatorname{Supp}_S V$. Since $V$ is finite over
$S$, Nakayama's lemma at $\mathfrak m$ gives $\mathfrak m V\ne V$.
The central action of $S$ makes $\mathfrak m V$ a $C$-submodule of
$V$. Hence $\mathfrak m V=0$ by simplicity, and
$\operatorname{Supp}_S V=\{\mathfrak m\}$. It follows that
\[
 \operatorname{Spec} S_0
 = V(\operatorname{Ann}_S C)
 = \operatorname{Supp}_S C
\]
consists of finitely many maximal ideals. Therefore $S_0$ is
Artinian.

Write $S_0=\prod_i S_i$ as a product of Artinian local rings.
The central idempotents of $S_0$ decompose $C$ and the tilting
complex accordingly. The opposite endomorphism algebra is the
product of those of the resulting tilting complexes, and
self-injectivity can be checked on each factor.

For each $i$, the Cohen structure theorem \cite{StacksCohen}
gives a surjection $Q_i\twoheadrightarrow S_i$ from a complete
regular local ring. Let $x_1,\ldots,x_d$ be a regular system of
parameters of $Q_i$. Since the maximal ideal of $S_i$ is
nilpotent, for sufficiently large $N$ the surjection factors
through
\[
 Q_i/(x_1^N,\ldots,x_d^N).
\]
This quotient is an Artinian local complete intersection, hence
Gorenstein. Restricting the action of $S_i$ along this surjection
leaves the algebra, its tilting complex and the opposite
endomorphism algebra unchanged. We may therefore assume that
$S$ is an Artinian local Gorenstein ring.
Write $D=\operatorname{Hom}_S(-,S)$. This is a duality on
finite $S$-modules. For a finite $S$-algebra $L$, the module
$DL$ is an injective cogenerator on either side. It is projective
if and only if $L$ is self-injective. Indeed, if it is projective,
the indecomposable injectives exhaust the isomorphism classes of
indecomposable projectives, since both sets are indexed by the
simple modules. In that case $DL$ is also a generator.

First suppose that $C$ is a finite ring, and let $E$ be the
endomorphism algebra of a tilting complex $T$ over $C$. We write
$D^b(L)$ for the bounded derived category of finite right
$L$-modules and put $\nu_L=-\otimes_L^{\mathbf L}DL$. On perfect complexes, the
$S$-linear equivalence afforded by $T$ identifies $\nu_C$ and
$\nu_E$. To see this, for a perfect complex $U$ and any bounded
complex $V$ of finite $L$-modules, the duality is given by
\[
 D\operatorname{Hom}_{D^b(L)}(U,V)
       \simeq\operatorname{Hom}_{D^b(L)}(V,\nu_L U).
\]
It follows by computing with a bounded projective representative of
$U$ and its termwise Nakayama image, a bounded injective complex.
The equivalence on perfect complexes is $S$-linear, as can be
seen from the dg endomorphism algebra of $T$ and its truncation
to $E$. Thus for perfect $V$ the object transported from the
$C$-side and $\nu_E U$ represent the same functor of $V$.
The transported object is perfect. Evaluating the natural
isomorphism at this object gives a morphism to $\nu_E U$,
which is an isomorphism since it is detected by the shifts of $E$.
In particular $\nu_E$ is a self-equivalence of the perfect
category.

The exact functor $\nu_C$ permutes the indecomposable finite
projectives. Consequently the iterates of $T$ under $\nu_C$ have
only finitely many possible graded terms up to isomorphism.
Since the homomorphism sets between finite $C$-modules are finite,
the iterates have only finitely many possible differentials on
these terms. Hence $\nu_E^s E\simeq E$ for some $s>0$.
We use the backward argument of
\cite[Theorem~2.1]{AlNofayeeRickard2013}.
Every $\nu_E^j E$, $j\geq0$, is perfect with vanishing positive
cohomology, since $\nu_E$ is derived tensor with a bimodule in
degree zero. If $\nu_E^{j+1} E$ is projective in degree zero,
take a minimal bounded projective representative of $\nu_E^j E$.
It has no positive terms, since otherwise the last differential
would be a split epimorphism. If its lowest nonzero term is in
negative degree, the lowest differential after applying $\nu_E$
termwise is injective because the image has zero cohomology
in that degree. It then splits since its source is injective.
This contradicts minimality, as Nakayama gives an equivalence
from finite projectives to finite injectives. Descending from
$j=s-1$ proves that $DE=\nu_E E$ is projective.

Next suppose the residue field of $S$ has characteristic zero.
A coefficient field makes $C$ and $E$ finite-dimensional
algebras over a field, so the field result applies. In residue
characteristic $p>0$ we first make a faithfully flat local
extension with algebraically closed residue field $k$. Such an
extension is obtained by adjoining successive roots of monic lifts
of irreducible residue polynomials and taking filtered unions.
At each step the maximal ideal is the extension of the original
maximal ideal $\mathfrak m$ and the layers of the
$\mathfrak m$-adic filtration are extended residue vector spaces.
The resulting ring is Artinian since the filtration is finite
and these layers are finite-dimensional over the new residue
field. It is Gorenstein since the socle commutes with this flat
extension and has dimension one. The endomorphism complex and
the dual $D$ commute with flat extension. Projectivity of a finite
module descends along it, and generation of the perfect category
ascends by extension of the finite construction of $C$ from $T$.

By the Cohen structure theorem the extended local ring is finite
over a truncated Witt ring $W_n(k)$
\cite{StacksCohen,SerreLocalFields}. We regard the two algebras
and the tilting complex as data over $W_n(k)$. Its own duality
$\operatorname{Hom}_{W_n(k)}(-,W_n(k))$ also detects
self-injectivity of a finite algebra by the first paragraph.
We now show that a counterexample over $W_n(k)$ would give one
over $W_n(k_0)$ for a finite field $k_0$.

Fix $n$, the cyclic module decompositions of $C$ and $E$ over
$W_n(k)$, and the finite ranks and degrees of the idempotent
matrices presenting the terms of $T$. Also fix the shape of a
finite construction of $C$ from $T$ by shifts, sums, cones and
direct summands. Each intermediate complex can be represented
by bounded finite projectives, and the isomorphisms or
retractions in the homotopy category are witnessed by chain
maps and homotopies. These witnesses involve only finitely
many additional matrices. With these choices the existence
of the counterexample is a first-order sentence in the
language of fields of characteristic $p$. More explicitly,
a cyclic module $W_n(k)/(p^i)$ is interpreted as the quotient
by the definable relation $u-v\in p^iW_n(k)$. Algebra structures
are specified by structure constants. Exactness in each
nonzero degree of the bounded endomorphism complex is expressed
by requiring every cycle to be a boundary. The identification
of $E$ with its zeroth cohomology is specified by cycles for
the module generators of $E$, with the module and ring relations
holding modulo boundaries, and by injectivity and surjectivity
of the resulting map. The dual of either algebra consists
of the tuples of values of linear functions on its cyclic
generators, subject to the corresponding annihilator
conditions. For a module of this type, the number of cyclic
generators bounds the number of generators of its dual over the
Witt ring and hence over the algebra. Projectivity of the dual
is therefore expressed by the existence of maps exhibiting it
as a summand of a free module of that fixed rank. To specify
a counterexample, require this splitting on the $C$-side and
its failure on the $E$-side. Witt addition and multiplication
are given by universal polynomials, so these are field formulas.

Completeness of the theory of algebraically closed fields of
characteristic $p$ now gives such data over $\overline{\mathbb F}_p$
\cite{BodirskyModelTheory}. Choose a finite subfield $k_0$
containing the coordinates of all structure constants, of the
complex and of the finite generation witnesses. The resulting
Witt ring extension is faithfully flat. Indeed, a residue
field basis lifts to a basis over $W_n(k_0)$ by the finite
$p$-adic filtration. Hence the endomorphism complex over
$W_n(k_0)$ has cohomology only in degree zero and its
endomorphism algebra extends to $E$ over $\overline{\mathbb F}_p$.
The finite generation witnesses descend by the choice of
$k_0$. Since duality commutes with this flat extension,
projectivity on the $C$-side descends, and projectivity
on the $E$-side would ascend. This produces a counterexample
over a finite ring, a contradiction.
\end{proof}

\begin{proof}[Proof of Theorem~\ref{thm:abe-hoshino-answer}]
The complex $\operatorname{Hom}_A^\bullet(P,P)$ commutes with
central localization. Thus
$\operatorname{Supp}_R A=\operatorname{Supp}_R B$:
if $A_{\mathfrak p}=0$, the localized endomorphism algebra
vanishes, and if $B_{\mathfrak p}=0$, the identity of
$P_{\mathfrak p}$ vanishes in the homotopy category, so its
thick closure is zero and $A_{\mathfrak p}=0$.
Projectivity of a finite $B$-module is detected at the maximal
ideals of $R$, by localizing the extension class of a finite
free presentation. The generator property is detected by
localizing the trace ideal. We may therefore assume that $R$ is
local Gorenstein and $A\ne0$.

By \cite[Proposition~2.9]{AbeHoshino2009}, choose an
$R$-regular sequence $\boldsymbol x$ of length $c$ in
$\operatorname{Ann}_R A$. It annihilates $B$ as it annihilates
$P$ termwise. Put $S=R/(\boldsymbol x)$. The change-of-rings
isomorphisms
\[
 \operatorname{Ext}_S^j(C,S)
       \simeq\operatorname{Ext}_R^{c+j}(C,R),
 \qquad C=A,B,
\]
are compatible with the bimodule actions. The two algebras are
maximal Cohen--Macaulay over the Gorenstein local ring $S$,
by local duality. Also $D_SA$ is a projective generator. We will
prove the same assertion for $D_SB$.

Choose a system of parameters of $S$. It is regular on $A$
and $B$. At a single step let $y$ be the next parameter, let
$\bar S=S/(y)$ and write $\bar C=C/yC$. The terms of
$H^\bullet=\operatorname{Hom}_A^\bullet(P,P)$ belong to
$\operatorname{add}(A_S)$. Multiplication by $y$ is injective
on each term. The quotient complex is the endomorphism complex
of $P/yP$ over $\bar A$, by projectivity of the terms of $P$.
Since $y$ is regular on $B=H^0(H^\bullet)$, the long exact
sequence gives $H^i(H^\bullet/yH^\bullet)=0$ for $i\ne0$ and
$H^0(H^\bullet/yH^\bullet)=\bar B$ as algebras. The functor
on bounded projective complexes given by reduction sends
$A$ to $\bar A$ and preserves the triangles and summands
which give generation. Thus $P/yP$ is tilting. In addition,
applying $\operatorname{Hom}_S(C,-)$ to the parameter sequence
gives
\[
 (D_SC)/y(D_SC)\simeq D_{\bar S}\bar C,\qquad C=A,B,
\]
since $\operatorname{Ext}^1_S(C,S)=0$. This identifies the
dual at each reduction, and reduction of $D_SA$ preserves its
projective generator property.

After all the parameters have been factored out, the base is
Artinian local Gorenstein, the reduction of $A$ is self-injective
and the reductions of $A,B$ are derived equivalent.
Lemma~\ref{lem:artin-self-injective-base} shows that the
reduced dual of $B$ is a projective generator.

For ascent through one parameter put $K=D_SB$ and take a finite
free presentation $0\to L\to F\to K\to 0$. Multiplication by
$y$ is injective on $K$ since it is injective on $S$, and on
$L$ since it is injective on $F$. Since $y$ is regular on
$B$ and $K$, reduction of a projective resolution of $K$ is a
projective resolution of $K/yK$. Consequently
\[
 \operatorname{Ext}_B^1(K,L/yL)
  \simeq \operatorname{Ext}_{\bar B}^1(K/yK,L/yL)=0.
\]
The last equality uses the projectivity of the reduced dual.
The sequence for multiplication by $y$ on $L$ gives
$y\operatorname{Ext}_B^1(K,L)=\operatorname{Ext}_B^1(K,L)$.
Nakayama's lemma makes this finite $S$-module zero, so the
presentation splits. Now that $K$ is projective, its trace
ideal reduces to the trace ideal of $K/yK$ by base change
for $\operatorname{Hom}_B(K,B)$. The latter ideal is $\bar B$.
Nakayama's lemma gives the generator property of $K$.
Iterating ascent and using the change-of-rings isomorphism
shows that $\operatorname{Ext}^c_R(B,R)$ is a projective
generator. This proves the theorem.
\end{proof}

\begin{corollary}\label{cor:abe-hoshino-nakayama}
Under the hypotheses of Theorem~\ref{thm:abe-hoshino-answer}, put
$\Omega_A=\operatorname{Ext}_R^c(A,R)$ and
$\nu_A=-\otimes_A^{\mathbf L}\Omega_A$. Then
$\operatorname{add}(\nu_A P)=\operatorname{add}(P)$.
\end{corollary}

\begin{proof}
Apply \cite[Corollary~4.4]{AbeHoshino2009} to the conclusion of
Theorem~\ref{thm:abe-hoshino-answer}.
\end{proof}

\section{Almost \texorpdfstring{$\nu$}{nu}-stable derived equivalences over an Artin base}
\label{anu-sec}

In this section, we extend the bimodule identities for almost
$\nu$-stable derived equivalences to algebras projective over
a commutative Artin ring $R$. All algebras are finite over $R$
and all bimodules are $R$-central. Put
$A_R^e=A\otimes_R A^{\mathrm{op}}$. We say that a pair
${}_A M_B,{}_B N_A$ satisfies the \emph{Morita-type bimodule identities}
if its members are finitely generated projective on each side and
\begin{equation}\label{anu-eq-identities}
 M\otimes_BN\simeq A\oplus P,
 \qquad N\otimes_AM\simeq B\oplus Q,
\end{equation}
where $P$ is $A_R^e$-projective and $Q$ is $B_R^e$-projective.
By Section~\ref{sec:relative-and-abe}, these identities
give equivalences of the relative stable categories.  Over a general
Artin base, they need not give equivalences of the ordinary stable
module categories.

We use the definition of an almost $\nu$-stable equivalence from
Hu--Xi \cite{HuXi2010}.  Thus, for an almost $\nu$-stable derived equivalence
$F:D^b(A\text{-}\mathrm{mod})\longrightarrow
D^b(B\text{-}\mathrm{mod})$ with quasi-inverse $G$, the associated
tilting complexes have radical projective representatives
\[
 T^\bullet\simeq G(B),\quad T^i=0\ (i\notin[-n,0]),
 \qquad
 U^\bullet\simeq F(A),\quad U^i=0\ (i\notin[0,n]).
\]
Set $T^- =\bigoplus_{i<0}T^i$ and
$U^+=\bigoplus_{i>0}U^i$.  The required conditions are
\begin{equation}\label{anu-eq-almost-nu}
 \operatorname{add}(T^-)=\operatorname{add}(\nu_A T^-),
 \qquad
 \operatorname{add}(U^+)=\operatorname{add}(\nu_B U^+).
\end{equation}
Here $\nu_A=D\operatorname{Hom}_A(-,A)$, with $D$ the usual
duality for finite modules over the Artinian base.

\begin{theorem}\label{anu-thm-raw-semt}
Suppose that $A$ and $B$ are finitely generated projective $R$-modules.
If there is an $R$-linear almost $\nu$-stable derived equivalence
between $A$ and $B$, then there is a pair of $R$-central bimodules
satisfying \eqref{anu-eq-identities}.
The assertion applies to any such $R$-linear equivalence.
\end{theorem}

For a field, Hu--Xi \cite[Theorem~5.3, pp.~141--143]{HuXi2010}
prove the stronger assertion identifying the induced ordinary stable
equivalence with a tensor functor.  Their Question~(2), p.~151, asks
about $R$-projective Artin algebras.  Theorem~\ref{anu-thm-raw-semt} establishes the bimodule conclusion
\eqref{anu-eq-identities}.  The proof controls the outer tensor products over $R$ and
then adjoins projective bimodule summands to obtain the two
isomorphisms in \eqref{anu-eq-identities}.

\subsection{Preliminary lemmas}

\begin{lemma}\label{anu-lem-standard-replacement}
Under the hypotheses of Theorem~\ref{anu-thm-raw-semt}, there exists
an $R$-linear standard almost $\nu$-stable derived equivalence between
$A$ and $B$.
Only the existence of a standard replacement is asserted.
\end{lemma}

\begin{proof}
We first establish a property of minimal projective complexes.  Let $H$
be a $t$-exact $R$-linear autoequivalence of $D^b(B\text{-}\mathrm{mod})$,
and write $E$ for its restriction to the heart.  Thus $E$ is a Morita
autoequivalence.  If $P^\bullet$ is a radical bounded complex of
projectives and $V^\bullet$ is the radical projective representative
of $H(P^\bullet)$, then
\begin{equation}\label{anu-eq-minimal-terms}
 V^i\simeq E(P^i)\qquad\text{for every }i.
\end{equation}
Indeed, for any simple $B$-module $S$, radicality makes the
differentials in $\operatorname{Hom}_B(P^\bullet,S)$ zero.  Hence
\[
 \operatorname{Hom}_{D^b(B)}(P^\bullet,S[j])
 \simeq \operatorname{Hom}_B(P^{-j},S).
\]
Applying $H$ gives the corresponding space
$\operatorname{Hom}_B(V^{-j},E(S))$, compatibly with the division
rings $\operatorname{End}_B(S)$ and $\operatorname{End}_B(E(S))$.
Their dimensions over the respective division rings determine the
multiplicities of the indecomposable projective covers, proving \eqref{anu-eq-minimal-terms}.  Existence of the projective
representative follows because $H$ preserves the thick subcategory
generated by the projectives.  An $R$-linear Morita equivalence also commutes with the Nakayama
functor on projectives, by the natural isomorphism
\[
 D\operatorname{Hom}_B(P,X)
 \simeq \operatorname{Hom}_B(X,\nu_B P)
\]
and the Yoneda lemma.  Hence the additive Nakayama-stability
conditions are preserved under \eqref{anu-eq-minimal-terms}.

Now put $T^\bullet=G(B)$.  Derived Morita theory supplies an $R$-linear
standard equivalence $G_0:D^b(B)\longrightarrow D^b(A)$ with
$G_0(B)\simeq T^\bullet$; the $R$-algebra identification of its
endomorphism algebra is the one induced by $G$.  See
\cite[Corollary~3.5]{anu-Rickard1991}, or
\cite[Section~2.2]{RouquierZimmermann2003} for the
two-sided realization over a projective commutative base.
The autoequivalence $H=FG_0$ satisfies $H(B)\simeq B$.
It is $t$-exact, since the groups
$\operatorname{Hom}_{D^b(B)}(B,X[j])$ detect the cohomology of $X$.
Set $F_0=G_0^{-1}$.  The first tilting complex of $F_0$ is
$T^\bullet$.  Its other tilting complex is
$F_0(A)\simeq H^{-1}F(A)$; by
\eqref{anu-eq-minimal-terms}, its radical terms are obtained from
those of $F(A)$ by the Morita equivalence $E^{-1}$ on the heart.
Thus both conditions in \eqref{anu-eq-almost-nu} hold for $F_0$.
\end{proof}

\begin{lemma}\label{anu-lem-nonselfinjective-base}
Let $R$ be local and not self-injective, and let $A$ be a nonzero
Artin $R$-algebra which is projective over $R$.  Then $A$ has no
nonzero finitely generated projective-injective left or right module.
Consequently, any $R$-linear almost $\nu$-stable derived equivalence
between two such algebras has its associated tilting complexes
concentrated in degree zero, and the algebras are $R$-linearly Morita
equivalent.
\end{lemma}

\begin{proof}
Restriction to $R$ takes an $A$-projective module to an $R$-projective
module.  It also takes an $A$-injective module to an $R$-injective
module: it is the right adjoint of the exact functor
$A\otimes_R-$.  A nonzero finite projective module over the local
ring $R$ is a nonzero free module.  If it were injective, its direct
summand $R$ would be injective, a contradiction.  The same argument
works on the right.
Each term in $T^-$ or $U^+$ from \eqref{anu-eq-almost-nu} is both
projective and injective, so these terms must vanish.  A tilting
complex concentrated in degree zero is a projective generator,
which gives the stated Morita equivalence.
\end{proof}

The next lemma replaces stabilized bimodule identities by actual
direct-sum decompositions after adjoining projective bimodules.

\begin{lemma}\label{anu-lem-padding}
Let $R$ be a commutative local Artin ring, and let $A,B$ be nonzero finite
free $R$-algebras.  Suppose ${}_A M_B$ and ${}_B N_A$ are finite and
projective on each side, and there are bimodule isomorphisms
\begin{equation}\label{anu-eq-stabilized}
 (M\otimes_BN)\oplus W_A\simeq A\oplus V_A,
 \qquad
 (N\otimes_AM)\oplus W_B\simeq B\oplus V_B,
\end{equation}
where $V_A,W_A$ are $A_R^e$-projective and $V_B,W_B$ are
$B_R^e$-projective.  Then, for a sufficiently large integer $t$, the bimodules
\[
 \widetilde M=M\oplus(A\otimes_RB)^t,
 \qquad
 \widetilde N=N\oplus(B\otimes_RA)^t
\]
satisfy the actual identities \eqref{anu-eq-identities} with
projective error bimodules.
\end{lemma}

\begin{proof}
Choose complements
$W_A\oplus W_A'\simeq(A_R^e)^u$ and
$W_B\oplus W_B'\simeq(B_R^e)^v$.  Write
$a=\operatorname{rank}_R A>0$ and
$b=\operatorname{rank}_R B>0$, and choose $t$ with
$t^2b\geq u$ and $t^2a\geq v$.
Expanding $\widetilde M\otimes_B\widetilde N$, the two mixed
summands are finite sums of $M\otimes_RA$ and $A\otimes_RN$.
They are $A_R^e$-projective because ${}_AM$ and $N_A$ are projective.
The remaining added summand is
\[
 (A\otimes_RB)^t\otimes_B(B\otimes_RA)^t
 \simeq(A\otimes_RB\otimes_RA)^{t^2}
 \simeq(A_R^e)^{t^2b}
\]
as an $A$-$A$-bimodule.  By \eqref{anu-eq-stabilized},
\[
 (M\otimes_BN)\oplus(A_R^e)^u
 \simeq A\oplus V_A\oplus W_A'.
\]
This yields the first required identity.  The reverse product has
the free added summand $(B_R^e)^{t^2a}$ and gives the second one.
The padded bimodules remain projective on both sides because
$A$ and $B$ are finite free over $R$.
\end{proof}

\subsection{Outer products over a self-injective base}

We first prove an intersection formula for additive subcategories.
It holds over every commutative Artin base.

\begin{lemma}\label{anu-lem-intersection}
Let $A,B$ be Artin $R$-algebras.  Let ${}_BV$ and $W_A$ be finite
projective modules.  Then, as additive subcategories of
$(B\otimes_RA^{\mathrm{op}})\text{-}\mathrm{mod}$,
\[
 \operatorname{add}(B\otimes_RW)
 \cap\operatorname{add}(V\otimes_RA)
 =\operatorname{add}(V\otimes_RW).
\]
\end{lemma}

\begin{proof}
Choose idempotents $e\in B$ and $f\in A$ with
$\operatorname{add}(V)=\operatorname{add}(Be)$ and
$\operatorname{add}(W)=\operatorname{add}(fA)$.  Put
$C=B\otimes_RA^{\mathrm{op}}$, $u=e\otimes1$, and
$v=1\otimes f^{\mathrm{op}}$.
For a simple $C$-module $S$, we claim that
\[
 uvS\ne0\quad\Longleftrightarrow\quad
 uS\ne0\ \hbox{and}\ vS\ne0.
\]
Only the reverse implication needs proof.  If $vS\ne0$ and $uvS=0$,
then $u$ kills $vS$.  The $R$-submodule $vS$ is stable under the $B$-action,
and $u$ commutes with the $A^{\mathrm{op}}$-action.  It follows that
$u$ kills $C(vS)=S$, contradicting $uS\ne0$.
An indecomposable projective cover of $S$ occurs in $Cu$ precisely
when $uS\ne0$, and similarly for $Cv$ and $Cuv$.
The Krull--Schmidt theorem therefore gives
$\operatorname{add}(Cu)\cap\operatorname{add}(Cv)
=\operatorname{add}(Cuv)$, which is the desired equality.
\end{proof}

\begin{lemma}\label{anu-lem-outer-injective}
Suppose that $R$ is a commutative self-injective local Artin ring,
and that $A,B$ are finite projective over $R$.  If ${}_BV$ and
$W_A$ are projective-injective, then $V\otimes_RW$ is
projective-injective over $B\otimes_RA^{\mathrm{op}}$.
\end{lemma}

\begin{proof}
Here $D_R=\operatorname{Hom}_R(-,R)$ is the duality for finite
$R$-modules.  Both $V$ and $W$ are finite projective over $R$.
Their outer product is projective over
$B\otimes_RA^{\mathrm{op}}$, and
\[
 D_R(V\otimes_RW)\simeq D_R(W)\otimes_RD_R(V)
\]
is projective over the opposite enveloping algebra: the duals of
the two injective factors are projective.  Applying the duality
proves injectivity of the original outer product.
\end{proof}

For an Artin algebra $A$, let $E_A$ be the direct sum of one
representative of each indecomposable projective module all of whose
positive Nakayama iterates are projective.  Thus
$\operatorname{add}(E_A)=\operatorname{add}(\nu_AE_A)$, and every
object of this subcategory is projective-injective.  This is the
maximal Nakayama-stable projective-injective subcategory used in
\cite[Proposition~3.8 and Lemma~3.9, pp.~129--130]{HuXi2010}.
Under \eqref{anu-eq-almost-nu}, one has $T^-\in\operatorname{add}(E_A)$,
$U^+\in\operatorname{add}(E_B)$, and
\begin{equation}\label{anu-eq-E-restriction}
 F(K^b(\operatorname{add}E_A))=K^b(\operatorname{add}E_B)
 \quad\hbox{up to isomorphism.}
\end{equation}
This equality and its inverse follow from Hu--Xi's Lemma~3.9,
which holds for Artin algebras.

\subsection{The two-sided construction}

\begin{proof}[Proof of Theorem~\ref{anu-thm-raw-semt}]
Decompose $R$ as a finite product of commutative Artinian local
rings.  Since $F$ is $R$-linear, it preserves the corresponding
central-idempotent decompositions of the derived categories.
The zero components of $A$ and $B$ correspond, and both the
hypotheses and the required bimodule identities hold componentwise.
It therefore suffices to work over a local $R$ with $A,B\ne0$.
In this case they are finite free over $R$.
Lemma~\ref{anu-lem-nonselfinjective-base} handles the case in which
$R$ is not self-injective, with $P=Q=0$.
For the remainder of the proof suppose that $R$ is self-injective,
and use $D_R=\operatorname{Hom}_R(-,R)$.
By Lemma~\ref{anu-lem-standard-replacement}, we may replace $F$ by
a standard almost $\nu$-stable equivalence.

\paragraph{A nonnegative two-sided representative.}
Let ${}_BX_A$ be a two-sided tilting complex inducing $F$, with
inverse ${}_AY_B$.  These complexes can be represented by bounded
complexes whose terms are projective on each side.  Indeed, take a bounded-above projective resolution over the
corresponding tensor product algebra.  Its restriction to each
side represents a perfect complex.  A sufficiently far negative
syzygy is therefore projective on both sides, and truncation at
that syzygy gives the required bounded representative.
Restriction preserves projectivity on each side because $A$ and
$B$ are projective over $R$.

We claim that $X$ has a radical representative
\begin{equation}\label{anu-eq-Delta}
 \Delta^\bullet=(0\longrightarrow N\longrightarrow\Delta^1
 \longrightarrow\cdots\longrightarrow\Delta^n\longrightarrow0),
 \qquad
 \Delta^i\in\operatorname{add}(U^+\otimes_R(T^-)^*)
 \quad(i>0),
\end{equation}
where ${}^*=\operatorname{Hom}_A(-,A)$ on left $A$-projectives.
To prove this assertion, retain the right $A$-action and consider
the derived equivalence
\[
 X\otimes_A^{\mathbf L}-:
 D^b(A_R^e)\longrightarrow D^b(B\otimes_RA^{\mathrm{op}}).
\]
Its two one-sided tilting complexes are represented by
$T^\bullet\otimes_RA$ and $U^\bullet\otimes_RA$.
Retaining the left $B$-action gives
\[
 -\otimes_B^{\mathbf L}X:
 D^b(B_R^e)\longrightarrow D^b(B\otimes_RA^{\mathrm{op}}),
\]
whose corresponding complexes are
$B\otimes_R\operatorname{Hom}_B(U^\bullet,B)$ and
$B\otimes_R\operatorname{Hom}_A(T^\bullet,A)$.
The chosen bounded models are projective, hence $K$-flat, on the
side over which each tensor product is taken.  Derived tensor associativity and the inverse identities give
inverse equivalences with these additional actions.  The bounded
models preserve boundedness and finite generation.  The displayed
one-sided representatives follow from projectivity over $R$.
Removing contractible summands, if necessary, preserves the
indicated degrees and additive supports.

Apply the normal-form statement
\cite[Lemma~3.1]{HuXi2010} to the diagonal modules $A$ and $B$
in these two equivalences.  It gives two nonnegative radical
representatives of the same object $X$, with positive terms,
respectively, in
\[
 \operatorname{add}(U^+\otimes_RA)
 \quad\hbox{and}\quad
 \operatorname{add}(B\otimes_R(T^-)^*).
\]
The comparison from homotopy morphisms to derived morphisms is an
isomorphism for nonnegative complexes whose positive terms are
projective \cite[Lemma~2.2(1), p.~114]{HuXi2010}.
Consequently the two representatives are homotopy equivalent,
and radicality makes them isomorphic as complexes.
Lemma~\ref{anu-lem-intersection} now proves
\eqref{anu-eq-Delta}.
In particular,
\begin{equation}\label{anu-eq-positive-support}
 \Delta^i\in\operatorname{add}(E_B\otimes_RE_A^*)
 \qquad(i>0).
\end{equation}

\paragraph{Projectivity of the degree-zero bimodules.}
After forgetting the right $A$-action, $\Delta^\bullet$ is
isomorphic to $U^\bullet$ in the derived category.  After forgetting
the left $B$-action, it is isomorphic to
$\operatorname{Hom}_A(T^\bullet,A)$.
In both cases these isomorphisms lift to homotopy equivalences by
the same comparison lemma.  The homotopy-inverse equation in
degree zero expresses $1_N$ as a map factoring through the direct
sum of the degree-zero projective term of the comparison complex
and $\Delta^1$.  Hence $N$ is projective both as a left $B$-module
and as a right $A$-module.

Put
\[
 \Theta^\bullet=\operatorname{Hom}_B^\bullet(\Delta^\bullet,B),
 \qquad M=\Theta^0=\operatorname{Hom}_B(N,B).
\]
This is an inverse two-sided tilting complex, supported in
$[-n,0]$.  Finite $R$-projectivity gives, for $i>0$,
\[
 \Theta^{-i}\in
 \operatorname{add}\bigl(D_R((T^-)^*)\otimes_R(U^+)^*\bigr)
 =\operatorname{add}\bigl(\nu_AT^-\otimes_R(U^+)^*\bigr)
 \subseteq\operatorname{add}(E_A\otimes_RE_B^*).
\]
The modules $E_A^*$ and $E_B^*$ are also projective-injective:
for instance $D_R(E_A^*)=\nu_AE_A$ is projective, and $E_A^*$
is projective by definition.  Lemma~\ref{anu-lem-outer-injective}
therefore shows that all nonzero-degree terms of both complexes
are projective-injective over their respective tensor product
algebras.

The module $M$ is projective as a right $B$-module.  To prove that
it is projective as a left $A$-module, compare $\Theta^\bullet$ with $T^\bullet$ after
forgetting the right $B$-action.  Both are nonpositive, and their
negative terms are injective left $A$-modules.  Thus their derived
isomorphism lifts to a homotopy equivalence by
\cite[Lemma~2.2(2)]{HuXi2010}.
The degree-zero homotopy-inverse equation now factors $1_M$
through $T^0\oplus\Theta^{-1}$, which is projective.
Thus $M$ is projective on both sides.

\paragraph{The four additive-support relations.}
We need
\begin{equation}\label{anu-eq-four-supports}
 \begin{aligned}
 N\otimes_AE_A&\in\operatorname{add}(E_B),&
 E_A^*\otimes_AM&\in\operatorname{add}(E_B^*),\\
 M\otimes_BE_B&\in\operatorname{add}(E_A),&
 E_B^*\otimes_BN&\in\operatorname{add}(E_A^*).
 \end{aligned}
\end{equation}
First, $E_A^*\otimes_AE_A\simeq\operatorname{End}_A(E_A)$
is finite projective over $R$: it is a direct summand of a matrix
algebra over $A$.  The analogous statement holds for $E_B$.
The same conclusion holds after replacing either factor by an
object of its additive closure.

The complex $\Delta^\bullet\otimes_AE_A$ is a bounded complex of
projective $B$-modules, and its positive terms lie in
$\operatorname{add}(E_B)$.  By \eqref{anu-eq-E-restriction}, it is
homotopy equivalent to a complex in
$K^b(\operatorname{add}E_B)$.
The degree-zero homotopy-inverse equation factors
$1_{N\otimes_AE_A}$ through objects in $\operatorname{add}(E_B)$;
this proves the first relation.  Applying the same argument to
the nonpositive complex $\Theta^\bullet\otimes_BE_B$ proves the
third relation.  Tensor--Hom adjunction and finite projectivity give
\[
 E_A^*\otimes_AM
 \simeq\operatorname{Hom}_B(N\otimes_AE_A,B),
\]
so the first relation implies the second.

For the fourth relation, use the derived bimodule isomorphism
\[
 \Delta^\bullet\simeq
 \operatorname{Hom}_A^\bullet(\Theta^\bullet,A).
\]
After tensoring with $E_B^*$ this gives, in $D^b(A^{\mathrm{op}})$,
\[
 E_B^*\otimes_B\Delta^\bullet
 \simeq
 \operatorname{Hom}_A^\bullet
       (\Theta^\bullet\otimes_BE_B,A).
\]
The right-hand side is homotopy equivalent to a complex in
$K^b(\operatorname{add}E_A^*)$, by
\eqref{anu-eq-E-restriction}.  The left-hand side is a bounded
complex of projective right $A$-modules, is nonnegative, and has
positive terms in $\operatorname{add}(E_A^*)$.
Its degree-zero homotopy-inverse equation therefore proves the
last relation in \eqref{anu-eq-four-supports}.
Only the derived duality of these complexes is used here; the
argument does not require an isomorphism
$N\simeq\operatorname{Hom}_A(M,A)$ of degree-zero bimodules.

\paragraph{Splitting the product complexes and padding.}
Consider $C^\bullet=\Delta^\bullet\otimes_A\Theta^\bullet$.
It is isomorphic to the diagonal bimodule $B$ in the derived category.
By \eqref{anu-eq-positive-support} and
\eqref{anu-eq-four-supports}, every tensor term
$\Delta^i\otimes_A\Theta^j$ other than the term with $(i,j)=(0,0)$
belongs to $\operatorname{add}(E_B\otimes_RE_B^*)$.
For $i>0,j<0$ this uses the $R$-projectivity of
$E_A^*\otimes_AE_A$; for $i=0$ or $j=0$ it uses the first row of
\eqref{anu-eq-four-supports}.
Such terms are $B_R^e$-projective-injective by
Lemma~\ref{anu-lem-outer-injective}.  Thus
\[
 C^0=(N\otimes_AM)\oplus W_B
\]
with $W_B$ projective-injective, and every term of $C^\bullet$
in a nonzero degree is projective-injective.

Since the only cohomology of $C^\bullet$ is $B$ in degree zero,
we may split off contractible two-term complexes from both ends.
At the leftmost negative degree the differential is a monomorphism
from an injective module, so it splits.  At the rightmost positive
degree the differential is an epimorphism onto a projective module,
so it splits.  Repeating these two operations gives
$C^0\simeq B\oplus V_B$ with $V_B$ projective-injective.
Consequently
\[
 (N\otimes_AM)\oplus W_B\simeq B\oplus V_B.
\]
The reverse product
$\Theta^\bullet\otimes_B\Delta^\bullet\simeq A$
is treated similarly, using the second row of
\eqref{anu-eq-four-supports}.  It gives
\[
 (M\otimes_BN)\oplus W_A\simeq A\oplus V_A
\]
with $W_A,V_A$ projective over $A_R^e$.
Lemma~\ref{anu-lem-padding} converts these stabilized
isomorphisms into \eqref{anu-eq-identities}.  Finally, the finite
product of the pairs constructed over the local factors of $R$
gives the required global pair.
\end{proof}

\begin{remark}\label{anu-rem-cancellation}
Over a general Artin base, the diagonal bimodule of a nonsemisimple
indecomposable algebra may be projective over its relative
enveloping algebra.  For example,
if $A=R$ is a nonsemisimple commutative local Artinian ring, then
$A_R^e=R$ and the diagonal bimodule is free.  Thus the nonprojectivity of the diagonal of a nonsemisimple block,
used in the field case, is unavailable over an Artin base.
Lemma~\ref{anu-lem-padding} replaces that cancellation argument
without excluding these blocks.
\end{remark}

\subsection{Relative stable equivalences and Gorenstein dimensions}

\begin{corollary}\label{anu-cor-relative}
Under the hypotheses of Theorem~\ref{anu-thm-raw-semt}, let
$\mathcal L_R(A)$ denote the full subcategory of finite $A$-modules
which are projective over $R$, and define $\mathcal L_R(B)$ similarly.
The constructed tensor functors induce mutually inverse equivalences
\[
 \mathcal L_R(A)/[A\text{-}\mathrm{proj}]
 \simeq
 \mathcal L_R(B)/[B\text{-}\mathrm{proj}].
\]
They also induce mutually inverse equivalences between the additive
quotients of $A\text{-}\mathrm{mod}$ and $B\text{-}\mathrm{mod}$ by
their respective relatively $R$-projective modules, namely the
additive closures of modules $A\otimes_RV$ and $B\otimes_RV$,
with $V$ a finite $R$-module.
\end{corollary}

\begin{proof}
Since $N_A$ is projective, $N\otimes_AX$ is a direct summand of
a finite direct sum of copies of $X$ as an $R$-module.  Thus it
preserves $R$-projectivity, and similarly for $M\otimes_B-$.
Both functors preserve ordinary projectives because the bimodules
are projective on the target side.  If $X$ is $R$-projective and $P$ is
$A_R^e$-projective, then $P\otimes_AX$ is a direct summand of a
finite direct sum of $A\otimes_RX$, hence is $A$-projective.
The two identities \eqref{anu-eq-identities} therefore give the
first pair of inverse equivalences.

For an arbitrary finite module $X$, the same argument shows that
$P\otimes_AX$ is relatively $R$-projective.
Moreover $N\otimes_A(A\otimes_RV)=N\otimes_RV$ is a direct
summand of a finite sum of $B\otimes_RV$ because ${}_BN$ is
projective.  Thus both functors preserve relative projectives, and
\eqref{anu-eq-identities} gives the second assertion.
\end{proof}

\begin{corollary}\label{anu-cor-gpd}
The pair in Theorem~\ref{anu-thm-raw-semt} can be chosen so that
its tensor functors preserve the Gorenstein projective dimension
of every finite module on both sides.  Explicitly, for finite
left modules $X$ over $A$ and $Y$ over $B$,
\[
 \operatorname{Gpd}_B(N\otimes_AX)=\operatorname{Gpd}_A X,
 \qquad
 \operatorname{Gpd}_A(M\otimes_BY)=\operatorname{Gpd}_B Y,
\]
including the value $\infty$.  For finite right modules $X'$ over
$A$ and $Y'$ over $B$,
\[
 \operatorname{Gpd}_{B^{\mathrm{op}}}(X'\otimes_AM)
 =\operatorname{Gpd}_{A^{\mathrm{op}}}X',
 \qquad
 \operatorname{Gpd}_{A^{\mathrm{op}}}(Y'\otimes_BN)
 =\operatorname{Gpd}_{B^{\mathrm{op}}}Y'.
\]
\end{corollary}

\begin{proof}
On every nonself-injective local factor of $R$, choose the Morita
pair from Lemma~\ref{anu-lem-nonselfinjective-base}; Morita
equivalences preserve totally acyclic complexes and their
dimensions.  On every self-injective local factor use the pair
constructed in the theorem, including the padding if necessary.
It satisfies all the hypotheses of
Theorem~\ref{thm:gpd-base-injective-semt}, which gives the displayed
equalities on that factor and for the opposite algebras.
Gorenstein projective dimension over a finite product of algebras
is the maximum of the dimensions of the components.  Combining
the finitely many factors proves the assertion.
\end{proof}

\begin{remark}\label{anu-rem-ordinary-boundary}
The tensor functors in Theorem~\ref{anu-thm-raw-semt} need not
realize the ordinary stable functor of the given almost
$\nu$-stable derived equivalence.  If $X$ is not projective over
$R$, tensoring a positive-degree term of $\Delta^\bullet$ with
$X$ may give a nonprojective $B$-module.  Similarly,
$P\otimes_AX$ need not be projective when $P$ is
$A_R^e$-projective.  Corollary~\ref{anu-cor-relative} gives the
two quotient categories on which the constructed tensor functors
are mutually inverse equivalences.
\end{remark}

\begin{example}\label{anu-ex-noetherian-nakayama-boundary}
The Artin duality used in this section cannot be replaced by
$D_R=\operatorname{Hom}_R(-,R)$ without changing the argument.
To see this, let $R=A=B=k[[t]]$, where $k$ is a field, and set
\[
 \nu_R=D_R\operatorname{Hom}_R(-,R).
\]
Then $\nu_R$ is naturally isomorphic to the identity on finite
projective $R$-modules.  For the derived autoequivalence $F=[-1]$
with inverse $G=[1]$, the complexes $T=G(R)=R[1]$ and
$U=F(R)=R[-1]$ are concentrated in degrees $-1$ and $1$,
respectively.  Thus $T^-=R=U^+$ and both formal equalities in
\eqref{anu-eq-almost-nu} hold with this $\nu_R$.
Nevertheless, $F$ is not $t$-exact, although $R$ is not
self-injective, since $\operatorname{Ext}_R^1(k,R)\simeq k$.
Hence the degree-zero conclusion of
Lemma~\ref{anu-lem-nonselfinjective-base} does not follow after
this replacement of the duality.  Matlis duality also does not
give a duality on finite $R$-modules: it sends $R$ to the
non-finitely generated module $E_R(k)$.

Theorem~\ref{thm:rz-relative-derived-stable} gives a positive
construction over Noetherian bases under explicit hypotheses on
the base duals.  Its relative Frobenius structure is the one
appropriate to that construction.  The example above concerns
only the proposed replacement of the duality in the Artinian
argument.
\end{example}

\appendix
\section{Abe--Hoshino algebras and projectivity over the base}
\label{sec:abe-boundary}
\subsection{The two Gorenstein conventions}

Abe and Hoshino use two notions of a Gorenstein algebra over a
commutative base.  In \cite[Definition~3.1]{AbeHoshino2007}, $R$ is
Gorenstein, $A$ has Gorenstein dimension zero as an $R$-module,
and $D_RA$ is a projective generator as a right $A$-module.
This definition does not require projectivity over $R$.  The
criterion in \cite[Theorem~4.5]{AbeHoshino2007} characterizes
preservation of the Gorenstein-dimension-zero and dual-generator
conditions by
$\operatorname{add}(T)=\operatorname{add}(\nu T)$ for a tilting
complex $T$; \cite[Corollary~4.6]{AbeHoshino2007} gives the
version in which $A_R$ and $B_R$ are projective.  The latter
criterion holds over a commutative Noetherian base without
requiring the base itself to be Gorenstein.

We use right-module conventions in this paragraph.  If $A_R$ is
projective and $D_RA$ is a projective
generator, then for a tilting complex
$T\in K^b(\mathrm{proj}\text{-}A)$ and
$B=\operatorname{End}_{K^b(\mathrm{proj}\text{-}A)}(T)$,
Corollary~4.6 states that
\begin{equation}
 \operatorname{add}(T)=\operatorname{add}(\nu T),
 \qquad \nu T=T\otimes_A^{\mathbf L}D_RA,
 \label{eq:abe-full-nu-invariance}
\end{equation}
is equivalent to the corresponding projectivity and dual-generator
conditions for $B$.  For these $R$-reflexive algebras the
dual-generator condition is equivalent on the two sides
\cite[Remark~3.2]{AbeHoshino2007}.  Together with
Theorem~\ref{thm:rz-relative-derived-stable}, this
criterion gives a relative stable equivalence of Morita type for
the associated standard derived equivalence, in the stabilized
sense of \eqref{eq:relative-stabilized-mt}.
Corollary~\ref{cor:rz-henselian-adjoint-type} applies when its
additional hypotheses hold.  Condition
\eqref{eq:abe-full-nu-invariance} concerns the whole tilting
complex, whereas almost $\nu$-stability concerns its specified
nonzero-degree terms with respect to the Artin duality.  These
are different hypotheses, and the Rouquier--Zimmermann
construction does not require almost $\nu$-stability.

The more general definition in
\cite[Definition~2.1]{AbeHoshino2009} involves a codimension $c$:
$R_{\mathfrak p}$ is
Gorenstein on the support of $A$,
\[
 \operatorname{Ext}_R^i(A,R)=0\quad(i\ne c),
 \qquad \Omega=\operatorname{Ext}_R^c(A,R)
       \text{ is a right }A\text{-projective generator}.
\]
Their Corollary~4.4 uses
$\nu=-\otimes_A^{\mathbf L}\Omega$ to characterize when
a derived equivalent algebra remains in the same codimension
class. Theorem~\ref{thm:abe-hoshino-answer} proves that concentration
of the base-ring Ext groups for the target algebra suffices, and
Corollary~\ref{cor:abe-hoshino-nakayama} gives the resulting
Nakayama stability. These conditions differ from finiteness of the left and
right self-injective dimensions.  In particular, over an Artin
Gorenstein base their
codimension-zero condition forces ordinary self-injectivity:
$D_RA$ is injective by restriction--coinduction, and $A$ is a
summand of a finite direct sum of copies of $D_RA$; duality
gives the assertion on the opposite side.

The condition $\mathrm{G\!\operatorname{-}dim}_R A=0$ alone
cannot replace base projectivity in
Proposition~\ref{prop:relative-frobenius-criterion} or in the
quoted Rouquier--Zimmermann theorem.  The following example shows
that the corresponding $R$-split
exact category need not be Frobenius.

\subsection{A self-dual algebra that is not projective over the base}

\begin{example}\label{ex:abe-nonprojective-relative}
Let $K$ be a field and put
\[
 R=K[t]/(t^2),\qquad
 A=K[x,y]/(x^2,y^2),\qquad t\longmapsto xy.
\]
Write $J=(x,y)=\operatorname{rad}A$,
$C=A/(xy)$ and $S=A/J$.  Then
\begin{equation}
 A_R\simeq R\oplus K^{\oplus2},\qquad
 \mathrm{G\!\operatorname{-}dim}_R A=0,
 \qquad D_RA\simeq A\quad\text{as bimodules}.
 \label{eq:abe-example-duality}
\end{equation}
Nevertheless, $A_R$ is not projective, and
\begin{equation}
 \mathcal P_R(A)=\operatorname{add}(A,C),\qquad
 \mathcal I_R(A)=\operatorname{add}(A,J),
 \qquad \mathcal P_R(A)\ne\mathcal I_R(A).
 \label{eq:abe-example-classes}
\end{equation}
The relative syzygy functor on
$A\text{-}\underline{\mathrm{mod}}_R$ is not an equivalence.
\end{example}

\begin{proof}
The subspace spanned by $1,xy$ is a free rank-one $R$-summand,
while the subspaces $Kx$ and $Ky$ are annihilated by $t$.
This gives the first decomposition in
\eqref{eq:abe-example-duality}, which is not free over the
local ring $R$.  Equivalently, multiplication by $t$ has
$K$-rank one on $A$, whereas it has rank two on a free
rank-two $R$-module.

The periodic complex with all terms $R$ and all differentials
multiplication by $t$ is exact and remains exact under
$\operatorname{Hom}_R(-,R)$.  The cokernel of each differential is $K$.
Thus $K$ has Gorenstein dimension zero over $R$, as do $R$ and
their finite direct sums, proving the dimension assertion.

The functional on $R$ taking the coefficient of $t$ identifies
$R$ with its $K$-dual.  Therefore
\[
 \operatorname{Hom}_R(A,R)\simeq\operatorname{Hom}_K(A,K)
\]
as $A$-bimodules.  The functional on $A$ taking the coefficient
of $xy$ is a symmetrizing form, so the latter bimodule is
isomorphic to $A$.  This proves the duality assertion.  In particular, this algebra satisfies the Gorenstein-dimension-zero
and dual-generator conditions in \cite[Definition~3.1]{AbeHoshino2007}.

Every finite $R$-module is a direct sum of copies of $R$ and
$K$.  Hence
\[
 \mathcal P_R(A)
  =\operatorname{add}(A\otimes_RR,A\otimes_RK)
  =\operatorname{add}(A,C).
\]
Similarly, relative injectives are generated by
$\operatorname{Hom}_R(A,R)\simeq A$ and
\[
 \operatorname{Hom}_R(A,K)
      \simeq D_K(A/(xy))\simeq\operatorname{ann}_A(xy)=J.
\]
The last isomorphism follows from the symmetrizing form on $A$.
The module $C$ is indecomposable because
$\operatorname{End}_A(C)\simeq C$ is local.  The module $J$
is indecomposable because it has a one-dimensional socle.
Both have dimension three, but
\[
 \dim_K\operatorname{top}C=1,
 \qquad \dim_K\operatorname{top}J=2.
\]
Thus they are not isomorphic, and neither is isomorphic to
the four-dimensional indecomposable module $A$.
The Krull--Schmidt theorem gives the inequality in
\eqref{eq:abe-example-classes}.  The $R$-split exact category
is therefore not Frobenius.

There is an $R$-split exact sequence
\begin{equation}
 0\longrightarrow S^{\oplus2}\longrightarrow C
       \xrightarrow{\pi} S\longrightarrow0.
 \label{eq:abe-example-relative-cover}
\end{equation}
Indeed, $t$ acts as zero on $C$, and
$\operatorname{rad}C=Kx\oplus Ky$ is annihilated by $J$.
The epimorphism $\pi$ is right minimal: an endomorphism of
$C$ inducing the identity on its simple top is a unit of the
local algebra $\operatorname{End}_A(C)$.  Thus $\pi$ is a relative projective cover, and we may take
$\Omega_R S=S^{\oplus2}$ in the relative stable category.

No nonzero endomorphism of $S$ factors through
$\operatorname{add}(A,C)$.  Otherwise a nonzero scalar multiple
of $1_S$ would so factor, making $S$ a direct summand of a
relative projective, contrary to Krull--Schmidt.  The same
argument applied to matrix entries gives
\[
 \underline{\operatorname{End}}_{A,R}(S)=K,
 \qquad
 \underline{\operatorname{End}}_{A,R}(S^{\oplus2})
       =\operatorname{Mat}_2(K).
\]
An equivalence would preserve endomorphism rings, but these
rings have $K$-dimensions one and four.  Hence $\Omega_R$
is not an equivalence.
\end{proof}

The example shows that the Gorenstein-dimension-zero and
dual-generator conditions alone do not make relative syzygy
invertible.  Relative stable equivalences nevertheless follow
from \eqref{eq:relative-raw-mt} or
\eqref{eq:relative-stabilized-mt}, whenever such bimodule data
are available, by Proposition~\ref{prop:relative-morita-bimodules}.

\subsection{A Brauer tree example}

\begin{example}\label{ex:relative-brauer-tree}
Let $K$ be algebraically closed.  Let $A_0$ and $B_0$ be the
basic Brauer tree algebras for the three-edge star and the
three-edge line, both with exceptional multiplicity one.
They are symmetric and standardly derived equivalent by the
classical Brauer-tree result of Rickard
\cite[Theorem~4.2]{RickardDerived1989}; see also the explicit
two-sided discussion in \cite[Sections~2.1--2.2]{Kozakai2019}.
Their Cartan matrices are
\[
 C_{A_0}=\begin{pmatrix}2&1&1\\1&2&1\\1&1&2\end{pmatrix},
 \qquad
 C_{B_0}=\begin{pmatrix}2&1&0\\1&2&1\\0&1&2\end{pmatrix}.
\]
The entries count the common vertices of the corresponding tree
edges, with diagonal entries equal to two.

Put $R=K[\epsilon]/(\epsilon^2)$,
$A=R\otimes_KA_0$ and $B=R\otimes_KB_0$.  Extending inverse
two-sided tilting complexes from $K$ to $R$ gives inverse
two-sided tilting complexes for $A$ and $B$, since extension
of scalars over a field is exact and preserves their tensor
identities.  Both algebras are finite free and symmetric over
$R$, so Theorem~\ref{thm:rz-relative-derived-stable} supplies
a relative Morita-type stable equivalence.
On the other hand, their Cartan matrices as finite-dimensional
$K$-algebras are $2C_{A_0}$ and $2C_{B_0}$.  These cannot be
identified by a permutation of the simple modules, because
one has zero off-diagonal entries and the other does not.
Thus $A$ and $B$ are not Morita equivalent.  The construction is
obtained by extending scalars in the classical Brauer tree example.

The same scalar-extension construction with $R=K[[t]]$ gives
finite free symmetric algebras over a complete Noetherian local
base of dimension one.  Their full $R$-relative stable
categories are triangle equivalent by
Theorem~\ref{thm:rz-relative-derived-stable}.  They are not
$R$-linearly Morita equivalent: such an equivalence would
restrict to a Morita equivalence between the subcategories
annihilated by $t$, and hence between $A_0$ and $B_0$, contrary
to their displayed Cartan matrices.
\end{example}

\section{Separable tensor functors and Gorenstein projective modules}
\label{sec:separable-boundary}

In a separable equivalence the regular bimodules are direct
summands of the tensor products of bimodules projective on each
side. We give an example over a field in which the complementary
summands fail the enveloping-algebra projectivity hypothesis
and the associated tensor functor fails to preserve Gorenstein
projective modules. This shows the role of the error bimodules
in Proposition~\ref{prop:gpd-error-functor-criterion}.

\begin{example}
\label{ex:separable-pair-not-gp-preserving}
Let $k$ be a field and $C=k[t]/(t^2)$.  Define the bimodule
${}_CL_C$ by taking $L_C=C_C$ and letting the left action factor
through the augmentation $C\to k$.  Form the triangular algebra
\[
 D=\begin{pmatrix}C&L\\0&C\end{pmatrix},
 \qquad
 u=\begin{pmatrix}t&0\\0&t\end{pmatrix}.
\]
The homomorphism $C\to D$, $c\mapsto\operatorname{diag}(c,c)$,
defines a bimodule $V={}_DD_C$ which is free of rank one on the
left and of rank three on the right.  Let $S=C/(t)$.  Since
$C$ is self-injective, $S$ is Gorenstein projective.  Nevertheless,
\[
 Z=V\otimes_CS\simeq D/Du
\]
is not Gorenstein projective.

Indeed, right multiplication by $u$ gives a projective resolution
\[
 \cdots\xrightarrow{\cdot u}D\xrightarrow{\cdot u}D
 \xrightarrow{\cdot u}D\longrightarrow Z\longrightarrow0.
\]
It is exact because $D_C\simeq C^3$.  Applying
$\operatorname{Hom}_D(-,D)$ changes the differential into left
multiplication by $u$.  Write an element of $D$ as $(a,\ell,b)$.
Then
\[
 u(a,\ell,b)=(ta,0,tb),
\]
and hence
\[
 \ker(u\cdot-)=kt\oplus L\oplus kt,
 \qquad
 \operatorname{im}(u\cdot-)=kt\oplus0\oplus kt.
\]
Thus $\operatorname{Ext}_D^i(Z,D)\simeq L\neq0$ as $k$-vector
spaces for every $i\geq1$.  In particular
$\operatorname{Gpd}_D Z=\infty$.

Now put $\Lambda=C\times D$ and let $U$ be the $\Lambda$-bimodule
whose only nonzero component is
$e_DUe_C=V$.  It is finitely generated projective on both sides.
The bimodules
\[
 {}_\Lambda M_\Lambda={}_\Lambda\Lambda_\Lambda,
 \qquad
 {}_\Lambda N_\Lambda=\Lambda\oplus U
\]
satisfy
\[
 M\otimes_\Lambda N\simeq\Lambda\oplus U,
 \qquad
 N\otimes_\Lambda M\simeq\Lambda\oplus U,
\]
and hence define a separable equivalence.  However, the image of
the Gorenstein projective module $(S,0)$ is
\[
 N\otimes_\Lambda(S,0)\simeq(S,Z),
\]
whose Gorenstein projective dimension is infinite.

The error bimodule $U$ is projective on both sides but is not
projective over $\Lambda^e=\Lambda\otimes_k\Lambda^{\mathrm{op}}$.
Indeed, if $U$ were projective over $\Lambda^e$, then
$U\otimes_\Lambda(S,0)$ would be a projective left $\Lambda$-module,
contrary to the calculation above.  Thus this pair does not define
a stable equivalence of Morita type.
\end{example}

\begin{remark}
Theorem~A (Theorem~3.2) of Sun--Zhao \cite{SunZhaoSeparable2025}
asserts that a specified separable tensor pair preserves Gorenstein
projective modules. Under the bimodule definition above,
Example~\ref{ex:separable-pair-not-gp-preserving} shows that this
assertion does not hold in general. The example concerns preservation
by the prescribed tensor functors. It does not rule out a statement
with a Frobenius or biadjoint hypothesis, or an invariance statement
formulated in terms of the existence of an equivalence between the
algebras.
\end{remark}

\begingroup
\small\raggedright
\makeatletter
\renewcommand{\@biblabel}[1]{[#1]}
\makeatother
\providecommand{\bysame}{\leavevmode\hbox to3em{\hrulefill}\thinspace}
\providecommand{\MR}{\relax\ifhmode\unskip\space\fi MR }
% \MRhref is called by the amsart/book/proc definition of \MR.
\providecommand{\MRhref}[2]{%
  \href{http://www.ams.org/mathscinet-getitem?mr=#1}{#2}
}
\providecommand{\href}[2]{#2}

\endgroup
\end{document}